\documentclass[a4paper,12pt,reqno]{amsart}
\usepackage{amsmath}
\usepackage{amssymb}
\usepackage{amsthm}

\usepackage{amscd}
\usepackage{amsfonts}
\usepackage{setspace}
\usepackage{version}

\title[Low moments of random multiplicative functions]{Moments of random multiplicative functions, I: Low moments, better than squareroot cancellation, and critical multiplicative chaos}
\author{Adam J Harper}
\address{Mathematics Institute, Zeeman Building, University of Warwick, Coventry CV4 7AL, England}
\email{A.Harper@warwick.ac.uk}
\date{20th March 2017}
\thanks{When this work was started, the author was supported by a research fellowship at Jesus College, Cambridge. The work was completed while the author was in residence at the Mathematical Sciences Research Institute in Berkeley, California (supported by the National Science Foundation under Grant No. DMS-1440140), during the Spring 2017 semester.}

\numberwithin{equation}{section}
\theoremstyle{plain}

\onehalfspace
\newcommand{\N}{\mathbb{N}}
\newcommand{\R}{\mathbb{R}}
\newcommand{\E}{\mathbb{E}}
\newcommand{\p}{\mathbb{P}}
\newcommand{\Z}{\mathbb{Z}}

\newtheorem{thmlower1}{Theorem}
\newtheorem{thmlower2}[thmlower1]{Theorem}
\newtheorem{cor1}{Corollary}
\newtheorem{cor2}[cor1]{Corollary}

\newtheorem{numth1}{Number Theory Result}
\newtheorem{harman1}{Harmonic Analysis Result}
\newtheorem{prop1}{Proposition}
\newtheorem{prop2}[prop1]{Proposition}
\newtheorem{prop3}[prop1]{Proposition}
\newtheorem{prop4}[prop1]{Proposition}

\newtheorem{lem1}{Lemma}
\newtheorem{lem2}[lem1]{Lemma}
\newtheorem{lem3}[lem1]{Lemma}
\newtheorem{lem4}[lem1]{Lemma}
\newtheorem{probres1}{Probability Result}
\newtheorem{probres2}[probres1]{Probability Result}
\newtheorem{prop5}[prop1]{Proposition}
\newtheorem{lem5}[lem1]{Lemma}
\newtheorem{prop6}[prop1]{Proposition}

\newtheorem{keyprop1}{Key Proposition}
\newtheorem{keyprop2}[keyprop1]{Key Proposition}
\newtheorem{keyprop3}[keyprop1]{Key Proposition}
\newtheorem{keyprop4}[keyprop1]{Key Proposition}

\newtheorem{lem6}[lem1]{Lemma}
\newtheorem{lem7}[lem1]{Lemma}
\newtheorem{prop7}[prop1]{Proposition}
\newtheorem{keyprop5}[keyprop1]{Key Proposition}

\newtheorem{problem1}{Probability Lemma}

\begin{document}

\maketitle

\begin{abstract}
We determine the order of magnitude of $\E|\sum_{n \leq x} f(n)|^{2q}$, where $f(n)$ is a Steinhaus or Rademacher random multiplicative function, and $0 \leq q \leq 1$. In the Steinhaus case, this is equivalent to determining the order of $\lim_{T \rightarrow \infty} \frac{1}{T} \int_{0}^{T} |\sum_{n \leq x} n^{-it}|^{2q} dt$.

In particular, we find that $\E|\sum_{n \leq x} f(n)| \asymp \sqrt{x}/(\log\log x)^{1/4}$. This proves a conjecture of Helson that one should have better than squareroot cancellation in the first moment, and disproves counter-conjectures of various other authors. We deduce some consequences for the distribution and large deviations of $\sum_{n \leq x} f(n)$.

The proofs develop a connection between $\E|\sum_{n \leq x} f(n)|^{2q}$ and the $q$-th moment of a critical, approximately Gaussian, multiplicative chaos, and then establish the required estimates for that. We include some general introductory discussion about critical multiplicative chaos to help readers unfamiliar with that area.
\end{abstract}

\section{Introduction}
Let $(f(p))_{p \; \text{prime}}$ be a sequence of independent Steinhaus random variables, i.e. independent random variables distributed uniformly on the unit circle $\{|z|=1\}$. We define a {\em Steinhaus random multiplicative function} $f$, by setting $f(n) := \prod_{p^{a} || n} f(p)^{a}$ for all natural numbers $n$ (where $p^a || n$ means that $p^a$ is the highest power of the prime $p$ that divides $n$, so $n = \prod_{p^{a} || n} p^{a}$). Thus $f$ is a random function taking values in the complex unit circle, that is totally multiplicative. An alternative model is to let $(f(p))_{p \; \text{prime}}$ be independent Rademacher random variables, taking values $\pm 1$ with probability $1/2$ each. Then we define a {\em Rademacher random multiplicative function} $f$, supported on squarefree numbers $n$ (i.e. numbers $n$ not divisible by any squares of primes) only, by $f(n) := \prod_{p |n} f(p)$.

Random multiplicative functions have attracted quite a lot of attention as models for functions of number theoretic interest: for example, Rademacher random multiplicative functions were introduced by Wintner~\cite{wintner} as a model for the M\"{o}bius function $\mu(n)$. In the Steinhaus case, for any real $q \geq 0$ and any given $x$ we have
$$ \E|\sum_{n \leq x} f(n)|^{2q} = \lim_{T \rightarrow \infty} \frac{1}{T} \int_{0}^{T} |\sum_{n \leq x} n^{-it}|^{2q} dt = \lim_{\text{prime} \; p \rightarrow \infty} \frac{1}{p-1} \sum_{\chi \; \text{mod} \; p} |\sum_{n \leq x} \chi(n)|^{2q} , $$
where the final sum is over the $p-1$ Dirichlet characters $\chi$ mod $p$. Thus, questions about moments of a Steinhaus random multiplicative function are equivalent to questions about the limiting behaviour of so-called zeta sums, or of character sums.

Since we always have $\E f(p) = 0$, and since numbers have unique prime factorisations, it is easy to check that the random variables $f(n)$ are {\em orthogonal} (i.e. $\E f(n) \overline{f(m)} = \textbf{1}_{n=m} |f(n)|^2$ for all $n,m$, where $\textbf{1}$ denotes the indicator function), and so
\begin{equation}
\E|\sum_{n \leq x} f(n)|^2 = \left\{ \begin{array}{ll}
     \sum_{n \leq x} 1 = x + O(1) & \text{in the Steinhaus case} , \\
     \sum_{n \leq x, \; n \; \text{squarefree}} 1 = \frac{6}{\pi^2}x + O(\sqrt{x})  & \text{in the Rademacher case} .
\end{array} \right. \nonumber
\end{equation}
With this orthogonality in mind, as well as the fact that the $f(n)$ are built from ``quite a lot'' of independent random variables $(f(p))_{p \; \text{prime}}$, there are various natural questions one might investigate that compare the behaviour of $f(n)$ with a sequence of completely independent random variables:
\begin{enumerate}
\item does one have a central limit theorem for $\frac{\sum_{n \leq x} f(n)}{\sqrt{\E|\sum_{n \leq x} f(n)|^2}}$?

\item what is the size of the moments $\E|\sum_{n \leq x} f(n)|^{2q}$, for general $q \geq 0$?

\item how large are the tail probabilities $\p(|\sum_{n \leq x} f(n)| \geq \lambda\sqrt{x})$, for $\lambda$ large? In particular, do they decay rapidly (e.g. exponentially) with $\lambda$?

\item does one have ``law of the iterated logarithm'' type almost sure bounds for $\sum_{n \leq x} f(n)$?
\end{enumerate}

Since the values $f(n)$ have a rather intricate dependence structure, these problems are probabilistically interesting. There are also number theoretic and analytic motivations for them. For example, Helson~\cite{helson} conjectured in the Steinhaus case one should have $\E|\sum_{n \leq x} f(n)| = o(\sqrt{x})$, and observed that if this is true then a certain generalisation of Nehari's theorem from harmonic analysis is false. See Saksman and Seip's open problems paper~\cite{saksmanseip} for a functional analysis perspective on Helson's conjecture and related questions. On the number theoretic side, it is well known that the Riemann Hypothesis is true if and only if $\sum_{n \leq x} \mu(n) = O_{\epsilon}(x^{1/2+\epsilon})$ for all positive $\epsilon$. In fact, a conjecture of Gonek based on the distribution of zeros of the Riemann zeta function (see Ng's paper~\cite{ng}) asserts that this sum should be $O(\sqrt{x} (\log\log\log x)^{5/4})$, and infinitely often as large as that. It would be very interesting to have sharp almost sure bounds for $\sum_{n \leq x} f(n)$ in the Rademacher case, to compare with Gonek's conjecture.

\vspace{12pt}
In this paper we will answer the second question for $0 \leq q \leq 1$, and derive some consequences for the other questions. In a companion paper~\cite{harperrmfhigh}, we will also answer the second question for $q > 1$, up to factors of size $e^{O(q^2)}$. Before explaining our theorems, we briefly summarise some of the previous literature on these problems. 

Regarding the moments $\E|\sum_{n \leq x} f(n)|^{2q}$, when $q$ is a fixed natural number one can expand the $2q$-th power and reduce the expectation calculation to a number theoretic counting problem. This leads to an asymptotic for the moment as $x \rightarrow \infty$, in both the Steinhaus and Rademacher cases: see the papers of Harper, Nikeghbali and Radziwi{\l}{\l}~\cite{hnr} and Heap and Lindqvist~\cite{heaplindqvist} for such calculations, and further references. We simply note here that the moments are of the form $x^q (\log x)^{\Theta(q^2)}$, so grow rapidly with $q$ in a way very unlike a sum of independent random variables. When $q$ is not a natural number there are no obvious methods available, so much less is known, but Helson~\cite{helson} conjectured (in the Steinhaus case) that the first moment should satisfy $\E|\sum_{n \leq x} f(n)| = o(\sqrt{x})$ as $x \rightarrow \infty$. This seems surprising from a number theoretic perspective, since one rarely expects to achieve better than squareroot cancellation. Exploring the conjecture, Weber~\cite{weber} established various results; Bondarenko and Seip~\cite{BondarenkoSeip} proved a lower bound $\E|\sum_{n \leq x} f(n)| \gg \sqrt{x}/(\log x)^{\delta}$ for a certain small explicit $\delta > 0$; and Harper, Nikeghbali and Radziwi{\l}{\l}~\cite{hnr} proved a stronger lower bound $\E|\sum_{n \leq x} f(n)| \gg \sqrt{x}/(\log\log x)^{3+o(1)}$. They also conjectured, in opposition to Helson's conjecture, that $\E|\sum_{n \leq x} f(n)|^{2q} \sim C(q) x^{q}$ as $x \rightarrow \infty$, for each fixed $0 \leq q \leq 1$.

Turning to the almost sure behaviour of $\sum_{n \leq x} f(n)$, it is known in the Rademacher case that the sum is almost surely $O(\sqrt{x} (\log\log x)^{2+\epsilon})$ for each $\epsilon > 0$, due to work of Lau, Tenenbaum and Wu~\cite{tenenbaum}. It is also known that the sum is almost surely {\em not} $O(\sqrt{x}/(\log\log x)^{5/2+\epsilon})$, due to work of Harper~\cite{harpergp}. This builds on previous work of many people, most notably Hal\'{a}sz~\cite{halasz}. These results are described by Lau, Tenenbaum and Wu as ``qualitatively matching the law of iterated logarithm'', so one might expect the proofs to involve Gaussian-type decay estimates for the tail probabilities $\p(|\sum_{n \leq x} f(n)| \geq \lambda\sqrt{x})$. However, so far as the author is aware, no bound is known for these probabilities for moderately sized $\lambda$ that improves on the Chebychev upper bound $\ll 1/\lambda^2$. The Lau--Tenenbaum--Wu almost sure upper bound instead exploits the special fact that $(1/\sqrt{x}) \sum_{n \leq x} f(n)$ changes size extremely slowly.

Finally looking at distributional questions, in the Rademacher case it is a natural default conjecture that $\frac{\sum_{n \leq x} f(n)}{\sqrt{\E|\sum_{n \leq x} f(n)|^2}} \stackrel{d}{\rightarrow} N(0,1)$ as $x \rightarrow \infty$, but Chatterjee (see $\S 6$ of Hough~\cite{hough}) conjectured that this should not hold. Chatterjee's conjecture was proved by the author~\cite{harperlimits}, using a special conditioning argument. On the other hand, if one restricts to many natural subsums one does have a central limit theorem: see the papers of Chatterjee and Soundararajan~\cite{chatsound}, Harper~\cite{harperlimits} and Hough~\cite{hough} for examples of such theorems, proved using Stein's method, a martingale decomposition, and the method of moments, respectively. It has remained an open question whether $\frac{\sum_{n \leq x} f(n)}{\sqrt{\E|\sum_{n \leq x} f(n)|^2}}$ has a limit distribution, and if so what.

\subsection{Statement of results}
We shall prove the following theorems, that determine the order of magnitude of $\E|\sum_{n \leq x} f(n)|^{2q}$ for all $0 \leq q \leq 1$, in both the Steinhaus and Rademacher cases.

\begin{thmlower1}
If $f(n)$ is a Steinhaus random multiplicative function, then uniformly for all large $x$ and $0 \leq q \leq 1$ we have
$$ \E|\sum_{n \leq x} f(n)|^{2q} \asymp \left( \frac{x}{1 + (1-q)\sqrt{\log\log x}} \right)^q . $$
\end{thmlower1}

\begin{thmlower2}
If $f(n)$ is a Rademacher random multiplicative function, then uniformly for all large $x$ and $0 \leq q \leq 1$ we have
$$ \E|\sum_{n \leq x} f(n)|^{2q} \asymp \left( \frac{x}{1 + (1-q)\sqrt{\log\log x}} \right)^q . $$ 
\end{thmlower2}

In particular, we find that $\E|\sum_{n \leq x} f(n)| \asymp \frac{\sqrt{x}}{(\log\log x)^{1/4}}$, which proves Helson's~\cite{helson} somewhat surprising conjecture that the first moment should be $o(\sqrt{x})$, and disproves the counter-conjecture of Harper, Nikeghbali and Radziwi{\l}{\l}~\cite{hnr} (see also Conjecture 1 of Heap and Lindqvist~\cite{heaplindqvist}).

Theorem 1 also implies a negative answer to the so-called embedding problem for Dirichlet polynomials (see Question 2 of \cite{saksmanseipintmeans}, or Problem 2.1 of \cite{saksmanseip}) for all exponents $0 < 2q < 2$. For it is easy to check using Riemann--Stieltjes integration that
$$ \int_{0}^{1} |\sum_{n \leq x} \frac{1}{n^{1/2+it}}|^{2q} dt = \int_{0}^{1} |\frac{x^{1/2-it}}{1/2-it} + O(1)|^{2q} dt \asymp x^q , $$
and Theorem 1 implies that $\lim_{T \rightarrow \infty} \frac{1}{2T} \int_{-T}^{T} |\sum_{n \leq x} n^{-it}|^{2q} dt = \E|\sum_{n \leq x} f(n)|^{2q} = o(x^q)$ for any fixed $0 < q < 1$, as $x \rightarrow \infty$. Thus there cannot exist any universal constant $C_{2q}$ such that $\int_{0}^{1} |P(1/2+it)|^{2q} dt \leq C_{2q} \lim_{T \rightarrow \infty} \frac{1}{2T} \int_{-T}^{T} |P(it)|^{2q} dt$ for all Dirichlet polynomials $P(s)$, as the Dirichlet polynomials $P(s) = \sum_{n \leq x} \frac{1}{n^s}$ for growing $x$ provide a sequence of counterexamples.

\vspace{12pt}
The proofs of the theorems divide into two parts, which we try to explain now. There are technical differences between the Steinhaus and Rademacher cases, but for $0 \leq q \leq 1$ most of the behaviour and key proof ideas are identical, so we confine this introductory discussion to the Steinhaus case.

The key number theoretic tool for studying multiplicative functions is complex analysis of the corresponding Euler products and Dirichlet series, $F(s) = \prod_{p \leq x} (1 - \frac{f(p)}{p^{s}})^{-1} = \sum_{\substack{n=1, \\ p|n \Rightarrow p \leq x}}^{\infty} \frac{f(n)}{n^s}$. This tool is particularly appealing for random multiplicative $f(n)$, because in the Euler product the different factors $(1-\frac{f(p)}{p^s})^{-1}$ are independent, whereas in the sum $\sum_{n \leq x} f(n)$ the contributions from the underlying independent $f(p)$ are entangled with one another in a highly non-trivial way. Thus the work of Hal\'{a}sz~\cite{halasz} and Harper~\cite{harpergp} on almost sure lower bounds for $\sum_{n \leq x} f(n)$ relied on a connection with lower bounds for the Euler product, and the work of Hal\'{a}sz and of Lau, Tenenbaum and Wu~\cite{tenenbaum} on upper bounds can also partially be understood in that way, though it isn't presented like that. But passing directly from $\sum_{n \leq x} f(n)$ to $F(s)$, for example using Perron's formula, wastes logarithmic factors that would be fatal when trying to prove our theorems. Thus the first stage of our proofs is to pass from $\E|\sum_{n \leq x} f(n)|^{2q}$ to a corresponding expectation involving $F(s)$, in an efficient way. Having done this, the second stage is to analyse the expectation involving $F(s)$.

More precisely, we first show (very roughly speaking) that
\begin{equation}\label{introprod}
\E|\sum_{n \leq x} f(n)|^{2q} \approx x^{q} \E\Biggl(\frac{1}{\log x} \int_{-1/2}^{1/2} |F(1/2+it)|^2 dt \Biggr)^{q} .
\end{equation}
This is done in a few steps. We show that it suffices to prove a comparable statement for quantities like $\E|\sum_{n \leq x, P(n) > \sqrt{x}} f(n)|^{2q}$, where $P(n)$ denotes the largest prime factor of $n$. (For the lower bound this is literally true, since if $\sum_{n \leq x, P(n) > \sqrt{x}} f(n)$ is large then with positive conditional probability the complete sum will also be large. For the upper bound, one splits the sum into several pieces depending on the size of $P(n)$, on various ranges, and gives a separate but similarly-shaped upper bound for the expectation of each piece.) The advantage of these sums is that, by multiplicativity, we can write
$$ \sum_{n \leq x, P(n) > \sqrt{x}} f(n) = \sum_{\sqrt{x} < p \leq x \; \text{prime}} f(p) \sum_{m \leq x/p} f(m) , $$
and importantly since $x/p \leq \sqrt{x}$ the inner sums are independent of the outer random variables $(f(p))_{\sqrt{x} < p \leq x}$. Now we might expect a sum of independent $f(p)$ weighted by the ``coefficients'' $\sum_{m \leq x/p} f(m)$ to behave, very roughly speaking, like a Gaussian with mean 0 and variance $\sum_{\sqrt{x} < p \leq x \; \text{prime}} |\sum_{m \leq x/p} f(m)|^2$. If this were the case, we would have
$$ \E\Biggl|\sum_{n \leq x, P(n) > \sqrt{x}} f(n) \Biggr|^{2q} \approx \E \Biggl( \sum_{\sqrt{x} < p \leq x \; \text{prime}} \Biggl|\sum_{m \leq x/p} f(m) \Biggr|^2 \Biggr)^{q} , $$
and this statement follows rigorously from Khintchine's inequality or, for the upper bound on our range of $q$, a suitable application of H\"{o}lder's inequality. Now we have passed from examining $\E|\sum_{n \leq x} f(n)|^{2q}$ at a single point $x$, to examining the average behaviour of sums of $f(n)$ up to various points $x/p$. Since primes $p$ are quite well distributed, one can work to replace the sum over primes by a corresponding integral, and after changing variables reduce to studying
$$ \frac{x^{q}}{\log^{q}x} \E\left( \int_{1}^{\sqrt{x}} \left|\sum_{m \leq z} f(m) \right|^2 \frac{dz}{z^{2}} \right)^q . $$
Finally, a version of Parseval's identity applied to this integral brings us to \eqref{introprod}.

As noted above, the foregoing argument has various antecedents that we should mention. In the work of Harper, Nikeghbali and Radziwi{\l}{\l}~\cite{hnr} on lower bounds for $\E|\sum_{n \leq x} f(n)|$, one follows the same steps as far as the application of Khintchine's inequality (although only for the lower bound and for $q=1/2$), but then establishes a connection with $F(s)$ in a different and less efficient way. In the work of Hal\'{a}sz~\cite{halasz} and of Lau, Tenenbaum and Wu~\cite{tenenbaum} on almost sure upper bounds, one also pulls out the value of $f$ on large primes and ends up dealing with integral averages of $|\sum_{n \leq x} f(n)|^2$, although these arise in a different way because from the beginning of their problem one is averaging over a sequence of $x$ values. Those authors never apply Parseval's identity, which could easily be done for the integral expressions they arrive at, but wouldn't in itself improve their results.

Now we try to explain our approach to analysing $x^{q} \E\left(\frac{1}{\log x} \int_{-1/2}^{1/2} |F(1/2+it)|^2 dt \right)^{q}$. To set the scene we note that, by H\"{o}lder's inequality,
$$ x^{q} \E\Biggl(\frac{1}{\log x} \int_{-1/2}^{1/2} |F(1/2+it)|^2 dt \Biggr)^{q} \leq x^{q} \Biggl(\frac{1}{\log x} \int_{-1/2}^{1/2} \E |F(1/2+it)|^2 dt \Biggr)^{q} \;\;\; \forall \; 0 \leq q \leq 1 . $$
A standard calculation shows $\E |F(1/2+it)|^2 \asymp \log x$ (see \eqref{stprodusual}, below), and inserting this yields the trivial upper bound $\E|\sum_{n \leq x} f(n)|^{2q} \ll x^{q}$. The major contribution to this expected size of $|F(1/2+it)|^2$ comes from the fairly rare event that $|\log|F(1/2+it)| - \log\log x| \leq \sqrt{\log\log x}$, but if integrating over $[-1/2,1/2]$ roughly corresponded to taking $\log x$ independent samples of $|F(1/2+it)|$ (because $F(s)$ varies with $s$ on a scale of $1/\log x$), one might indeed typically find a few such values of $\log|F(1/2+it)|$ with $|t| \leq 1/2$. (This discussion of $\E |F(1/2+it)|^2$ is similar to the analysis of the second moment of the Riemann zeta function, and the values of zeta that make the major contribution to it. See the introduction to Soundararajan's paper~\cite{soundmoments}, for example.) So the essence of Theorems 1 and 2 is that, when looking at $\E\left(\frac{1}{\log x} \int_{-1/2}^{1/2} |F(1/2+it)|^2 dt \right)^{q}$ with $q$ a little smaller than 1, integrating over $[-1/2,1/2]$ {\em does not correspond} to taking $\log x$ independent samples of $|F(1/2+it)|$, so the above application of H\"{o}lder's inequality is wasteful.

It turns out that $\frac{1}{\log x} \int_{-1/2}^{1/2} |F(1/2+it)|^2 dt$ is fairly close to (the total mass of a truncation of) a probabilistic object called {\em critical multiplicative chaos}, and our analysis of it draws on ideas from that field. We shall sketch these ideas now without assuming familiarity with the area, and then remark on connections with the literature.

There is some dependence between the random products $(|F(1/2+it)|^2)_{|t| \leq 1/2}$, since if $t$ changes slightly then $p^{it} = e^{it\log p}$ only changes slightly. The largest primes involved in $F(1/2+it)$ have size $p \approx x$, which is why we expect the product to change little when $t$ varies by less than $1/\log x$. But many primes in the product are much smaller than $x$, so at least the subproduct over smaller primes will remain unchanged over wider $t$ intervals. In fact, to understand $(|F(1/2+it)|^2)_{|t| \leq 1/2}$ properly one should think of the product consisting of $\log\log x$ blocks or ``scales'' of primes of comparable logarithmic size, each of which remains constant on a different $t$ scale.

In view of these non-trivial dependencies amongst the $(|F(1/2+it)|^2)_{|t| \leq 1/2}$, there are certain events involving the size of different subproducts of $|F(1/2+it)|^2$ that occur with probability close to 1, but would not do so if the products behaved independently. Now if $\mathcal{G}$ is some event, and if we let $q'=(1+q)/2$, we have that $\E\left(\frac{1}{\log x} \int_{-1/2}^{1/2} |F(1/2+it)|^2 dt \right)^{q}$ is
\begin{eqnarray}
& = & \E\Biggl(\textbf{1}_{\mathcal{G}} \frac{1}{\log x} \int_{-1/2}^{1/2} |F(1/2+it)|^2 \Biggr)^{q} + \E\Biggl(\textbf{1}_{\mathcal{G} \; \text{fails}} \frac{1}{\log x} \int_{-1/2}^{1/2} |F(1/2+it)|^2 \Biggr)^{q} \nonumber \\
& \leq & \E\Biggl(\textbf{1}_{\mathcal{G}} \frac{1}{\log x} \int_{-1/2}^{1/2} |F(\frac{1}{2}+it)|^2 \Biggr)^{q} + \p(\mathcal{G} \; \text{fails})^{\frac{q'-q}{q'}} \Biggl( \E\Biggl( \frac{1}{\log x} \int_{-1/2}^{1/2} |F(\frac{1}{2}+it)|^2 \Biggr)^{q'} \Biggr)^{\frac{q}{q'}} , \nonumber
\end{eqnarray}
the second line following by H\"{o}lder's inequality. We may now apply H\"{o}lder's inequality to the first term, as above, but hope to obtain additional savings because of the indicator function $\textbf{1}_{\mathcal{G}}$. In the second term, the prefactor $\p(\mathcal{G} \; \text{fails})^{\frac{q'-q}{q'}} = \p(\mathcal{G} \; \text{fails})^{\frac{1-q}{2q'}}$ can supply a saving if the event $\mathcal{G}$ is sufficiently probable, in particular if it occurs with probability larger than $1 - e^{-K/(1-q)}$, for some large $K$.

It turns out that if $\mathcal{G}$ is the event, not just that $\log|F(1/2+it)| \leq \log\log x + K/(1-q)$, but that all subproducts of $F(1/2+it)$ obey a comparable bound (with $\log\log x$ replaced by the number of ``scales'' involved in the subproduct), then the indicator function $\textbf{1}_{\mathcal{G}}$ produces a saving factor of the shape $\frac{K}{(1-q)\sqrt{\log\log x}}$. Meanwhile, $\mathcal{G}$ does occur with probability larger than $1 - e^{-K/(1-q)}$. This is essentially the argument that leads to the upper bounds in Theorems 1 and 2. To make things work, one needs to develop results that allow the estimation of $\E \textbf{1}_{\mathcal{G}} |F(1/2+it)|^2$ (an analogue of Girsanov's theorem from the theory of Gaussian random variables). Some technical work is also required to allow the estimation of $\p(\mathcal{G} \; \text{fails})$ (a discretisation argument) and to handle the term $\E( \frac{1}{\log x} \int_{-1/2}^{1/2} |F(\frac{1}{2}+it)|^2 )^{q'}$ that emerged from H\"{o}lder's inequality (we use an iterative procedure, an alternative would be to consider a sequence of different $K$ values).

The lower bounds in Theorems 1 and 2 are proved by comparing $\E \frac{1}{\log x} \int_{\mathcal{L}} |F(1/2+it)|^2 dt$ with $\E\left(\frac{1}{\log x} \int_{\mathcal{L}} |F(1/2+it)|^2 dt \right)^{2}$, where $\mathcal{L} \subseteq [-1/2,1/2]$ is a certain random subset chosen so that this second moment remains roughly the same size as the square of the first moment. In fact, we essentially choose $\mathcal{L}$ as the set of points at which $\log|F(1/2+it)| \leq \log\log x + 1/(1-q)$, and all subproducts of $F(1/2+it)$ obey a comparable bound. The idea here is that when one expands out the second moment, one needs to estimate terms $\E |F(1/2+it_1)|^2 |F(1/2+it_2)|^2$, and ideally one wants the answer to be approximately the product $\E |F(1/2+it_1)|^2 \E |F(1/2+it_2)|^2$. One cannot achieve exactly this when $t_1 - t_2$ is small because of dependencies between terms in the two products, but by restricting $t_1 , t_2$ to the set $\mathcal{L}$ one can ensure $\E |F(1/2+it_1)|^2 |F(1/2+it_2)|^2$ only blows up slowly as $t_1 - t_2$ becomes small, so one retains control when integrating over $t_1 , t_2$. The reader might think of this calculation as the underlying motivation for our choice of $\mathcal{G}$ in our upper bound proof, as well. The lower bound proofs are conceptually easier than the upper bounds, since we don't need any discretisation argument and don't need to control the contribution from points outside the nice set $\mathcal{L}$. However, on the probabilistic side we require a two dimensional Girsanov-type theorem allowing us to estimate $\E \textbf{1}_{A} |F(1/2+it_1)|^2 |F(1/2+it_2)|^2$, for certain events $A$ related to the definition of $\mathcal{L}$, and it requires quite a lot of work to set this up carefully.

\vspace{12pt}
Now we make a few remarks about multiplicative chaos: see the survey of Rhodes and Vargas~\cite{rhodesvargas} for more about the field. The theory begins with a random function $h(\cdot)$ on some space, say the interval $[-1/2,1/2]$. We construct a family of random measures on that space, depending on a parameter $\gamma > 0$, by defining the measure of a subset to be the integral of $e^{\gamma h(\cdot)}$ over the subset. In particular, the total measure $\int_{-1/2}^{1/2} e^{\gamma h(t)} dt$ is a random variable. This description is inaccurate in several respects, notably that $h(\cdot)$ is actually taken as a random {\em generalised} function, and one ``truncates'' or ``regularises'' $h$ on different scales (depending on a parameter $\epsilon$, say) to form a genuine random function $h_{\epsilon}$, and then studies the limit behaviour of the random measure as $\epsilon \rightarrow 0$. Lots of the probabilistic attention is devoted to showing that, for suitable $h(\cdot)$, different regularisation procedures result in the same limit behaviour. The most studied situation is Gaussian multiplicative chaos, where $h(t)$ is a mean zero Gaussian generalised function, with (regularised) variance approximately the same for each $t$ . To obtain non-trivial behaviour one needs to rescale $e^{\gamma h(t)}$ at each $t$ by its expected value $e^{(1/2) \gamma^2 \E h(t)^2}$. One also assumes that the collection of $h(t)$, or in fact their regularisations, have a certain logarithmic kind of covariance structure.

Let us compare with our situation, where we are interested in $\int_{-1/2}^{1/2} |F(1/2+it)|^2 dt = \int_{-1/2}^{1/2} e^{2\log|F(1/2+it)|} dt$. We can think of the length $x$ of the random product $F(1/2+it)$ as a regularisation parameter, but we are not here interested in constructing or analysing limit objects, so much as obtaining moment information about $\int_{-1/2}^{1/2} e^{2\log|F(1/2+it)|} dt$ uniformly in $x$. Our random function $\log|F(1/2+it)| = \sum_{p \leq x} \log|1 - \frac{f(p)}{p^{1/2+it}}|^{-1}$ is not Gaussian, but is approximately Gaussian as it is a sum of many independent components. When we rescale $\int_{-1/2}^{1/2} |F(1/2+it)|^2 dt$ by dividing by $\log x \asymp \E|F(1/2+it)|^2$, this directly corresponds to rescaling $e^{\gamma h(t)}$ by $e^{(1/2) \gamma^2 \E h(t)^2}$. It also turns out that the random variables $\log|F(1/2+it)|$ have a logarithmic covariance structure. This is why ideas from the theory of multiplicative chaos are relevant to $\int_{-1/2}^{1/2} |F(1/2+it)|^2 dt$.

Now the behaviour of Gaussian multiplicative chaos changes as $\gamma$ increases. There is a {\em critical value} $\gamma_c$ (depending e.g. on the dimension of the space one is working in and on the variance of $h(t)$) such that $\int e^{\gamma h_{\epsilon}(t) - (1/2) \gamma^2 \E h_{\epsilon}(t)^2}$ converges to a non-trivial limit measure as $\epsilon \rightarrow 0$ when $\gamma < \gamma_c$, but converges to the zero measure when $\gamma = \gamma_c$. It is natural to expect such a transition, because as $\gamma$ increases the dominant contribution to $\int e^{\gamma h_{\epsilon}(t) - (1/2) \gamma^2 \E h_{\epsilon}(t)^2}$ comes from larger values of $h_{\epsilon}(t)$, so when $\gamma$ is large enough one expects with high probability never to find such values on the whole range of integration. Number theory readers may again recognise a parallel with the analysis of moments of the Riemann zeta function. Critical multiplicative chaos, where $\gamma = \gamma_c$, is more difficult to analyse than the subcritical case, but recently Duplantier, Rhodes, Sheffield and Vargas~\cite{duprsvnorm} showed in some generality that if one replaces $\int e^{\gamma_c h_{\epsilon}(t) - (1/2) \gamma_{c}^2 \E h_{\epsilon}(t)^2}$ by $\int \sqrt{\E h_{\epsilon}(t)^2} e^{\gamma_c h_{\epsilon}(t) - (1/2) \gamma_{c}^2 \E h_{\epsilon}(t)^2}$, this converges to a non-zero limit measure. In our case the exponent 2 of $|F(1/2+it)|$ corresponds to critical $\gamma$ (see our earlier discussion about the main contribution to our integral coming when $\log|F(1/2+it)| \approx \log\log x$), and  we have $\E (\log|F(1/2+it)|)^2 \asymp \log\log x$, so the factor $\sqrt{\E h_{\epsilon}(t)^2}$ directly suggests the factor $\sqrt{\log\log x}$ in the denominator in Theorems 1 and 2, for $q$ away from 1.

We end with specific connections between our problem and the multiplicative chaos literature. The most relevant work is a preprint of Saksman and Webb~\cite{saksmanwebb}, see also \cite{saksmanwebb2}, showing a random model for $\log\zeta(1/2+it)$ is well approximated by a perturbation of a Gaussian random field, so one can apply many results about Gaussian multiplicative chaos to the multiplicative chaos arising from that random model as well. The model for $\log\zeta(1/2+it)$ is very close to $\log|F(1/2+it)|$ in the Steinhaus case, so it is possible that combining (the rigorous version of) \eqref{introprod}, Saksman and Webb's approximation, and the results of Duplantier, Rhodes, Sheffield and Vargas~\cite{duprsvnorm} about moments of the total measure of critical Gaussian chaos (which ultimately stem from Kahane's convexity inequality and results of Hu and Shi~\cite{hushi} for branching random walk), one could get another proof of Theorem 1 for $q$ bounded away from 1. It is not clear whether one could prove the full uniform version of Theorem 1, where $q=q(x)$ may be close to 1, in this way. In any case, our proofs here are self contained. One inspiration for our arguments is a beautiful short paper of Berestycki~\cite{berestycki}, proving convergence results for subcritical Gaussian chaos using Girsanov's theorem and restricting attention to high probability good events. We are in the more delicate critical case, in a non-Gaussian setting, and trying to prove different kinds of results, but several features carry over.

We should also mention the connection between multiplicative chaos and the maxima of logarithmically correlated random processes, corresponding to $h(t)$ in our earlier discussion. This interaction is well known to workers in the area (see section 6.4 of Rhodes and Vargas~\cite{rhodesvargas}, for example), and the arguments of Berestycki~\cite{berestycki} build on approaches to analysing the maxima of such processes. Roughly speaking, one could perhaps say that multiplicative chaos is a bit easier to analyse (because exponentiating boosts the role of very large values and somewhat dilutes the interaction at nearby $t$), but on the other hand one is interested in proving sharper results about multiplicative chaos (e.g. bounds that are sharp up to constants). See the paper of Arguin, Belius and Harper~\cite{abh} for results on the maximum of a random model for $\log\zeta(1/2+it)$, which again is close to $\log|F(1/2+it)|$ in the Steinhaus case.

\subsection{Some corollaries}
The following is an immediate corollary of our theorems.
\begin{cor1}
If $f(n)$ is a Steinhaus or Rademacher random multiplicative function, then
$$ \frac{\sum_{n \leq x} f(n)}{\sqrt{\E|\sum_{n \leq x} f(n)|^2}} \stackrel{p}{\rightarrow} 0 \;\;\; \text{as} \; x \rightarrow \infty , $$
where $\stackrel{p}{\rightarrow}$ denotes convergence in probability.
\end{cor1}

\begin{proof}[Proof of Corollary 1]
As noted earlier, we have $\sqrt{\E|\sum_{n \leq x} f(n)|^2} \asymp \sqrt{x}$, whereas we have $\E|\sum_{n \leq x} f(n)| \asymp \sqrt{x}/(\log\log x)^{1/4}$ by Theorem 1 in the Steinhaus case, or Theorem 2 in the Rademacher case. Thus $\frac{\sum_{n \leq x} f(n)}{\sqrt{\E|\sum_{n \leq x} f(n)|^2}}$ converges to 0 in $L^{1}$, which is a stronger statement than convergence to 0 in probability.
\end{proof}

Corollary 1 resolves the question of the limiting distribution of $\frac{\sum_{n \leq x} f(n)}{\sqrt{\E|\sum_{n \leq x} f(n)|^2}}$, which previously has generated quite a lot of discussion, exploration of the behaviour when one conditions on $f(p)$ on small primes $p$, and numerical simulations (see the papers of Chatterjee and Soundararajan~\cite{chatsound} and Hough~\cite{hough}, for example), as well as Harper's~\cite{harperlimits} proof that the distribution is not $N(0,1)$ in the Rademacher case. One now has an obvious and interesting follow-up question, namely what is the distribution of the properly renormalised sum $\frac{\sum_{n \leq x} f(n)}{\sqrt{x}/(\log\log x)^{1/4}}$? Our arguments here may be a reasonable starting point for investigating this, since they show that matters substantially reduce to understanding the distribution of $\frac{1}{\log x} \int_{-1/2}^{1/2} |F(1/2+it)|^2 dt$, which we can try to access by understanding the distribution of the total mass of critical multiplicative chaos.

\vspace{12pt}
By Chebychev's inequality, we immediately have
$$ \p(|\sum_{n \leq x} f(n)| \geq \lambda\sqrt{x}) \leq \frac{\E|\sum_{n \leq x} f(n)|^2}{(\lambda \sqrt{x})^2} \leq \frac{1}{\lambda^2} \;\;\; \forall \; \lambda > 0 . $$
It is natural to want to improve this, and one might even hope to obtain exponential decay for large $\lambda$, by analogy with large deviation estimates for sums of independent random variables. However, so far as the author is aware, no improved estimate whatsoever is known for fixed $x$ and moderately large $\lambda$ (for $\lambda$ larger than a suitable power of $\log x$, one gets a better bound by looking e.g. at the fourth moment $\E|\sum_{n \leq x} f(n)|^4$). In the following corollary we obtain a small improvement of the Chebychev upper bound. We also show, perhaps surprisingly, that this is close to best possible on a wide range of $\lambda$.

\begin{cor2}
Let $x$ be large, and let $f(n)$ be a Steinhaus or Rademacher random multiplicative function. For all $\lambda \geq 2$, we have
$$ \p(|\sum_{n \leq x} f(n)| \geq \lambda\frac{\sqrt{x}}{(\log\log x)^{1/4}}) \ll \frac{\min\{\log\lambda, \sqrt{\log\log x}\}}{\lambda^2} . $$

In addition, for all $2 \leq \lambda \leq e^{\sqrt{\log\log x}}$ we have
$$ \p(|\sum_{n \leq x} f(n)| \geq \lambda\frac{\sqrt{x}}{(\log\log x)^{1/4}}) \gg \frac{1}{\lambda^2 (\log\log x)^{O(1)}} . $$
\end{cor2}

The upper bound here is a rather direct corollary of Theorems 1 and 2, and we shall prove it immediately.
\begin{proof}[Proof of Corollary 2, upper bound]
If $\lambda \geq e^{\sqrt{\log\log x}}$, then the result follows from Chebychev's inequality applied to $\E|\sum_{n \leq x} f(n)|^2$.

For smaller $\lambda$, for any $q \leq 1$ we have
$$ \p(|\sum_{n \leq x} f(n)| \geq \lambda\frac{\sqrt{x}}{(\log\log x)^{1/4}}) \leq \frac{\E|\sum_{n \leq x} f(n)|^{2q}}{(\lambda \frac{\sqrt{x}}{(\log\log x)^{1/4}})^{2q}} \ll \frac{1}{\lambda^{2q}} \frac{(\log\log x)^{q/2}}{(1 + (1-q)\sqrt{\log\log x})^q} . $$
If we set $q = 1-\delta$ with $0 \leq \delta \leq 1$, then the right hand side is
$$ \frac{1}{\lambda^2} \frac{\lambda^{2\delta} (\log\log x)^{q/2}}{(1 + \delta\sqrt{\log\log x})^q} \leq \frac{1}{\lambda^2} \frac{\lambda^{2\delta}}{\delta} = \frac{\log\lambda}{\lambda^2} \frac{e^{2\delta \log\lambda}}{\delta \log\lambda}  . $$
Choosing $\delta = \frac{1}{2\log\lambda}$ yields the claimed upper bound.
\end{proof}

Proving the lower bound in Corollary 2 requires some additional ideas, so its proof is deferred to section \ref{corsec}. Very roughly speaking, since the value of $|F(1/2+it)|^2$ doesn't change much on $t$ intervals of length $1/\log x$, and since we have something like $|\sum_{n \leq x} f(n)| \approx \sqrt{x} (\frac{1}{\log x} \int_{-1/2}^{1/2} |F(1/2+it)|^2 dt)^{1/2}$, it will suffice to prove that
$$ \p(\max_{|t| \leq 1/2} |F(1/2+it)|^2 \geq \lambda^2 \log^{2}x) \gg \frac{1}{\lambda^2 (\log\log x)^{O(1)}} \;\;\; \forall \; 2 \leq \lambda \leq e^{\sqrt{\log\log x}} . $$
This can be done by approximating the random variables $(|F(1/2+it)|^2)_{|t| \leq 1/2}$ by the exponentials of certain correlated Gaussian random variables, and applying known results about the maximum of such Gaussians. 

It seems another interesting open question to determine the exact magnitude of the probabilities in Corollary 2. They may be $\asymp 1/\lambda^2$ for any fixed $\lambda$, since Barral, Kupiainen, Nikula, Saksman and Webb~\cite{barknsw} have shown that the limiting total measure of critical multiplicative chaos has upper tails of the shape $\asymp 1/\lambda$, which would suggest tails $\asymp 1/\lambda^2$ here because of the power 1/2 in the approximation of $|\sum_{n \leq x} f(n)|$ above. On a wide range of $\lambda \geq e^{\sqrt{\log\log x}}$, the proof of the lower bound in Corollary 2 still yields a result, of the form $e^{-(\log\lambda)^{2}/\log\log x}/(\lambda^2 (\log\log x)^{O(1)})$. This order of magnitude may be essentially the correct answer for $\lambda \geq e^{\sqrt{\log\log x}}$. It is unclear what one should expect as the correct answer between these ranges.

\subsection{Organisational remarks}
In section \ref{seceulerprod}, we make rigorous the statement that $\E|\sum_{n \leq x} f(n)|^{2q} \approx x^{q} \E(\frac{1}{\log x} \int_{-1/2}^{1/2} |F(1/2+it)|^2 dt)^{q}$. In section \ref{secprobcalc}, which is the longest section, we obtain various probabilistic estimates for $\E|F(1/2+it)|^2$ and for $\E \textbf{1}_{A} |F(1/2+it)|^2$, for certain kinds of events $A$. Section \ref{secmainupper} contains the fairly quick deduction of the upper bound parts of Theorems 1 and 2. Section \ref{secmainlower} includes some further probabilistic estimates, generalising the work from section \ref{secprobcalc} to two dimensions, and then the deduction of the lower bound parts of Theorems 1 and 2. Finally, section \ref{corsec} proves the lower bound from Corollary 2, and the appendix gives proofs of two Probability Results on Gaussian random walks, deferred from section \ref{secprobcalc}.

Note that it will suffice to prove Theorems 1 and 2 for $2/3 \leq q \leq 1$. For the upper bounds, if we know that $\E|\sum_{n \leq x} f(n)|^{4/3} \ll \frac{x^{2/3}}{(\log\log x)^{1/3}}$ (the result with $q=2/3$) then H\"{o}lder's inequality yields
$$ \E|\sum_{n \leq x} f(n)|^{2q} \leq \Biggl( \E|\sum_{n \leq x} f(n)|^{4/3} \Biggr)^{3q/2} \ll \frac{x^{q}}{(\log\log x)^{q/2}} \;\;\; \forall \; 0 \leq q \leq 2/3 , $$
as desired. For the lower bounds, if we know that $\E|\sum_{n \leq x} f(n)|^{4/3} \asymp \frac{x^{2/3}}{(\log\log x)^{1/3}}$ and that $\E|\sum_{n \leq x} f(n)|^{3/2} \asymp \frac{x^{3/4}}{(\log\log x)^{3/8}}$ then, by H\"{o}lder's inequality,
\begin{eqnarray}
\frac{x^{2/3}}{(\log\log x)^{1/3}} \asymp \E|\sum_{n \leq x} f(n)|^{4/3} & = & \E|\sum_{n \leq x} f(n)|^{\frac{q/3}{3/2 - 2q}} |\sum_{n \leq x} f(n)|^{\frac{2 - 3q}{3/2 - 2q}} \nonumber \\
& \leq & \Biggl( \E|\sum_{n \leq x} f(n)|^{2q} \Biggr)^{\frac{1}{6(3/2 - 2q)}} \Biggl( \E|\sum_{n \leq x} f(n)|^{3/2} \Biggr)^{\frac{4-6q}{3(3/2 - 2q)}} \nonumber \\
& \ll & \Biggl( \E|\sum_{n \leq x} f(n)|^{2q} \Biggr)^{\frac{1}{6(3/2 - 2q)}} \Biggl( \frac{x}{\sqrt{\log\log x}} \Biggr)^{\frac{2-3q}{3 - 4q}} \nonumber
\end{eqnarray}
for all $0 \leq q \leq 2/3$, which implies that $\E|\sum_{n \leq x} f(n)|^{2q} \gg (\frac{x}{\sqrt{\log\log x}})^{q}$. Restricting to $2/3 \leq q \leq 1$ will be a useful simplification in our main arguments, making various series involving $q$ converge and giving us access to Minkowski's inequality in certain places.

We finish with a remark on notation and references. As usual, we will say a number $n$ is $y$-smooth if all prime factors of $n$ are $\leq y$. We will generally use the letter $p$ to denote primes. Unless mentioned otherwise, the letters $c, C$ will be used to denote positive constants, $c$ usually being a small constant and $C$ a large one. We write $f(x) = O(g(x))$ and $f(x) \ll g(x)$, both of which mean that there exists $C$ such that $|f(x)| \leq Cg(x)$, for all $x$. In a few places this notation will be adorned with a subscript parameter (e.g. $O_{\epsilon}(\cdot)$ and $\ll_{\delta}$), meaning that the implied constant $C$ is allowed to depend on that parameter. We write $f(x) \asymp g(x)$ to mean that $g(x) \ll f(x) \ll g(x)$, in other words that $c g(x) \leq |f(x)| \leq C g(x)$ for some $c,C$, for all $x$.

Less standard results that we use are stated explicitly as results in the text, and we try also to give references in line for most of the standard material in number theory and probability, to help readers who are less familiar with one area or the other. The books of Gut~\cite{gut} and of Montgomery and Vaughan~\cite{mv} may be consulted as excellent general references for probabilistic and number theoretic background.

\section{The reduction to Euler products}\label{seceulerprod}
In this section we shall prove four Propositions that make precise and rigorous the assertion in \eqref{introprod}, that $\E|\sum_{n \leq x} f(n)|^{2q}$ may be bounded by studying integrals of Euler products.

\subsection{Some lemmas}
We begin by recording two fairly straightforward lemmas we need.

\begin{numth1}[See Lemma 2.1 of Lau, Tenenbaum and Wu~\cite{tenenbaum}]
Let $0 < \delta < 1$, let $m \geq 1$, and suppose that $\max\{3, 2m\} \leq y \leq z \leq y^{10}$ and that $1 < u \leq v(1 - y^{-\delta})$. As usual, let $\Omega(d)$ denote the total number of prime factors of $d$ (counted with multiplicity). Then
$$ \sum_{\substack{u \leq d \leq v, \\ p | d \Rightarrow y \leq p \leq z}} m^{\Omega(d)} \ll_{\delta} \frac{(v-u) m}{\log y} \prod_{y \leq p \leq z} \left(1 - \frac{m}{p}\right)^{-1} . $$

In particular, for any large $x$ and any $k \leq \log\log x - 5$ we have
$$ \sum_{\substack{d \leq x, \\ p | d \Rightarrow x^{e^{-(k+1)}} \leq p \leq x^{e^{-k}}}} 1 \ll 2^{-e^{k}} \frac{x}{\log x} . $$
\end{numth1}

The first statement here is a slight generalisation of Lemma 2.1 of Lau, Tenenbaum and Wu~\cite{tenenbaum} (see also Lemma 3 of Hal\'{a}sz~\cite{halasz}), in which the sum over $d$ was restricted to squarefree numbers. The proof is short, so we give it in full.
\begin{proof}[Proof of Number Theory Result 1]
We have
$$ \sum_{\substack{u \leq d \leq v, \\ p | d \Rightarrow y \leq p \leq z}} m^{\Omega(d)} \leq m \sum_{\substack{y \leq q \leq z, \\ q \; \text{prime}}} \sum_{\substack{u/q \leq d \leq v/q, \\ p | d \Rightarrow y \leq p \leq z}} m^{\Omega(d)} \leq m \sum_{\substack{u/z \leq d \leq v/y, \\ p | d \Rightarrow y \leq p \leq z}} m^{\Omega(d)} \sum_{\substack{u/d \leq q \leq v/d, \\ q \; \text{prime}}} 1 . $$
Now $(v-u)/d \geq y(v-u)/v \geq y^{1-\delta}$ here, whilst $u/d \leq z \leq y^{10}$, so the Brun--Titchmarsh upper bound for primes in intervals is available (see e.g. Theorem 3.9 of Montgomery and Vaughan~\cite{mv}) and implies that $\sum_{\substack{u/d \leq q \leq v/d, \\ q \; \text{prime}}} 1 \ll_{\delta} (v-u)/(d \log y)$. The first part of the result now follows on inserting this estimate and noting that
$$ \sum_{\substack{u/z \leq d \leq v/y, \\ p | d \Rightarrow y \leq p \leq z}} \frac{m^{\Omega(d)}}{d} \leq \sum_{\substack{d : \\ p | d \Rightarrow y \leq p \leq z}} \frac{m^{\Omega(d)}}{d} = \prod_{y \leq p \leq z} \left(1 - \frac{m}{p}\right)^{-1} , $$
the product certainly being well defined in view of our assumption that $y \geq 2m$.

To deduce the second part of the result, we note that
$$ \sum_{\substack{d \leq x, \\ p | d \Rightarrow x^{e^{-(k+1)}} \leq p \leq x^{e^{-k}}}} 1 \leq \sqrt{x} + 5^{-e^{k}/2}  \sum_{\substack{\sqrt{x} \leq d \leq x, \\ p | d \Rightarrow x^{e^{-(k+1)}} \leq p \leq x^{e^{-k}}}} 5^{\Omega(d)} , $$
since a number $d \geq \sqrt{x}$ with all its prime factors smaller than $x^{e^{-k}}$ must have at least $e^{k}/2$ such factors (counted with multiplicity). We can apply the first part of the result with $\delta = 1/2$, say, and deduce that
\begin{eqnarray}
5^{-e^{k}/2}  \sum_{\substack{\sqrt{x} \leq d \leq x, \\ p | d \Rightarrow x^{e^{-(k+1)}} \leq p \leq x^{e^{-k}}}} 5^{\Omega(d)} & \ll & 5^{-e^{k}/2} \frac{x}{\log(x^{e^{-(k+1)}})} \prod_{x^{e^{-(k+1)}} \leq p \leq x^{e^{-k}}}\left(1 - \frac{5}{p}\right)^{-1} \nonumber \\
& \ll & 5^{-e^{k}/2} e^{k} \frac{x}{\log x} \ll 2^{-e^{k}} \frac{x}{\log x} . \nonumber
\end{eqnarray}
Here the Mertens estimate (see e.g. Theorem 2.7 of Montgomery and Vaughan~\cite{mv}) implies the product over primes is $\ll 1$. Since $2^{-e^{k}} \frac{x}{\log x}$ is always larger than $\sqrt{x}$ on our range of $k$, the result follows.
\end{proof}

We will need the following version of Parseval's identity for Dirichlet series.
\begin{harman1}[See (5.26) in sec. 5.1 of Montgomery and Vaughan~\cite{mv}]
Let $(a_n)_{n=1}^{\infty}$ be any sequence of complex numbers, and let $A(s) := \sum_{n=1}^{\infty} \frac{a_n}{n^s}$ denote the corresponding Dirichlet series, and $\sigma_c$ denote its abscissa of convergence. Then for any $\sigma > \max\{0,\sigma_c \}$, we have
$$ \int_{0}^{\infty} \frac{|\sum_{n \leq x} a_n |^2}{x^{1 + 2\sigma}} dx = \frac{1}{2\pi} \int_{-\infty}^{\infty} \left|\frac{A(\sigma + it)}{\sigma + it}\right|^2 dt . $$
\end{harman1}

\subsection{Upper bounds: statement of the propositions}
We will need a little notation. Given a random multiplicative function $f(n)$ (either Steinhaus or Rademacher, depending on the context), and an integer $0 \leq k \leq \log\log x$, let $F_k$ denote the partial Euler product of $f(n)$ over $x^{e^{-(k+1)}}$-smooth numbers. Thus for all complex $s$ with $\Re(s) > 0$, we have
$$ F_{k}(s) = \prod_{p \leq x^{e^{-(k+1)}}} \left(1 - \frac{f(p)}{p^s}\right)^{-1} = \sum_{\substack{n=1, \\ n \; \text{is} \; x^{e^{-(k+1)}} \; \text{smooth}}}^{\infty} \frac{f(n)}{n^s} $$
in the Steinhaus case, and
$$ F_{k}(s) = \prod_{p \leq x^{e^{-(k+1)}}} \left(1 + \frac{f(p)}{p^s}\right) = \sum_{\substack{n=1, \\ n \; \text{is} \; x^{e^{-(k+1)}} \; \text{smooth}}}^{\infty} \frac{f(n)}{n^s} $$
in the Rademacher case (the product taking a different form because $f(n)$ is only supported on squarefree numbers in that case).

\begin{prop1}
Let $f(n)$ be a Steinhaus random multiplicative function, let $x$ be large, and set $\mathcal{K} := \lfloor \log\log\log x \rfloor$. Uniformly for all $2/3 \leq q \leq 1$, we have
$$ \vert\vert \sum_{n \leq x} f(n) \vert\vert_{2q} \ll \sqrt{\frac{x}{\log x}} \sum_{0 \leq k \leq \mathcal{K}} \vert\vert \int_{-1/2}^{1/2} |F_k(1/2 - \frac{k}{\log x} + it)|^2 dt \vert\vert_{q}^{1/2} + \sqrt{\frac{x}{\log x}} . $$
\end{prop1}

It is perhaps worth emphasising that $\vert\vert \cdot \vert\vert_{r} := (\E|\cdot|^{r})^{1/r}$ is a genuine norm when $r \geq 1$, but not for smaller $r$. Thus $\vert\vert \cdot \vert\vert_{2q}$ is a genuine norm on our range of $q$ (so we may apply Minkowski's inequality to it, as we shall in the proof), but $\vert\vert \cdot \vert\vert_{q}$ is not and Minkowski's inequality is not applicable.

\begin{prop2}
Let $f(n)$ be a Rademacher random multiplicative function, let $x$ be large, and set $\mathcal{K} := \lfloor \log\log\log x \rfloor$. Uniformly for all $2/3 \leq q \leq 1$, we have
\begin{eqnarray}
\vert\vert \sum_{n \leq x} f(n) \vert\vert_{2q} & \ll & \sqrt{\frac{x}{\log x}} \sum_{0 \leq k \leq \mathcal{K}} \max_{N \in \Z} \frac{1}{(|N|+1)^{1/8}} \vert\vert \int_{N-1/2}^{N+1/2} |F_k(1/2 - \frac{k}{\log x} + it)|^2 dt \vert\vert_{q}^{1/2} \nonumber \\
&& + \sqrt{\frac{x}{\log x}} . \nonumber
\end{eqnarray}
\end{prop2}

In these bounds we expect, and it will turn out to be the case, that the main contribution comes from small $k$ and small $N$. In the Steinhaus case, for any fixed $t \in \R$ the distribution of $(f(n)n^{it})$ is the same as the distribution of $(f(n))$, which is why one doesn't need to deal with translates by $N$ in the Steinhaus case. This is one of a few differences between Rademacher and Steinhaus random multiplicative functions that will recur a number of times in our analysis.

\subsection{Lower bounds: statement of the propositions}
For our work on lower bounds, we again connect the size of $\vert\vert \sum_{n \leq x} f(n) \vert\vert_{2q}$ with a certain integral average, and thence with random Euler products. Let $F$ denote the partial Euler product of $f(n)$, either Steinhaus or Rademacher, over $x$-smooth numbers. (Thus $F = F_{-1}$, if we slightly abuse our earlier notation).

\begin{prop3}
There exists a large absolute constant $C > 0$ such that the following is true. If $f(n)$ is a Steinhaus random multiplicative function, and $x$ is large, then uniformly for all $2/3 \leq q \leq 1$ we have
$$ \vert\vert \sum_{n \leq x} f(n) \vert\vert_{2q} \gg \sqrt{\frac{x}{\log x}} \vert\vert \int_{1}^{\sqrt{x}} \left|\sum_{m \leq z} f(m) \right|^2 \frac{dz}{z^{2}} \vert\vert_{q}^{1/2} - C\sqrt{\frac{x}{\log x}} . $$

In particular, for any large quantity $V$ we have that $\vert\vert \sum_{n \leq x} f(n) \vert\vert_{2q}$ is
$$ \gg \sqrt{\frac{x}{\log x}} \Biggl( \vert\vert \int_{-1/2}^{1/2} |F(1/2 + \frac{4V}{\log x} + it)|^2 dt \vert\vert_{q}^{1/2} - \frac{C}{e^V} \vert\vert \int_{-1/2}^{1/2} |F(1/2 + \frac{2V}{\log x} + it)|^2 dt \vert\vert_{q}^{1/2}  - C \Biggr) . $$
\end{prop3}

\begin{prop4}
If $f(n)$ is a Rademacher random multiplicative function, the first bound in Proposition 10 continues to hold, and the second bound may be replaced by the statement that
\begin{eqnarray}
\vert\vert \sum_{n \leq x} f(n) \vert\vert_{2q} & \gg & \sqrt{\frac{x}{\log x}} \Biggl( \vert\vert \int_{-1/2}^{1/2} |F(1/2 + \frac{4V}{\log x} + it)|^2 dt \vert\vert_{q}^{1/2} - \nonumber \\
&& - \frac{C}{e^V} \max_{N \in \Z} \frac{1}{(|N|+1)^{1/8}} \vert\vert \int_{N-1/2}^{N+1/2} |F(1/2 + \frac{2V}{\log x} + it)|^2 dt \vert\vert_{q}^{1/2}  - C \Biggr) . \nonumber
\end{eqnarray}
\end{prop4}

When we come to apply these Propositions, we will choose $V$ to be a sufficiently large fixed constant that, because of the factor $C/e^{V}$, the second subtracted Euler product integral is negligible compared with the first. 

\subsection{Proof of Propositions 1 and 2}
We begin with Proposition 1. Let $P(n)$ denote the largest prime factor of $n$, and let $\Psi(x,y) := \#\{n \leq x : n \; \text{is} \; y \; \text{smooth}\} = \#\{n \leq x : P(n) \leq y \}$. Then by Minkowski's inequality, for all $2/3 \leq q \leq 1$ we have
$$ \vert\vert \sum_{n \leq x} f(n) \vert\vert_{2q} \leq \sum_{0 \leq k \leq \mathcal{K}} \vert\vert \sum_{\substack{n \leq x, \\ x^{e^{-(k+1)}} < P(n) \leq x^{e^{-k}}}} f(n) \vert\vert_{2q} + \vert\vert \sum_{\substack{n \leq x, \\ P(n) \leq x^{e^{-(\mathcal{K}+1)}} }} f(n) \vert\vert_{2q} . $$
Furthermore, by H\"{o}lder's inequality we have
$$ \vert\vert \sum_{\substack{n \leq x, \\ P(n) \leq x^{e^{-(\mathcal{K}+1)}} }} f(n) \vert\vert_{2q} \leq \vert\vert \sum_{\substack{n \leq x, \\ P(n) \leq x^{e^{-(\mathcal{K}+1)}} }} f(n) \vert\vert_{2} = \Psi(x,x^{e^{-(\mathcal{K}+1)}})^{1/2} , $$
and recalling that $\mathcal{K} := \lfloor \log\log\log x \rfloor$ as well as standard estimates for smooth numbers (see Theorem 7.6 of Montgomery and Vaughan~\cite{mv}, for example), the above is $\leq \Psi(x,x^{1/\log\log x})^{1/2} \ll \sqrt{x} (\log x)^{-c\log\log\log x}$. This contribution is more than acceptable.

If we let $\E^{(k)}$ denote expectation conditional on $(f(p))_{p \leq x^{e^{-(k+1)}}}$, and use H\"{o}lder's inequality and a mean square calculation, we see $\sum_{0 \leq k \leq \mathcal{K}} \vert\vert \sum_{\substack{n \leq x, \\ x^{e^{-(k+1)}} < P(n) \leq x^{e^{-k}}}} f(n) \vert\vert_{2q}$ is
\begin{eqnarray}
& = & \sum_{0 \leq k \leq \mathcal{K}} \vert\vert \sum_{\substack{1 < m \leq x , \\ p|m \Rightarrow x^{e^{-(k+1)}} < p \leq x^{e^{-k}}}} f(m) \sum_{\substack{n \leq x/m, \\ n \; \text{is} \; x^{e^{-(k+1)}} \text{-smooth}}} f(n)  \vert\vert_{2q} . \nonumber \\
& = & \sum_{0 \leq k \leq \mathcal{K}} \Biggl( \E \E^{(k)}\Biggl|\sum_{\substack{1 < m \leq x , \\ p|m \Rightarrow x^{e^{-(k+1)}} < p \leq x^{e^{-k}}}} f(m) \sum_{\substack{n \leq x/m, \\ n \; \text{is} \; x^{e^{-(k+1)}} \text{-smooth}}} f(n) \Biggr|^{2q} \Biggr)^{1/2q} \nonumber \\
& \leq & \sum_{0 \leq k \leq \mathcal{K}} \Biggl( \E( \E^{(k)}\Biggl|\sum_{\substack{1 < m \leq x , \\ p|m \Rightarrow x^{e^{-(k+1)}} < p \leq x^{e^{-k}}}} f(m) \sum_{\substack{n \leq x/m, \\ n \; \text{is} \; x^{e^{-(k+1)}} \text{-smooth}}} f(n) \Biggr|^2)^{q} \Biggr)^{1/2q} \nonumber \\
& = & \sum_{0 \leq k \leq \mathcal{K}} \vert\vert \sum_{\substack{1 < m \leq x , \\ p|m \Rightarrow x^{e^{-(k+1)}} < p \leq x^{e^{-k}}}} \Biggl|\sum_{\substack{n \leq x/m, \\ n \; \text{is} \; x^{e^{-(k+1)}} \text{-smooth}}} f(n) \Biggr|^2 \vert\vert_{q}^{1/2} . \nonumber
\end{eqnarray}

To proceed further, we want to replace $\Biggl|\sum_{\substack{n \leq x/m, \\ n \; \text{is} \; x^{e^{-(k+1)}} \text{-smooth}}} f(n) \Biggr|^2$ in the above by a smoothed version. Set $X = e^{\sqrt{\log x}}$, say, and note that (uniformly for any $2/3 \leq q \leq 1$ and any $0 \leq k \leq \mathcal{K}$) we have
\begin{eqnarray}
&& \E \Biggl(\sum_{\substack{1 < m \leq x , \\ p|m \Rightarrow x^{e^{-(k+1)}} < p \leq x^{e^{-k}}}} \Biggl| \sum_{\substack{n \leq x/m, \\ n \; \text{is} \; x^{e^{-(k+1)}} \text{-smooth}}} f(n) \Biggr|^2 \Biggr)^{q} \nonumber \\
& \ll & \E \Biggl( \sum_{\substack{1 < m \leq x , \\ p|m \Rightarrow x^{e^{-(k+1)}} < p \leq x^{e^{-k}}}} \frac{X}{m} \int_{m}^{m(1 + \frac{1}{X})} \Biggl| \sum_{\substack{n \leq x/t, \\ x^{e^{-(k+1)}} \text{-smooth}}} f(n) \Biggr|^2 dt \Biggr)^{q} + \nonumber \\
&& + \E \Biggl(\sum_{\substack{1 < m \leq x , \\ p|m \Rightarrow x^{e^{-(k+1)}} < p \leq x^{e^{-k}}}} \frac{X}{m} \int_{m}^{m(1 + \frac{1}{X})} \Biggl| \sum_{\substack{x/t < n \leq x/m, \\ x^{e^{-(k+1)}} \text{-smooth}}} f(n) \Biggr|^2 dt \Biggr)^{q} . \nonumber
\end{eqnarray}
Using H\"{o}lder's inequality and a mean square calculation, the second term is at most
$$ \Biggl(\sum_{\substack{1 < m \leq x , \\ p|m \Rightarrow x^{e^{-(k+1)}} < p \leq x^{e^{-k}}}} \frac{X}{m} \int_{m}^{m(1 + \frac{1}{X})} \E \Biggl| \sum_{\substack{x/t < n \leq x/m, \\ x^{e^{-(k+1)}} \text{-smooth}}} f(n) \Biggr|^2 dt \Biggr)^{q} \leq \Biggl(\sum_{\substack{1 < m \leq x , \\ p|m \Rightarrow x^{e^{-(k+1)}} < p \leq x^{e^{-k}}}} (1 + \frac{x}{mX}) \Biggr)^{q} . $$
In particular, using the second part of Number Theory Result 1 we obtain that the sum $\sum_{\substack{1 < m \leq x , \\ p|m \Rightarrow x^{e^{-(k+1)}} < p \leq x^{e^{-k}}}} 1 \ll 2^{-e^{k}} \frac{x}{\log x}$, and we have $\sum_{m \leq x} \frac{x}{mX} \ll \frac{x\log x}{X} \ll 2^{-e^{k}} \frac{x}{\log x}$ on our range $0 \leq k \leq \mathcal{K}$ as well. So taking $2q$-th roots and summing over $0 \leq k \leq \mathcal{K}$ leads to an acceptable overall contribution $\ll \sqrt{x/\log x}$ in Proposition 1.

Meanwhile, by swapping the order of the sum and integral we see the first term in the above is
$$ \E \Biggl( \int_{x^{e^{-(k+1)}}}^{x} \Biggl| \sum_{\substack{n \leq x/t, \\ x^{e^{-(k+1)}} \text{-smooth}}} f(n) \Biggr|^2 \sum_{\substack{t/(1+1/X) \leq m \leq t , \\ p|m \Rightarrow x^{e^{-(k+1)}} < p \leq x^{e^{-k}}}} \frac{X}{m} dt \Biggr)^{q} , $$
and a standard sieve estimate (sieving out by all primes in $[2,t^{1/10}] \backslash (x^{e^{-(k+1)}}, x^{e^{-k}}]$, say) shows this has order at most
$$ \E \Biggl( \int_{x^{e^{-(k+1)}}}^{x} \Biggl| \sum_{\substack{n \leq x/t, \\ x^{e^{-(k+1)}} \text{-smooth}}} f(n) \Biggr|^2 \frac{dt}{\log t} \Biggr)^q = x^q \E \Biggl( \int_{1}^{x^{1 - e^{-(k+1)}}} \Biggl| \sum_{\substack{n \leq z, \\ x^{e^{-(k+1)}} \text{-smooth}}} f(n) \Biggr|^2 \frac{dz}{z^2 \log(\frac{x}{z})} \Biggr)^q . $$
To obtain the second expression here, we made a substitution $z = x/t$ in the integral.

To obtain a satisfactory dependence on $k$ in our final estimations, we now note that if $z \leq \sqrt{x}$ we have $\log(x/z) \gg \log x$, whereas if $\sqrt{x} < z \leq x^{1 - e^{-(k+1)}}$ we have $\log(x/z) \gg e^{-k} \log x$. Thus in any case we have $\log(x/z) \gg z^{-2k/\log x} \log x$, so the above is
$$ \ll \frac{x^q}{\log^{q}x} \E \Biggl( \int_{1}^{x^{1 - e^{-(k+1)}}} \Biggl| \sum_{\substack{n \leq z, \\ x^{e^{-(k+1)}} \text{-smooth}}} f(n) \Biggr|^2 \frac{dz}{z^{2- 2k/\log x}} \Biggr)^q . $$
Finally, recalling that $F_k$ denotes the partial Euler product of $f(n)$ over $x^{e^{-(k+1)}}$-smooth numbers, and that $2/3 \leq q \leq 1$, we can apply Harmonic Analysis Result 1 to deduce the expectation above is
\begin{eqnarray}
\leq \E \Biggl( \int_{-\infty}^{\infty} \frac{|F_{k}(\frac{1}{2} - \frac{k}{\log x} + it)|^2}{|1/2 - \frac{k}{\log x} + it|^2} dt \Biggr)^q & \leq & \sum_{n \in \Z} \E \Biggl( \int_{n-1/2}^{n+1/2} \frac{|F_{k}(\frac{1}{2} - \frac{k}{\log x} + it)|^2}{|1/2 - \frac{k}{\log x} + it|^2} dt \Biggr)^q \nonumber \\
& \ll & \sum_{n \in \Z} \frac{1}{|n|^{2q} + 1} \E \Biggl( \int_{n-1/2}^{n+1/2} |F_{k}(\frac{1}{2} - \frac{k}{\log x} + it)|^2 dt \Biggr)^q . \nonumber
\end{eqnarray}
In the Steinhaus case, since the law of the random function $f(n)$ is the same as the law of $f(n)n^{it}$ for any fixed $t \in \R$ we have
$$  \E \Biggl( \int_{n-1/2}^{n+1/2} |F_{k}(1/2 - \frac{k}{\log x} + it)|^2 dt \Biggr)^q = \E \Biggl( \int_{-1/2}^{1/2} |F_{k}(1/2 - \frac{k}{\log x} + it)|^2 dt \Biggr)^q \;\;\; \forall \; n . $$
Proposition 1 now follows on putting everything together.
\qed

\vspace{12pt}
The proof of Proposition 2, covering the Rademacher case, is extremely similar to the Steinhaus case. Indeed, the only non-trivial change comes at the very end, where (since it is no longer the case that the law of the random function $f(n)$ is the same as the law of $f(n)n^{it}$) we apply the bound
\begin{eqnarray}
&& \sum_{n \in \Z} \frac{1}{|n|^{2q} + 1} \E \Biggl( \int_{n-1/2}^{n+1/2} |F_{k}(1/2 - \frac{k}{\log x} + it)|^2 dt \Biggr)^q \nonumber \\
& \ll & \max_{N \in \Z} \frac{1}{(|N|+1)^{1/4}} \E \Biggl( \int_{N-1/2}^{N+1/2} |F_{k}(1/2 - \frac{k}{\log x} + it)|^2 dt \Biggr)^q . \nonumber
\end{eqnarray}
Proposition 2 follows on inserting this into the argument, and putting everything together and taking $2q$-th roots.
\qed

\subsection{Proof of Propositions 3 and 4}
Again we let $P(n)$ denote the largest prime factor of $n$, and we introduce an auxiliary Rademacher random variable $\epsilon$ (independent of everything else). Proceeding similarly as in section 2.2 of Harper, Nikeghbali and Radziwi\l\l~\cite{hnr}, we find that if $2/3 \leq q \leq 1$ then
\begin{eqnarray}
\E \Biggl|\sum_{\substack{n \leq x, \\ P(n) > \sqrt{x}}} f(n) \Biggr|^{2q} & = & \frac{1}{2^{2q}} \E \Biggl|\sum_{\substack{n \leq x, \\ P(n) > \sqrt{x}}} f(n) + \sum_{\substack{n \leq x, \\ P(n) \leq \sqrt{x}}} f(n) + \sum_{\substack{n \leq x, \\ P(n) > \sqrt{x}}} f(n) - \sum_{\substack{n \leq x, \\ P(n) \leq \sqrt{x}}} f(n) \Biggr|^{2q} \nonumber \\
& \leq & \E \Biggl|\sum_{\substack{n \leq x, \\ P(n) > \sqrt{x}}} f(n) + \sum_{\substack{n \leq x, \\ P(n) \leq \sqrt{x}}} f(n) \Biggr|^{2q} + \E \Biggl|\sum_{\substack{n \leq x, \\ P(n) > \sqrt{x}}} f(n) - \sum_{\substack{n \leq x, \\ P(n) \leq \sqrt{x}}} f(n) \Biggr|^{2q} \nonumber \\
& = & 2 \E \Biggl|\epsilon \sum_{\substack{n \leq x, \\ P(n) > \sqrt{x}}} f(n) + \sum_{\substack{n \leq x, \\ P(n) \leq \sqrt{x}}} f(n) \Biggr|^{2q} = 2 \E|\sum_{n \leq x} f(n)|^{2q} , \nonumber
\end{eqnarray}
since the law of $\epsilon \sum_{\substack{n \leq x, \\ P(n) > \sqrt{x}}} f(n) = \epsilon \sum_{\sqrt{x} < p \leq x} f(p) \sum_{m \leq x/p} f(m)$ conditional on the values $(f(p))_{p \leq \sqrt{x}}$ is the same as the law of $\sum_{\substack{n \leq x, \\ P(n) > \sqrt{x}}} f(n)$. We can rewrite this as
$$ \vert\vert \sum_{n \leq x} f(n) \vert\vert_{2q} \gg \vert\vert \sum_{\substack{n \leq x, \\ P(n) > \sqrt{x}}} f(n) \vert\vert_{2q} . $$

Now in the decomposition $\sum_{\substack{n \leq x, \\ P(n) > \sqrt{x}}} f(n) = \sum_{\sqrt{x} < p \leq x} f(p) \sum_{m \leq x/p} f(m)$, the inner sums are determined by the values $(f(p))_{p \leq \sqrt{x}}$, which are independent of the outer random variables $(f(p))_{\sqrt{x} < p \leq x}$. So conditioning on the values $(f(p))_{p \leq \sqrt{x}}$ determining the inner sums and applying the lower bound part of Khintchine's inequality (see e.g. Lemma 3.8.1 of Gut~\cite{gut} for the Rademacher case of this, the Steinhaus case may be proved similarly), it follows that
$$ \E|\sum_{\substack{n \leq x, \\ P(n) > \sqrt{x}}} f(n)|^{2q} \gg \E \left(\sum_{\sqrt{x} < p \leq x} \left|\sum_{m \leq x/p} f(m) \right|^2 \right)^{q} \geq \frac{1}{\log^{q}x} \E \left(\sum_{\sqrt{x} < p \leq x} \log p \left|\sum_{m \leq x/p} f(m) \right|^2 \right)^{q} . $$

Next there comes a smoothing step, similarly as in our work on upper bounds. Recall that we let $X = e^{\sqrt{\log x}}$. It is easy to check that we always have $|a + b|^2 \geq (1/4)|a|^2 - \min\{|b|^2, |a/2|^2 \} \geq 0$, say, and therefore
\begin{eqnarray}
&& \sum_{\sqrt{x} < p \leq x} \log p \Biggl|\sum_{m \leq x/p} f(m) \Biggr|^2 = \sum_{\sqrt{x} < p \leq x} \log p \frac{X}{p} \int_{p}^{p(1+1/X)} \Biggl|\sum_{m \leq x/p} f(m) \Biggr|^2 dt \nonumber \\
& \geq & \frac{1}{4} \sum_{\sqrt{x} < p \leq x} \log p \frac{X}{p} \int_{p}^{p(1+1/X)} \Biggl|\sum_{m \leq x/t} f(m) \Biggr|^2 dt \nonumber \\
&& - \sum_{\sqrt{x} < p \leq x} \log p \frac{X}{p} \int_{p}^{p(1+1/X)} \min\{\Biggl|\sum_{x/t < m \leq x/p} f(m) \Biggr|^2, \frac{1}{4} \Biggl|\sum_{m \leq x/t} f(m) \Biggr|^2 \} dt . \nonumber
\end{eqnarray}
Thus when $2/3 \leq q \leq 1$, we get
\begin{eqnarray}
\E|\sum_{\substack{n \leq x, \\ P(n) > \sqrt{x}}} f(n)|^{2q} & \gg & \frac{1}{\log^{q}x} \E \Biggl( \frac{1}{4} \sum_{\sqrt{x} < p \leq x} \log p \frac{X}{p} \int_{p}^{p(1+1/X)} \Biggl|\sum_{m \leq x/t} f(m) \Biggr|^2 dt \Biggr)^{q} \nonumber \\
&& - \frac{1}{\log^{q}x} \E\Biggl( \sum_{\sqrt{x} < p \leq x} \log p \frac{X}{p} \int_{p}^{p(1+1/X)} \Biggl|\sum_{x/t < m \leq x/p} f(m) \Biggr|^2 dt \Biggr)^{q} . \nonumber
\end{eqnarray}

To complete the proof of the first part of Proposition 3, we apply H\"{o}lder's inequality and obtain that $\E\Biggl( \sum_{\sqrt{x} < p \leq x} \log p \frac{X}{p} \int_{p}^{p(1+1/X)} \Biggl|\sum_{x/t < m \leq x/p} f(m) \Biggr|^2 dt \Biggr)^q$ is
\begin{eqnarray}
& \leq & \left(\sum_{\sqrt{x} < p \leq x} \log p \frac{X}{p} \int_{p}^{p(1+1/X)} \E \Biggl|\sum_{x/t < m \leq x/p} f(m) \Biggr|^2 dt \right)^q \nonumber \\
& \ll & \left( \sum_{\sqrt{x} < p \leq x} \log p (\frac{x}{p X} + 1) \right)^q \ll \left(\frac{x\log x}{X} + x \right)^q \ll x^q , \nonumber
\end{eqnarray}
and also that the main term $\E\Biggl( \sum_{\sqrt{x} < p \leq x} \log p \frac{X}{p} \int_{p}^{p(1+1/X)} \left|\sum_{m \leq x/t} f(m) \right|^2 dt \Biggr)^q$ is
\begin{eqnarray}
& = & \E\left( \int_{\sqrt{x}}^{x} \sum_{t/(1 + 1/X) < p \leq t} \log p \frac{X}{p} \left|\sum_{m \leq x/t} f(m) \right|^2 dt \right)^q \nonumber \\
& \gg & \E\left( \int_{\sqrt{x}}^{x} \left|\sum_{m \leq x/t} f(m) \right|^2 dt \right)^q = x^{q} \E\left( \int_{1}^{\sqrt{x}} \left|\sum_{m \leq z} f(m) \right|^2 \frac{dz}{z^{2}} \right)^q . \nonumber
\end{eqnarray}
The first part of Proposition 3 follows on putting all our calculations together.

To deduce the second part of Proposition 3, we simply note that for any large $V$ and any $2/3 \leq q \leq 1$ we have
\begin{eqnarray}
&& \E\left( \int_{1}^{\sqrt{x}} \left|\sum_{m \leq z} f(m) \right|^2 \frac{dz}{z^{2}} \right)^q \geq \E\left( \int_{1}^{\sqrt{x}} \left|\sum_{\substack{m \leq z, \\ x \text{-smooth}}} f(m) \right|^2 \frac{dz}{z^{2 + 8V/\log x}} \right)^q \nonumber \\
& \geq & \E\left( \int_{1}^{\infty} \left|\sum_{\substack{m \leq z, \\ x \text{-smooth}}} f(m) \right|^2 \frac{dz}{z^{2 + 8V/\log x}} \right)^q - \E\left( \int_{\sqrt{x}}^{\infty} \left|\sum_{\substack{m \leq z, \\ x \text{-smooth}}} f(m) \right|^2 \frac{dz}{z^{2 + 8V/\log x}} \right)^q \nonumber \\
& \geq & \E\left( \int_{1}^{\infty} \left|\sum_{\substack{m \leq z, \\ x \text{-smooth}}} f(m) \right|^2 \frac{dz}{z^{2 + 8V/\log x}} \right)^q - \frac{1}{e^{2Vq}} \E\left( \int_{1}^{\infty} \left|\sum_{\substack{m \leq z, \\ x \text{-smooth}}} f(m) \right|^2 \frac{dz}{z^{2 + 4V/\log x}} \right)^q . \nonumber
\end{eqnarray}
By Harmonic Analysis Result 1, the first term here is $\gg \E(\int_{-1/2}^{1/2} |F(1/2 + \frac{4V}{\log x} + it)|^2 dt )^q$ and the subtracted second term is $\ll e^{-2Vq} \E(\int_{-\infty}^{\infty} \frac{|F(1/2 + \frac{2V}{\log x} + it)|^2}{|1/2 + \frac{2V}{\log x} + it|^2} dt )^q$, which in the Steinhaus case is $\ll e^{-2Vq} \E(\int_{-1/2}^{1/2} |F(1/2 + \frac{2V}{\log x} + it)|^2 dt )^q $ by ``translation invariance in law''. Putting everything together, this finishes the proof of Proposition 3.
\qed

\vspace{12pt}
The arguments in the Rademacher case are, once again, exactly the same until the final line, where we don't have ``translation invariance'' so we must upper bound $\E(\int_{-\infty}^{\infty} \frac{|F(1/2 + \frac{2V}{\log x} + it)|^2}{|1/2 + \frac{2V}{\log x} + it|^2} dt )^q$ by $\max_{N \in \Z} \frac{1}{(|N|+1)^{1/4}} \E ( \int_{N-1/2}^{N+1/2} |F(1/2 + \frac{2V}{\log x} + it)|^2 dt )^q$.
\qed

\section{Probabilistic calculations}\label{secprobcalc}
In this section we collect together various probabilistic calculations, that we wish to isolate in advance of our main proofs. These are of two basic kinds: firstly the estimation of the mean square, and related quantities, for some random Euler products; and secondly, estimation of the probability of certain events where the measure is weighted by the mean square of random Euler products. The latter calculations provide analogues of Girsanov's theorem in our setting, where the logarithms of our random products are not actually Gaussian random variables, but only approximately so.

\subsection{The mean square of random Euler products}

\begin{lem1}
If $f$ is a Steinhaus random multiplicative function, then for any real $t$ and $u$, any real $400(1 + u^2) \leq x \leq y$, and any real $\sigma \geq - 1/\log y$, we have
$$ \E \prod_{x < p \leq y} \left|1 - \frac{f(p)}{p^{1/2+\sigma}}\right|^{-2} \left|1 - \frac{f(p)}{p^{1/2+\sigma+it}}\right|^{-iu} = \exp\{\sum_{x < p \leq y} \frac{1 + iu\cos(t\log p) - u^{2}/4}{p^{1 + 2\sigma}} + T(u)\} , $$
where $T(u) = T_{x,y,\sigma,t}(u)$ satisfies $|T(u)| \ll \frac{1 + |u|^3}{\sqrt{x} \log x}$, and its derivative satisfies $|T'(u)| \ll \frac{1}{\sqrt{x} \log x}$ when $|u| \leq 1$.
\end{lem1}

\begin{proof}[Proof of Lemma 1]
To simplify the writing of the proof, we temporarily set $R_{p}(t) := -\Re\log(1 - \frac{f(p)}{p^{1/2+\sigma+it}})$. Thus we may rewrite
$$ \left|1 - \frac{f(p)}{p^{1/2+\sigma}}\right|^{-2} \left|1 - \frac{f(p)}{p^{1/2+\sigma+it}}\right|^{-iu} = \exp\{2 R_{p}(0) +iu R_{p}(t) \} = 1 + \sum_{j=1}^{\infty} \frac{(2 R_{p}(0) +iu R_{p}(t))^j}{j!} . $$

Now using the Taylor expansion of the logarithm, we have $R_{p}(t) = \sum_{k=1}^{\infty} \frac{\Re(f(p)p^{-it})^k}{k p^{k(1/2+\sigma)}} = \frac{\Re f(p) p^{-it}}{p^{1/2+\sigma}} + O(\frac{1}{p^{1+2\sigma}})$. In particular, by symmetry we have $\E \Re (f(p)p^{-it})^k = \E \Re f(p)^k = 0$ for all $k \geq 1$, and therefore $\E R_{p}(t) = 0$. We also have
$$ \E R_{p}(t)^2 = \E\frac{(\Re f(p)p^{-it})^2}{p^{1 + 2\sigma}} + O(\frac{1}{p^{3/2+3\sigma}}) = \frac{1}{2 p^{1 + 2\sigma}} + O(\frac{1}{p^{3/2+3\sigma}}) , $$
as well as
$$ \E R_{p}(0) R_{p}(t) = \E\frac{(\Re f(p))(\Re f(p)p^{-it})}{p^{1 + 2\sigma}} + O(\frac{1}{p^{3/2+3\sigma}}) = \frac{\cos(t\log p)}{2 p^{1 + 2\sigma}} + O(\frac{1}{p^{3/2+3\sigma}}) . $$
For $j \geq 3$ we can use the trivial bound $|R_{p}(t)^j| \leq (\sum_{k=1}^{\infty} \frac{1}{p^{k(1/2+\sigma)}})^j = \frac{1}{(p^{1/2+\sigma}-1)^j}$.

Next, we note that for primes $y \geq p > x \geq 400(1+u^2)$ we have $\frac{1}{p^{1/2+\sigma}} = \frac{e^{-\sigma\log p}}{p^{1/2}} \leq \frac{e}{p^{1/2}}$, and therefore $(2+|u|)/p^{1/2 + \sigma} \leq e/5$. So putting things together, for such primes we have
\begin{eqnarray}
&& \E\left|1 - \frac{f(p)}{p^{1/2+\sigma}}\right|^{-2} \left|1 - \frac{f(p)}{p^{1/2+\sigma+it}}\right|^{-iu} \nonumber \\
& = & 1 + \frac{(4\E R_{p}(0)^2 + 4iu\E R_{p}(0)R_{p}(t) - u^{2} \E R_{p}(t)^2)}{2} + \E \sum_{j=3}^{\infty} \frac{(2 R_{p}(0) +iu R_{p}(t))^j}{j!} \nonumber \\
& = & 1 + \frac{(1+iu\cos(t\log p) -u^{2}/4)}{p^{1+2\sigma}} + O(\sum_{j=3}^{\infty} \frac{(2+|u|)^j}{j! (p^{1/2+\sigma}-1)^j} ) \nonumber \\
& = & 1 + \frac{(1+iu\cos(t\log p) -u^{2}/4)}{p^{1+2\sigma}} + D_{p}(u) , \nonumber
\end{eqnarray}
where $D_{p}(u)$ satisfies $|D_{p}(u)| \ll \frac{1+|u|^3}{p^{3/2+3\sigma}} \ll \frac{1+|u|^3}{p^{3/2}}$, and its derivative satisfies $|D_{p}'(u)| \ll \frac{1}{p^{3/2+3\sigma}} \ll \frac{1}{p^{3/2}}$ for $|u| \leq 1$.

Finally, note that we may rewrite
\begin{eqnarray}
\E\left|1 - \frac{f(p)}{p^{1/2+\sigma}}\right|^{-2} \left|1 - \frac{f(p)}{p^{1/2+\sigma+it}}\right|^{-iu} & = & 1 + \frac{(1+iu\cos(t\log p) -u^{2}/4)}{p^{1+2\sigma}} + D_{p}(u) \nonumber \\
& = & \exp\{\frac{(1+iu\cos(t\log p) -u^{2}/4)}{p^{1+2\sigma}} + T_{p}(u)\} , \nonumber
\end{eqnarray}
where $T_{p}(u)$ again satisfies $|T_{p}(u)| \ll \frac{1+|u|^3}{p^{3/2}}$, and also $|T_{p}'(u)| \ll \frac{1}{p^{3/2}}$ for $|u| \leq 1$. Since $f$ is independent on distinct primes, we then deduce
$$ \E \prod_{x < p \leq y} \left|1 - \frac{f(p)}{p^{1/2+\sigma}}\right|^{-2} \left|1 - \frac{f(p)}{p^{1/2+\sigma+it}}\right|^{-iu}  = \exp\{\sum_{x < p \leq y} \frac{1 + iu\cos(t\log p) - u^{2}/4}{p^{1 + 2\sigma}} + \sum_{x < p \leq y} T_{p}(u)\} , $$
which implies Lemma 1 in view of the standard Chebychev-type estimate $\sum_{p > x} 1/p^{3/2} \ll 1/(\sqrt{x} \log x)$.
\end{proof}

We note in particular that, by Lemma 1, for any $400 \leq x \leq y$ and any $\sigma \geq -1/\log y$ and for Steinhaus random multiplicative $f$ we have
\begin{equation}\label{stprodusual}
\E \prod_{x < p \leq y} \left|1 - \frac{f(p)}{p^{1/2+\sigma}}\right|^{-2} = \exp\{\sum_{x < p \leq y} \frac{1}{p^{1 + 2\sigma}} + O(\frac{1}{\sqrt{x}\log x}) \} .
\end{equation}
Moreover, under the same conditions the same proof shows that $\E \prod_{x < p \leq y} \left|1 - \frac{f(p)}{p^{1/2+\sigma}}\right|^{2}$ satisfies this estimate too (note the change of the exponent from $-2$ to 2).

\vspace{12pt}
We will need an analogue of Lemma 1 for the Rademacher case. Both the formulation and proof of this are slightly more complicated. This is because the distribution of a Rademacher random multiplicative function is not invariant under shifts by $n^{it}$, and also because $f(p)^2 \equiv 1$ in the Rademacher case, and the Euler products we deal with here have slightly different forms.

\begin{lem2}
If $f$ is a Rademacher random multiplicative function, then for any real $t_1 , t_2$ and $u$, any real $400(1 + u^2) \leq x \leq y$, and any real $\sigma \geq - 1/\log y$, we have
\begin{eqnarray}
&& \E \prod_{x < p \leq y} \left|1 + \frac{f(p)}{p^{1/2+\sigma+it_1}}\right|^{2} \left|1 + \frac{f(p)}{p^{1/2+\sigma+i(t_1 + t_2)}}\right|^{iu} \nonumber \\
& = & \exp\{\sum_{x < p \leq y} \frac{1 + iuc(t_1 , t_2 , p) - (u^{2}/4)(1+\cos(2(t_1 + t_2)\log p))}{p^{1 + 2\sigma}} + T(u)\} . \nonumber
\end{eqnarray}
Here $c(t_1 , t_2 , p) = 2\cos(t_1 \log p) \cos((t_1 + t_2)\log p) - (1/2)\cos(2(t_1+t_2)\log p)$, and $T(u) = T_{x,y,\sigma,t_1,t_2}(u)$ satisfies $|T(u)| \ll \frac{1 + |u|^3}{\sqrt{x} \log x}$, and its derivative satisfies $|T'(u)| \ll \frac{1}{\sqrt{x} \log x}$ when $|u| \leq 1$. 
\end{lem2}

\begin{proof}[Proof of Lemma 2]
Let us temporarily set $R_{p}(t) := \Re\log(1 + \frac{f(p)}{p^{1/2+\sigma+it}})$, so we may rewrite
\begin{eqnarray}
\left|1 + \frac{f(p)}{p^{1/2+\sigma+it_1}}\right|^{2} \left|1 + \frac{f(p)}{p^{1/2+\sigma+i(t_1 + t_2)}}\right|^{iu} & = & \exp\{2 R_{p}(t_1) +iu R_{p}(t_1 + t_2) \} \nonumber \\
& = & 1 + \sum_{j=1}^{\infty} \frac{(2 R_{p}(t_1) +iu R_{p}(t_1 + t_2))^j}{j!} . \nonumber
\end{eqnarray}

Using Taylor expansion (and the fact that $f(p) \in \{\pm 1\}$), we obtain that $R_{p}(t) = \sum_{k=1}^{\infty} (-1)^{k-1} \frac{\Re(f(p)p^{-it})^k}{k p^{k(1/2+\sigma)}} = \frac{f(p) \cos(t\log p)}{p^{1/2+\sigma}} - \frac{\cos(2t\log p)}{2p^{1+2\sigma}} + O(\frac{1}{p^{3/2+3\sigma}})$. In particular, this implies that $\E R_{p}(t) = - \frac{\cos(2t\log p)}{2p^{1+2\sigma}} + O(\frac{1}{p^{3/2+3\sigma}})$. We also have
$$ \E R_{p}(t)^2 = \frac{\cos^{2}(t\log p)}{p^{1 + 2\sigma}} + O(\frac{1}{p^{3/2+3\sigma}}) = \frac{(1+\cos(2t\log p))}{2 p^{1 + 2\sigma}} + O(\frac{1}{p^{3/2+3\sigma}}) , $$
as well as
$$ \E R_{p}(t_1) R_{p}(t_1 + t_2) = \frac{\cos(t_1 \log p) \cos((t_1 + t_2)\log p)}{p^{1 + 2\sigma}} + O(\frac{1}{p^{3/2+3\sigma}}) . $$
For $j \geq 3$ we again have $|R_{p}(t)^j| \leq (\sum_{k=1}^{\infty} \frac{1}{p^{k(1/2+\sigma)}})^j = \frac{1}{(p^{1/2+\sigma}-1)^j}$.

As in Lemma 1, for primes $y \geq p > x \geq 400(1+u^2)$ we have $(2+|u|)/p^{1/2 + \sigma} \leq e/5$, so for such primes we get
\begin{eqnarray}
&& \E\left|1 + \frac{f(p)}{p^{1/2+\sigma+it_1}}\right|^{2} \left|1 + \frac{f(p)}{p^{1/2+\sigma+i(t_1 + t_2)}}\right|^{iu} \nonumber \\
& = & 1 - \frac{\cos(2t_1 \log p)}{p^{1+2\sigma}} - \frac{(iu/2)\cos(2(t_1 + t_2)\log p)}{p^{1+2\sigma}} + O(\frac{1 + |u|}{p^{3/2+3\sigma}}) + \nonumber \\
&& + \frac{(1+\cos(2t_1 \log p)) + 2iu \cos(t_1 \log p) \cos((t_1 + t_2)\log p) - \frac{u^{2}}{4}(1+\cos(2(t_1 + t_2)\log p))}{p^{1+2\sigma}} \nonumber \\
&& + O(\frac{1 + u^2}{p^{3/2+3\sigma}}) + \E\sum_{j=3}^{\infty} \frac{(2 R_{p}(t_1) +iu R_{p}(t_1 + t_2))^j}{j!} \nonumber \\
& = & 1 + \frac{1+iuc(t_1 , t_2 , p) - (u^{2}/4)(1+\cos(2(t_1 + t_2)\log p))}{p^{1+2\sigma}} + D_{p}(u) , \nonumber
\end{eqnarray}
where $D_{p}(u)$ satisfies $|D_{p}(u)| \ll \frac{1+|u|^3}{p^{3/2+3\sigma}} \ll \frac{1+|u|^3}{p^{3/2}}$, and its derivative satisfies $|D_{p}'(u)| \ll \frac{1}{p^{3/2+3\sigma}} \ll \frac{1}{p^{3/2}}$ for $|u| \leq 1$.

From this point, the proof concludes exactly as the proof of Lemma 1.
\end{proof}

Again, Lemma 2 implies that for any $400 \leq x \leq y$ and $\sigma \geq -1/\log y$ and real $t_1$; and for Rademacher random multiplicative $f$; we have
\begin{equation}\label{radprodusual}
\E \prod_{x < p \leq y} \left|1 + \frac{f(p)}{p^{1/2+\sigma+it_1}}\right|^{2} = \exp\{\sum_{x < p \leq y} \frac{1}{p^{1 + 2\sigma}} + O(\frac{1}{\sqrt{x}\log x}) \} .
\end{equation}
Under these conditions, we can use the same proof to estimate $\E \prod_{x < p \leq y} \left|1 + \frac{f(p)}{p^{1/2+\sigma+it_1}}\right|^{-2}$ (note the change of the exponent from 2 to $-2$), and in the Rademacher case this slightly changes the answer: we find
\begin{equation}\label{radprodinv}
\E \prod_{x < p \leq y} \left|1 + \frac{f(p)}{p^{1/2+\sigma+it_1}}\right|^{-2} = \exp\{\sum_{x < p \leq y} \frac{1 + 2\cos(2 t_1 \log p)}{p^{1 + 2\sigma}} + O(\frac{1}{\sqrt{x}\log x}) \} .
\end{equation}

\subsection{Girsanov-type calculations in the Steinhaus case}\label{subsecstgirsanov}
Let $f(n)$ be a Steinhaus random multiplicative function, and let $x$ be large and $-1/100 \leq \sigma \leq 1/100$, say. Later we will impose some further restrictions on $\sigma$. Let us introduce a new ``tilted'' probability measure $\tilde{\p} = \tilde{\p}_{x,\sigma}$ by setting
$$ \tilde{\p}(A) := \frac{\E \textbf{1}_{A} \prod_{p \leq x^{1/e}} \left|1 - \frac{f(p)}{p^{1/2+\sigma}}\right|^{-2}}{\E \prod_{p \leq x^{1/e}} \left|1 - \frac{f(p)}{p^{1/2+\sigma}}\right|^{-2}} $$
for each event $A$, where $\textbf{1}$ denotes the indicator function. We will also sometimes write $\tilde{\E}$ to denote expectation (i.e. integration) with respect to the measure $\tilde{\p}$. The exact choice of the range of $p$ in the definition of $\tilde{\p}$ is not too important, since the independence of the $f(p)$ means that if the event $A$ doesn't involve a particular prime, the expectation of that term will factor out and cancel between the numerator and denominator.

Furthermore, for each $l \in \N \cup \{0\}$ set $I_l(s) := \prod_{x^{e^{-(l+2)}} < p \leq x^{e^{-(l+1)}}} (1 - \frac{f(p)}{p^s})^{-1}$, the $l$-th ``increment'' of the Euler product corresponding to Steinhaus $f$.

Our goal here is to estimate $\tilde{\p}(A)$ for certain events $A$ corresponding to restrictions on partial Euler products. We wish to show these probabilities are essentially the same as they would be if the increments $\log|I_{l}(s)|$ were Gaussian. In that case the very useful Girsanov theorem shows the tilted probabilities are themselves Gaussian probabilities for certain shifted Gaussians, so can be understood and estimated fairly straightforwardly. We will achieve all this by combining Lemma 1 (which corresponds to a characteristic function calculation for the tilted measure $\tilde{\p}$) with the Berry--Esseen theorem on distributional approximation, and then develop the necessary Gaussian estimates. Since we must apply the Berry--Esseen theorem multiple times, and it only supplies an absolute rather than a relative error, we must be careful to control the sizes of things to make this all work.

\vspace{12pt}
Let $(l_j)_{j=1}^{n}$ denote a {\em strictly decreasing} sequence of non-negative integers, with $l_1 \leq \log\log x - 2$, and define a corresponding increasing sequence of real numbers $(x_j)_{j=1}^{n}$ by setting $x_j := x^{e^{-(l_j + 1)}}$. 

\begin{lem3}
Let the situation be as above, and suppose that $x_1$ is sufficiently large and that $|\sigma| \leq 1/\log x_n$. Suppose further that $(v_j)_{j=1}^{n}$ is any sequence of real numbers satisfying
 $$ |v_j| \leq (1/40)\sqrt{\log x_j} + 2 \;\;\; \forall 1 \leq j \leq n , $$
and $(t_j)_{j=1}^{n}$ is any sequence of real numbers.

Then we have
\begin{eqnarray}
&& \tilde{\p}(v_j \leq \log|I_{l_j}(1/2 + \sigma + it_j)| \leq v_j + 1/j^2 \; \forall 1 \leq j \leq n) \nonumber \\
& = & \left(1+O\left(\frac{1}{x_{1}^{1/100}}\right) \right) \p(v_j \leq N_j \leq v_j + 1/j^2 \; \forall 1 \leq j \leq n) , \nonumber
\end{eqnarray}
where $N_j$ are independent Gaussian random variables with mean $\sum_{x_{j}^{1/e} < p \leq x_j} \frac{\cos(t_j \log p)}{p^{1 + 2\sigma}}$ and variance $\sum_{x_{j}^{1/e} < p \leq x_j} \frac{1}{2p^{1 + 2\sigma}}$.
\end{lem3}

\begin{proof}[Proof of Lemma 3]
By independence, both the probability on the left and the one on the right factor as a product over $j$, so it will suffice to prove that
$$ \tilde{\p}(v_j \leq \log|I_{l_j}(1/2 + \sigma + it_j)| \leq v_j + 1/j^2) = (1+O(\frac{1}{x_{j}^{1/100}}))\p(v_j \leq N_j \leq v_j + 1/j^2) \;\;\; \forall j \leq n . $$
We also note at the outset that, because of our assumption on $\sigma$ as well as the standard Mertens estimate for sums over primes (see e.g. Theorem 2.7 of Montgomery and Vaughan~\cite{mv}), we have
$$ e^{-2}(1 + O(\frac{1}{\log x_j})) \leq \sum_{x_{j}^{1/e} < p \leq x_j} \frac{1}{p^{1 + 2\sigma}} \leq e^{2} \sum_{x_{j}^{1/e} < p \leq x_j} \frac{1}{p} = e^{2}(1 + O(\frac{1}{\log x_j})) \;\;\; \forall 1 \leq j \leq n . $$

It is a standard calculation that the characteristic function
\begin{eqnarray}
\E e^{iuN_j} & = & \E \exp\{iu\sum_{x_{j}^{1/e} < p \leq x_j} \frac{\cos(t_j \log p)}{p^{1 + 2\sigma}} + iu\sqrt{\sum_{x_{j}^{1/e} < p \leq x_j} \frac{1}{2p^{1 + 2\sigma}}} N(0,1)\} \nonumber \\
& = & \exp\{\sum_{x_{j}^{1/e} < p \leq x_j} \frac{iu\cos(t_j \log p) - u^{2}/4}{p^{1 + 2\sigma}} \} . \nonumber
\end{eqnarray}
Combining this with Lemma 1 (applied with $x,y$ replaced by $x_{j}^{1/e}$ and $x_j$), we see that for any $|u| \leq x_{j}^{1/20}$ we have
\begin{eqnarray}
|\tilde{\E} e^{iu\log|I_{l_j}(1/2 + \sigma + it_j)|} - \E e^{iuN_j}| & = & \exp\{ - \frac{u^2}{4} \sum_{x_{j}^{1/e} < p \leq x_j} \frac{1}{p^{1 + 2\sigma}} \} |e^{T(u) - T(0)} - 1| \nonumber \\
& \ll & \exp\{ - \frac{u^2}{4} \sum_{x_{j}^{1/e} < p \leq x_{j}} \frac{1}{p^{1 + 2\sigma}} \} \frac{|u| + |u|^3}{\sqrt{x_{j}^{1/e}} \log x_j} . \nonumber
\end{eqnarray}
Then by the Berry--Esseen theorem (see e.g. Lemma 7.6.1 of Gut~\cite{gut}), we have
\begin{eqnarray}
&& |\tilde{\p}(v_j \leq \log|I_{l_j}(1/2 + \sigma + it_j)| \leq v_j + 1/j^2) - \p(v_j \leq N_j \leq v_j + 1/j^2)| \nonumber \\
& \ll & \int_{- x_{j}^{1/20}}^{x_{j}^{1/20}} \Biggl|\frac{\tilde{\E} e^{iu\log|I_{l_j}(1/2 + \sigma + it_j)|} - \E e^{iuN_j}}{u} \Biggr| du + \frac{1}{x_{j}^{1/20}} \nonumber \\
& \ll & \int_{- x_{j}^{1/20}}^{x_{j}^{1/20}} \frac{1+u^2}{\sqrt{x_{j}^{1/e}} \log x_j} du + \frac{1}{x_{j}^{1/20}} \ll \frac{1}{x_{j}^{1/60}} . \nonumber
\end{eqnarray}

Finally, since $N_j$ has mean $O(1)$ and variance $\frac{1}{2} \sum_{x_{j}^{1/e} < p \leq x_j} \frac{1}{p^{1 + 2\sigma}} \geq \frac{e^{-2} + O(1/\log x_j)}{2}$, and since $|v_j| \leq (1/40)\sqrt{\log x_j} + 2$, we have
$$ \p(v_j \leq N_j \leq v_j + 1/j^2) \gg (1/j^2) e^{-(|v_j| + O(1))^{2}/(e^{-2} + O(1/\log x_j))} \gg \frac{1}{j^{2} x_{j}^{9/1600}} \gg \frac{1}{x_{j}^{1/150}} . $$
So we may rewrite the error $1/x_{j}^{1/60}$ from the Berry--Esseen Theorem as $O(\frac{\p(v_j \leq N_j \leq v_j + 1/j^2)}{x_{j}^{1/60 - 1/150}})$, which is $O(\frac{\p(v_j \leq N_j \leq v_j + 1/j^2)}{x_{j}^{1/100}})$.
\end{proof}

Using Lemma 3 and a suitable slicing argument, we can obtain a similar approximation result for the probabilities of slightly more complicated events.

\begin{lem4}
Suppose $(u_j)_{j=1}^{n}$ and $(v_j)_{j=1}^{n}$ are sequences of real numbers satisfying
$$ -(1/80)\sqrt{\log x_j} \leq u_j \leq v_j \leq (1/80)\sqrt{\log x_j} \;\;\; \forall 1 \leq j \leq n , $$
and otherwise let the situation be as in Lemma 3. Then we have
\begin{eqnarray}
\p(u_j + 2 \leq \sum_{m=1}^{j} N_m \leq v_j - 2 \; \forall j \leq n) & \ll & \tilde{\p}(u_j \leq \sum_{m=1}^{j} \log|I_{l_m}(\frac{1}
{2} + \sigma + it_m)| \leq v_j \; \forall j \leq n) \nonumber \\
& \ll & \p(u_j - 2 \leq \sum_{m=1}^{j} N_m \leq v_j + 2 \; \forall j \leq n) , \nonumber
\end{eqnarray}
where the Gaussian random variables $N_m$ are also as in Lemma 3.

In addition, if the numbers $(t_j)_{j=1}^{n}$ satisfy $|t_j| \leq \frac{1}{j^{2/3} \log x_j}$ then we have
\begin{eqnarray}
&& \p((u_j - j) + O(1) \leq \sum_{m=1}^{j} G_m \leq (v_j-j) - O(1) \; \forall j \leq n) \nonumber \\
& \ll & \tilde{\p}(u_j \leq \sum_{m=1}^{j} \log|I_{l_m}(1/2 + \sigma + it_m)| \leq v_j \; \forall j \leq n) \nonumber \\
& \ll & \p((u_j - j) - O(1) \leq \sum_{m=1}^{j} G_m \leq (v_j-j) + O(1) \; \forall j \leq n) , \nonumber
\end{eqnarray}
where $G_m$ are independent Gaussian random variables, each having mean 0 and variance $\sum_{x_{m}^{1/e} < p \leq x_m} \frac{1}{2p^{1 + 2\sigma}}$.
\end{lem4}

In the second part of Lemma 4, the assumption that $|t_m| \leq \frac{1}{m^{2/3} \log x_m}$ together with our standing assumption $|\sigma| \leq \frac{1}{\log x_n}$ implies that $\E N_m = \sum_{x_{m}^{1/e} < p \leq x_m} \frac{\cos(t_m \log p)}{p^{1 + 2\sigma}} \approx \sum_{x_{m}^{1/e} < p \leq x_m} \frac{1}{p} \approx 1$. Subtracting these means for each $m$ yields the mean zero random variables $G_m$, and produces the subtracted term $-j$ in the upper and lower bounds of our events. As we shall see, this subtracted term in the upper bound will be at the heart of everything, since it can reduce a sequence of fairly large terms $v_j \approx j$ that one would expect to impose very little constraint on our partial sums, to a smaller sequence that does impose a non-trivial constraint.

\begin{proof}[Proof of Lemma 4]
As remarked previously, to prove the first part (involving the $N_m$) we shall simply approximate the probability we are interested in by a sum of probabilities of the form treated in Lemma 3.

Note that if we have
$$ -(1/80)\sqrt{\log x_j} \leq \sum_{m=1}^{j} \log|I_{l_m}(1/2 + \sigma + it_m)| \leq (1/80)\sqrt{\log x_j} , $$
and the analogous bounds for the sum up to $j-1$, then we must have
$$ |\log|I_{l_j}(1/2 + \sigma + it_j)|| \leq (1/80)\sqrt{\log x_{j-1}} + (1/80)\sqrt{\log x_{j}} \leq (1/40)\sqrt{\log x_{j}} . $$
So if we let $\mathcal{R}_1 := \{r \in \Z : |r| \leq (1/40)\sqrt{\log x_1} + 2\}$ (say), and more generally let $\mathcal{R}_j := \{r \in (1/j^2)\Z : |r| \leq (1/40)\sqrt{\log x_j} + 2\}$, then in order to have $u_j \leq \sum_{m=1}^{j} \log|I_{l_m}(1/2 + \sigma + it_m)| \leq v_j$ for all $1 \leq j \leq n$ we must have
$$ r_j \leq \log|I_{l_j}(1/2 + \sigma + it_j)| \leq r_j + 1/j^2 \;\;\; \forall 1 \leq j \leq n , $$
for some $r_1 \in \mathcal{R}_1 , ..., r_n \in \mathcal{R}_n$ satisfying $u_j - \sum_{m=1}^{j} 1/m^2 \leq \sum_{m=1}^{j} r_m \leq v_j$ for all $1 \leq j \leq n$. Thus, applying Lemma 3, we obtain
\begin{eqnarray}
&& \tilde{\p}(u_j \leq \sum_{m=1}^{j} \log|I_{l_m}(1/2 + \sigma + it_m)| \leq v_j \; \forall 1 \leq j \leq n) \nonumber \\
& \ll & \sum_{\substack{r_1 \in \mathcal{R}_1 , ... , r_n \in \mathcal{R}_n , \\ u_j - \sum_{m=1}^{j} 1/m^2 \leq \sum_{m=1}^{j} r_m \leq v_j \; \forall 1 \leq j \leq n}} \p(r_j \leq N_j \leq r_j + 1/j^2 \; \forall 1 \leq j \leq n) \nonumber \\
& \leq & \p(u_j - \sum_{m=1}^{j} \frac{1}{m^2} \leq \sum_{m=1}^{j} N_m \leq v_j + \sum_{m=1}^{j} \frac{1}{m^2} \; \forall j) \leq \p(u_j - 2 \leq \sum_{m=1}^{j} N_m \leq v_j + 2 \; \forall j \leq n) . \nonumber
\end{eqnarray}
This is the desired upper bound, and the proof of the corresponding lower bound is exactly similar.

The second part of the lemma follows by setting $G_m = N_m - \E N_m = N_m - \sum_{x_{m}^{1/e} < p \leq x_m} \frac{\cos(t_m \log p)}{p^{1 + 2\sigma}}$, and noting that
\begin{eqnarray}
\sum_{x_{m}^{1/e} < p \leq x_m} \frac{\cos(t_m \log p)}{p^{1 + 2\sigma}} & = & \sum_{x_{m}^{1/e} < p \leq x_m} \frac{1}{p^{1 + 2\sigma}} + O(\sum_{x_{m}^{1/e} < p \leq x_m} \frac{(|t_m|\log p)^2}{p^{1 + 2\sigma}}) \nonumber \\
& = & \sum_{x_{m}^{1/e} < p \leq x_m} \frac{1}{p} + O(\sum_{x_{m}^{1/e} < p \leq x_m} \frac{|\sigma| \log p + (|t_m|\log p)^2}{p}) \nonumber \\
& = & 1 + O(\frac{1}{\log x_m} + \frac{\log x_m}{\log x_n} + \frac{1}{m^{4/3}}) = 1 + O(\frac{1}{e^{n-m}} + \frac{1}{m^{4/3}}) , \nonumber
\end{eqnarray}
under our conditions that $|\sigma| \leq 1/\log x_n$ and $|t_m| \leq \frac{1}{m^{2/3}\log x_m}$. Here we used the standard Mertens and Chebychev estimates, namely $\sum_{x_{m}^{1/e} < p \leq x_m} \frac{1}{p} = 1 + O(\frac{1}{\log x_m})$ and $\sum_{x_{m}^{1/e} < p \leq x_m} \frac{\log p}{p} \ll \log x_m$. We also used the facts that $x_m \geq x_{1}^{e^{m-1}}$ and, more generally, $x_n \geq x_{m}^{e^{n-m}}$, both of which follow from the definition of the sequence $(x_j)_{j=1}^{n}$.
\end{proof}

\vspace{12pt}
At this point it will be helpful to examine the kinds of events we shall actually be interested in for our application, and their probabilities in the Gaussian case.

\begin{probres1}
Let $a$ and $n$ be large, and let $G_1 , ..., G_n$ be independent Gaussian random variables, each having mean zero and variance between $1/20$ and $20$ (say). Then uniformly for any function $h(j)$ satisfying $|h(j)| \leq 10\log j$, we have
$$ \p(\sum_{m=1}^{j} G_m \leq a + h(j) \; \forall 1 \leq j \leq n) \asymp \min\{1, \frac{a}{\sqrt{n}}\} . $$
\end{probres1}

Since Probability Result 1 is a purely probabilistic statement about Gaussian random walks, we postpone its proof to the appendix. We remark that such a statement with $h(j) \equiv 0$ is standard. It requires a little more work to obtain the more general estimate, but it is natural to think such a result should hold because one typically thinks of random variables fluctuating on the order of their standard deviation (here $\asymp \sqrt{j}$ up to step $j$), so perturbing an event by a term much smaller than this shouldn't alter its probability very much. Having a more general estimate will be very useful later because, on taking exponentials, the function $h(j)$ can supply extra savings of powers of $j$ in various places.

We will also need a small variant of Probability Result 1, where the event involves a non-trivial but relaxed lower barrier and a slightly tightened upper barrier as well.
\begin{probres2}
There is a large absolute constant $B$ such that the following is true. Let $a$ and $n$ be large, and let $G_1 , ..., G_n$ be independent Gaussian random variables, each having mean zero and variance between $1/20$ and $20$. Then uniformly for any function $h(j)$ satisfying $|h(j)| \leq 10\log j$, and any function $g(j)$ satisfying $g(j) \leq -Bj$, we have
$$ \p(g(j) \leq \sum_{m=1}^{j} G_m \leq \min\{a,Bj\} + h(j) \; \forall 1 \leq j \leq n) \asymp \min\{1, \frac{a}{\sqrt{n}}\} . $$
\end{probres2}

Again, it is natural to think such a result should hold because, if the sum up to $j$ typically fluctuates on the scale of $\sqrt{j}$, then the sums should be rather insensitive to a much more relaxed barrier of the shape $Bj$ or $-Bj$. We postpone the proof to the appendix.

\vspace{12pt}
By combining Lemma 4 with Probability Results 1 and 2, we can finally prove the $\tilde{\p}$ probability estimate we shall need to obtain Theorem 1. Again, this requires a little care to ensure the size restrictions on $u_j, v_j$ in Lemma 4 are respected when it is applied.

\begin{prop5}
There is a large natural number $B$ such that the following is true.

Let $n \leq \log\log x - (B+1)$ be large, and define the decreasing sequence $(l_j)_{j=1}^{n}$ of non-negative integers by $l_j := \lfloor \log\log x \rfloor - (B+1) - j$. Suppose that $|\sigma| \leq \frac{1}{e^{B+n+1}}$, and that $(t_j)_{j=1}^{n}$ is a sequence of real numbers satisfying $|t_j| \leq \frac{1}{j^{2/3} e^{B+j+1}}$ for all $j$.

Then uniformly for any large $a$ and any function $h(n)$ satisfying $|h(n)| \leq 10\log n$, and with $I_{l}(s)$ denoting the increments of the Euler product corresponding to a Steinhaus random multiplicative function (as before), we have
$$ \tilde{\p}(-a -Bj \leq \sum_{m=1}^{j} \log|I_{l_m}(1/2 + \sigma + it_m)| \leq a + j + h(j) \; \forall j \leq n) \asymp \min\{1,\frac{a}{\sqrt{n}}\} . $$
\end{prop5}

\begin{proof}[Proof of Proposition 5]
We note first that $\log x_j = e^{-(l_j + 1)}\log x = e^{B+j} \frac{\log x}{e^{\lfloor \log\log x \rfloor}}$ here, so our assumptions imply that $|\sigma| \leq 1/\log x_n$ and that $|t_j| \leq 1/(j^{2/3} \log x_j)$ for all $1 \leq j \leq n$. Furthermore, provided $B$ is fixed large enough then $x_1$ will be sufficiently large that Lemmas 3 and 4 may legitimately be applied.

We can lower bound the probability we are interested in by
$$ \tilde{\p}(-Bj \leq \sum_{m=1}^{j} \log|I_{l_m}(1/2 + \sigma + it_m)| \leq \min\{a,Bj\} + j + h(j) \; \forall j \leq n) , $$
and here since $\log x_j \geq e^{B+j}$ we have $Bj \leq (1/80)\sqrt{\log x_j}$ and $\min\{a,Bj\} + j + h(j) \leq (1/80)\sqrt{\log x_j}$ for all $1 \leq j \leq n$, provided $B$ is fixed large enough. Thus Lemma 4 is applicable, and yields that our probability is
$$ \gg \p(-(B+1)j + O(1) \leq \sum_{m=1}^{j} G_m \leq \min\{a,Bj\} + h(j) - O(1) \; \forall j \leq n) , $$
where $G_j$ are independent Gaussians with mean 0 and variance $\sum_{x_{j}^{1/e} < p \leq x_j} \frac{1}{2p^{1 + 2\sigma}}$. In particular, the variance here is $\leq \frac{e^{(2\log x_j)/\log x_n}}{2} \sum_{x_{j}^{1/e} < p \leq x_j} \frac{1}{p}  \leq 5$, and similarly the variance is $\geq 1/20$, so we can use Probability Result 2 to deduce the desired lower bound
$$ \p(-(B+1)j + O(1) \leq \sum_{m=1}^{j} G_m \leq \min\{a,Bj\} + h(j) - O(1) \; \forall j \leq n) \gg \min\{1,\frac{a}{\sqrt{n}}\} . $$

Obtaining a matching upper bound for our probability will be slightly more involved, though not too much so. If $a \geq \sqrt{n}$ then the upper bound is trivial, so we may assume instead that $a < \sqrt{n}$. If our event occurs, we must have
$$ -(B+1)a \leq \sum_{m \leq a} \log|I_{l_m}(1/2 + \sigma + it_m)| \leq 2a + h(\lfloor a \rfloor) , $$
and this will be independent of the behaviour of the summands with $m > a$. So we may upper bound the probability we are interested in by
$$ \tilde{\p}(-3a - h(\lfloor a \rfloor) - Bj \leq \sum_{a < m \leq j} \log|I_{l_m}(1/2 + \sigma + it_m)| \leq (B+2)a + j + h(j) \; \forall a < j \leq n) . $$
Now the point is that we have $3a + h(\lfloor a \rfloor) + Bj \leq (1/80)\sqrt{\log x_j}$ and $(B+2)a + j + h(j) \leq (1/80)\sqrt{\log x_j}$ for $j > a$, so Lemma 4 is now applicable and we may upper bound our probability by
$$ \p(-2a -(B+1)j -h(\lfloor a \rfloor) - O(1) \leq \sum_{a < m \leq j} G_m \leq (B+2)a + a + h(j) + O(1) \; \forall a < j \leq n) . $$
Notice here that Lemma 4 is applied with $j$ replaced by $j - \lfloor a \rfloor$, hence the added $a$ terms in the upper and lower limits. Finally Probability Result 1 is applicable, and since $B$ is an absolute constant and $n-a \geq n - \sqrt{n} \gg n$ we obtain our desired upper bound $\ll \min\{1, \frac{Ba}{\sqrt{n}}\} \ll \min\{1, \frac{a}{\sqrt{n}}\}$ for the probability. 
\end{proof}

\subsection{Girsanov-type calculations in the Rademacher case}\label{subsecradgirsanov}
In this subsection, we develop an analogue of Proposition 5 to cover Rademacher random multiplicative functions. The set-up and arguments are very similar as in the Steinhaus case, broadly speaking, but there is the usual complication that $f(n)n^{it}$ no longer has the same law as $f(n)$, and both the formulation and proofs must be adjusted to address this. One already sees this issue when comparing the Euler product calculations in Lemma 2 with those in Lemma 1.

Thus for each $t \in \R$, and for large $x$ and $-1/100 \leq \sigma \leq 1/100$, we define a tilted probability measure $\tilde{\p}_{t}^{\text{Rad}} = \tilde{\p}_{x,\sigma,t}^{\text{Rad}}$ by setting
$$ \tilde{\p}_{t}^{\text{Rad}}(A) := \frac{\E \textbf{1}_{A} \prod_{p \leq x^{1/e}} \left|1 + \frac{f(p)}{p^{1/2+\sigma+it}}\right|^{2}}{\E \prod_{p \leq x^{1/e}} \left|1 + \frac{f(p)}{p^{1/2+\sigma+it}}\right|^{2}} $$
for each event $A$, where $f$ is a Rademacher random multiplicative function. Again, we will sometimes write $\tilde{\E}_{t}^{\text{Rad}}$ to denote expectation (i.e. integration) with respect to the measure $\tilde{\p}_{t}^{\text{Rad}}$.

Now for each $l \in \N \cup \{0\}$ we set $I_l(s) := \prod_{x^{e^{-(l+2)}} < p \leq x^{e^{-(l+1)}}} (1 + \frac{f(p)}{p^s})$, the $l$-th ``increment'' of the Euler product corresponding to $f$. There should be no confusion with the corresponding notation from the Steinhaus case, since one simply takes increments of the appropriate Euler product for the kind of random multiplicative function one is working with. As before, we let $(l_j)_{j=1}^{n}$ denote a {\em strictly decreasing} sequence of non-negative integers, with $l_1 \leq \log\log x - 2$, and define a corresponding increasing sequence of real numbers $(x_j)_{j=1}^{n}$ by setting $x_j := x^{e^{-(l_j + 1)}}$. 

\vspace{12pt}
The following result is the Rademacher analogue of Lemma 4. Note that in Lemma 4 we had $t=0$ (although, because of translation invariance in law, this wasn't actually a restriction), so the condition on $|t_j - t|$ below is analogous to the condition on $|t_j|$ that we had there.

\begin{lem5}
Let $t \in \R$, and let the situation be as above.

Suppose that $x_1 \geq \max\{e^{C/|t|},e^{C\log^{2}|t|}\}$ is large, and that $|\sigma| \leq 1/\log x_n$, and that the numbers $(t_j)_{j=1}^{n}$ satisfy $|t_j - t| \leq \frac{1}{j^{2/3} \log x_j}$. Suppose $(u_j)_{j=1}^{n}$ and $(v_j)_{j=1}^{n}$ are sequences of real numbers satisfying
$$ -(1/80)\sqrt{\log x_j} \leq u_j \leq v_j \leq (1/80)\sqrt{\log x_j} \;\;\; \forall 1 \leq j \leq n . $$
Then we have
\begin{eqnarray}
&& \p((u_j - j) + O(1) \leq \sum_{m=1}^{j} G_m \leq (v_j-j) - O(1) \; \forall j \leq n) \nonumber \\
& \ll & \tilde{\p}_{t}^{\text{Rad}}(u_j \leq \sum_{m=1}^{j} \log|I_{l_m}(1/2 + \sigma + it_m)| \leq v_j \; \forall j \leq n) \nonumber \\
& \ll & \p((u_j - j) - O(1) \leq \sum_{m=1}^{j} G_m \leq (v_j-j) + O(1) \; \forall j \leq n) , \nonumber
\end{eqnarray}
where $G_m$ are independent Gaussian random variables, each having mean 0 and variance $\sum_{x_{m}^{1/e} < p \leq x_m} \frac{1+\cos(2t_m \log p)}{2p^{1 + 2\sigma}}$.
\end{lem5}

\begin{proof}[Proof of Lemma 5]
Lemma 2 implies that for any $|u| \leq x_{j}^{1/20}$, say, the characteristic function $\tilde{\E}_{t}^{\text{Rad}} e^{iu\log|I_{l_j}(1/2 + \sigma + it_j)|}$ is
$$ = \exp\{\sum_{x_{j}^{1/e} < p \leq x_j} \frac{iuc(t , t_j - t , p) - (u^{2}/4)(1+\cos(2t_j \log p))}{p^{1 + 2\sigma}} + T(u) - T(0)\} , $$
where $c(t , t_j - t , p) = 2\cos(t \log p) \cos(t_j\log p) - (1/2)\cos(2 t_j \log p)$. Using a standard trigonometric identity, together with our assumption on $|t_j - t|$, we have
\begin{eqnarray}
2\cos(t \log p) \cos(t_j\log p) & = & \cos((t+t_j)\log p) + \cos((t-t_j)\log p) \nonumber \\
& = & \cos((t+t_j)\log p) + 1 + O(1/j^{4/3}) . \nonumber
\end{eqnarray}
As noted in the proof of Lemma 4, we have $\sum_{x_{j}^{1/e} < p \leq x_j} \frac{1}{p^{1+2\sigma}} = 1 + O(\frac{1}{\log x_j} + \frac{\log x_j}{\log x_n})$. Furthermore, since we assume that $x_1 \geq \max\{e^{C/|t|},e^{C\log^{2}|t|}\}$ is large, which implies that $x_{j}^{1/e} \geq x_{1}^{1/e} \geq \max\{e^{C/e|t|}, e^{(C/e)\log^{2}|t|}\}$, a strong form of the Prime Number Theorem implies that both
$$ \sum_{x_{j}^{1/e} < p \leq x_j} \frac{\cos((t+t_j)\log p)}{p^{1 + 2\sigma}} , \;\;\; \sum_{x_{j}^{1/e} < p \leq x_j} \frac{\cos(2 t_j \log p)}{p^{1 + 2\sigma}} \ll \frac{1}{|t|\log x_j} \ll \frac{1}{e^{j}|t| \log x_1} \ll \frac{1}{Ce^{j}} . $$
See section 6.1 of Harper~\cite{harpergp} for details of such calculations.

Having made these preliminary observations, we can check the proofs of Lemma 3 and Lemma 4, and see they carry over to the present case. Here the relevant Gaussian random variables $N_j$ will have mean $\sum_{x_{j}^{1/e} < p \leq x_j} \frac{c(t , t_j - t , p)}{p^{1 + 2\sigma}}$ and variance $\sum_{x_{j}^{1/e} < p \leq x_j} \frac{1+\cos(2t_j \log p)}{2p^{1 + 2\sigma}}$ (rather than mean $\sum_{x_{j}^{1/e} < p \leq x_j} \frac{\cos(t_j \log p)}{p^{1 + 2\sigma}}$ and variance $\sum_{x_{j}^{1/e} < p \leq x_j} \frac{1}{2p^{1 + 2\sigma}}$, in the Steinhaus case), but our assumptions on $t_j , x_j$ and the above calculations show we still have $\sum_{x_{j}^{1/e} < p \leq x_j} \frac{c(t , t_j - t , p)}{p^{1 + 2\sigma}} = 1 + O(\frac{1}{\log x_j} + \frac{\log x_j}{\log x_n} + \frac{1}{j^{4/3}}) = 1 + O(\frac{1}{e^{n-j}} + \frac{1}{j^{4/3}})$, as well as $\sum_{x_{j}^{1/e} < p \leq x_j} \frac{1+\cos(2t_j \log p)}{2p^{1 + 2\sigma}} = \sum_{x_{j}^{1/e} < p \leq x_j} \frac{1}{2p^{1 + 2\sigma}} + O(1/C)$. Thus the arguments of Lemma 3 and Lemma 4 go through to prove Lemma 5.
\end{proof}

Combining Lemma 5 with Probability Results 1 and 2 leads to the following proposition, which is a Rademacher analogue of Proposition 5.

\begin{prop6}
There is a large natural number $B$ such that the following is true.

Let $t \in \R$, and let $D \geq \max\{\log(1/|t|), 2\log\log(1+|t|)\} + (B+1)$ be any natural number. Let $n \leq \log\log x - D$ be large, and define the decreasing sequence $(l_j)_{j=1}^{n}$ of non-negative integers by $l_j := \lfloor \log\log x \rfloor - D - j$. Suppose that $|\sigma| \leq \frac{1}{e^{D+n}}$, and that $(t_j)_{j=1}^{n}$ is a real sequence satisfying $|t_j - t| \leq \frac{1}{j^{2/3} e^{D+j}}$ for all $j$.

Then uniformly for any large $a$ and any function $h(n)$ satisfying $|h(n)| \leq 10\log n$, and with $I_{l}(s)$ denoting the increments of the Euler product corresponding to a Rademacher random multiplicative function (as before), we have
$$ \tilde{\p}_{t}^{\text{Rad}}(-a -Bj \leq \sum_{m=1}^{j} \log|I_{l_m}(1/2 + \sigma + it_m)| \leq a + j + h(j) \; \forall j \leq n) \asymp \min\{1,\frac{a}{\sqrt{n}}\} . $$
\end{prop6}

\begin{proof}[Proof of Proposition 6]
The only non-trivial change from the proof of Proposition 5 is that we must verify the extra condition $x_1 \geq \max\{e^{C/|t|},e^{C\log^{2}|t|}\}$ in Lemma 5. However, since we have $\log x_j = e^{-(l_j + 1)} \log x = e^{D+j-1} \frac{\log x}{e^{\lfloor \log\log x \rfloor}}$, in particular we have $\log x_1 \geq e^{D} \geq e^{B+1} \max\{\frac{1}{|t|}, \log^{2}(1+|t|)\}$, so the condition will be satisfied provided $B$ is fixed large enough in terms of the absolute constant $C$ in Lemma 5.
\end{proof}

\section{Proofs of the upper bounds in Theorems 1 and 2}\label{secmainupper}

\subsection{The upper bound in the Steinhaus case}\label{subsecsteinupper}
For each $|t| \leq 1/2$, set $t(-1) = t$, and then iteratively for each $0 \leq j \leq \log\log x - 2$ define
$$ t(j) := \max\{u \leq t(j-1): u = \frac{n}{((\log x)/e^{j+1}) \log((\log x)/e^{j+1})} \; \text{for some} \; n \in \Z\} . $$
Thus the points $t(j)$ form a sequence of approximations to $t$, in which for each $j$ we have
\begin{eqnarray}\label{approxdist}
|t-t(j)| = t - t(j) = \sum_{l=0}^{j} t(l-1)-t(l) & \leq & \sum_{l=0}^{j} \frac{1}{((\log x)/e^{l+1}) \log((\log x)/e^{l+1})} \nonumber \\
& \leq & \frac{2}{((\log x)/e^{j+1}) \log((\log x)/e^{j+1})} . 
\end{eqnarray}
With our choice of labelling these approximations become coarser (i.e. worse) as $j$ increases, but recall that in our labelling the lengths of our Euler products $F_{j}(s)$ also go down as $j$ increases.

Furthermore, recall that we let $I_l(s) = \prod_{x^{e^{-(l+2)}} < p \leq x^{e^{-(l+1)}}} (1 - \frac{f(p)}{p^s})^{-1}$, the $l$-th ``increment'' of the Euler product corresponding to Steinhaus $f$.

Given this notation, let $B$ be the large fixed natural number from Proposition 5, and let $\mathcal{G}(k)$ denote the event that for all $|t| \leq 1/2$ and all $k \leq j \leq \log\log x - B - 2$, we have
$$ \Biggl( \frac{\log x}{e^{j+1}} e^{g(x,j)} \Biggr)^{-1} \leq \prod_{l = j}^{\lfloor \log\log x \rfloor - B - 2} |I_{l}(1/2 - \frac{k}{\log x} + it(l))| \leq \frac{\log x}{e^{j+1}} e^{g(x,j)} , $$
where $g(x,j) :=  C\min\{\sqrt{\log\log x}, \frac{1}{1-q} \} + 2\log\log(\frac{\log x}{e^{j+1}})$ for a large constant $C$. Thus $\mathcal{G}(k)$ is the event that our Euler product isn't too large or too small on any ``scale'' $k \leq j \leq \log\log x - B - 2$. The fact that one has a different point $t(l)$ in each part $I_l$ of the Euler product is a little inelegant: it would be nice to keep $t$ itself everywhere, in which case for each $j$ the product would essentially just be $|F_{j}(1/2 - k/\log x + it)|$. But this technical device of varying $t$ slightly, and thereby discretising the set of $t$, will make it easier to prove an important estimate we shall need (Key Proposition 2, below). Similarly, the lower bound condition in the definition of $\mathcal{G}(k)$ is really a technical device to make Proposition 5 applicable: it is the upper bound condition that imposes a significant constraint.

\vspace{12pt}
With the above preparations we can state two key estimates, from which we will deduce the upper bound part of Theorem 1.

\begin{keyprop1}
For all large $x$, and uniformly for $0 \leq k \leq \mathcal{K} = \lfloor \log\log\log x \rfloor$ and $2/3 \leq q \leq 1$, we have
$$ \E( \textbf{1}_{\mathcal{G}(k)} \int_{-1/2}^{1/2} |F_{k}(1/2 - \frac{k}{\log x} + it)|^2 dt )^q \ll \Biggl( \frac{\log x}{e^{k}} C \min\{1, \frac{1}{(1-q)\sqrt{\log\log x}} \} \Biggr)^q , $$
where $\textbf{1}$ denotes the indicator function.
\end{keyprop1}

\begin{keyprop2}
For all large $x$, and uniformly for $0 \leq k \leq \mathcal{K} = \lfloor \log\log\log x \rfloor$ and $2/3 \leq q \leq 1$, we have
$$ \p(\mathcal{G}(k) \; \text{fails}) \ll e^{-2C \min\{\sqrt{\log\log x}, \frac{1}{1-q} \}} . $$
\end{keyprop2}

\begin{proof}[Proof of the upper bound in Theorem 1, assuming Key Propositions 1 and 2]
If $q$ satisfies $1 - \frac{1}{\sqrt{\log\log x}} \leq q \leq 1$ then the upper bound we need to prove is $\E|\sum_{n \leq x} f(n)|^{2q} \ll x^q$, which is trivial in view of H\"{o}lder's inequality and the estimate $\E|\sum_{n \leq x} f(n)|^2 \leq x$. In view of this and of Proposition 1, the upper bound in Theorem 1 will follow if we can show that
$$ \E(\frac{e^{k} (1-q) \sqrt{\log\log x}}{\log x} \int_{-1/2}^{1/2} |F_{k}(1/2 - \frac{k}{\log x} + it)|^2 dt )^q \ll 1 , $$
uniformly for all $0 \leq k \leq \mathcal{K} = \lfloor \log\log\log x \rfloor$ and $2/3 \leq q \leq 1 - \frac{1}{\sqrt{\log\log x}}$.

To prove this, for $1/\sqrt{\log\log x} \leq \delta \leq 1/6$ (say) we define
$$ R(\delta) = R(\delta,k,x) := \sup_{1-2\delta \leq q \leq 1-\delta} \E(\frac{e^{k} (1-q) \sqrt{\log\log x}}{\log x} \int_{-1/2}^{1/2} |F_{k}(1/2 - \frac{k}{\log x} + it)|^2 dt )^q . $$
We have
\begin{eqnarray}
&& \E( \int_{-1/2}^{1/2} |F_{k}(1/2 - \frac{k}{\log x} + it)|^2 dt )^q \nonumber \\
& \leq & \E( \textbf{1}_{\mathcal{G}(k)} \int_{-1/2}^{1/2} |F_{k}(\frac{1}{2} - \frac{k}{\log x} + it)|^2 dt )^q + \E( \textbf{1}_{\mathcal{G}(k) \; \text{fails}} \int_{-1/2}^{1/2} |F_{k}(\frac{1}{2} - \frac{k}{\log x} + it)|^2 dt )^q , \nonumber
\end{eqnarray}
so applying Key Proposition 1 with $C$ a large fixed constant, we obtain
\begin{eqnarray}
R(\delta) & \ll & C^{q} + \sup_{1-2\delta \leq q \leq 1-\delta} \E(\textbf{1}_{\mathcal{G}(k) \; \text{fails}} \frac{e^{k} (1-q) \sqrt{\log\log x}}{\log x} \int_{-1/2}^{1/2} |F_{k}(\frac{1}{2} - \frac{k}{\log x} + it)|^2 dt )^q \nonumber \\
& \leq & C + \sup_{1-2\delta \leq q \leq 1-\delta} \E(\textbf{1}_{\mathcal{G}(k) \; \text{fails}} \frac{e^{k} (1-q) \sqrt{\log\log x}}{\log x} \int_{-1/2}^{1/2} |F_{k}(\frac{1}{2} - \frac{k}{\log x} + it)|^2 dt )^q . \nonumber 
\end{eqnarray}

Next, for each $1-2\delta \leq q \leq 1-\delta$ we set $q' = (1+q)/2$, so that $1-\delta \leq q' \leq 1-\delta/2$. Then by H\"{o}lder's inequality with exponents $q'/(q'-q)$ and $q'/q$, we have
\begin{eqnarray}
&& \E(\textbf{1}_{\mathcal{G}(k) \; \text{fails}} \frac{e^{k} (1-q) \sqrt{\log\log x}}{\log x} \int_{-1/2}^{1/2} |F_{k}(1/2 - \frac{k}{\log x} + it)|^2 dt )^q \nonumber \\
& \leq & \E(\textbf{1}_{\mathcal{G}(k) \; \text{fails}})^{(q'-q)/q'} \Biggl(\E(\frac{e^{k} (1-q) \sqrt{\log\log x}}{\log x} \int_{-1/2}^{1/2} |F_{k}(1/2 - \frac{k}{\log x} + it)|^2 dt )^{q'} \Biggr)^{q/q'} \nonumber \\
& \ll & \p(\mathcal{G}(k) \; \text{fails})^{\delta/2} \Biggl(\E(\frac{e^{k} (1-q') \sqrt{\log\log x}}{\log x} \int_{-1/2}^{1/2} |F_{k}(1/2 - \frac{k}{\log x} + it)|^2 dt )^{q'} \Biggr)^{q/q'} . \nonumber
\end{eqnarray}
The point is that, by Key Proposition 2, we have $\p(\mathcal{G}(k) \; \text{fails}) \ll e^{-2C \min\{\sqrt{\log\log x}, \frac{1}{1-q} \}}\leq e^{-C/\delta}$. So substituting back above, noting that always $q/q' \leq 1$, we deduce the recursive bound
$$ R(\delta) \ll C + e^{-C/2} (1 + R(\delta/2)) \ll C + e^{-C/2}R(\delta/2) . $$

Iterating the recursive bound (with $C$ fixed sufficiently large to compensate for the implicit constant there), replacing $\delta$ by $\delta/2, \delta/4, \delta/8$, etc., we see that uniformly for $1/\sqrt{\log\log x} \leq \delta \leq 1/6$ we have
$$ R(\delta) \ll 1 + R(1/\sqrt{\log\log x}) . $$
However, using H\"{o}lder's inequality (and using the fact that $f(n)n^{it}$ has the same law as $f(n)$ in the Steinhaus case, and using \eqref{stprodusual}) we have trivially that
\begin{eqnarray}
R(1/\sqrt{\log\log x}) & \ll & \sup_{1-2/\sqrt{\log\log x} \leq q \leq 1-1/\sqrt{\log\log x}} \E(\frac{e^{k}}{\log x} \int_{-1/2}^{1/2} |F_{k}(\frac{1}{2} - \frac{k}{\log x} + it)|^2 dt )^q \nonumber \\
& \leq & \sup_{1-2/\sqrt{\log\log x} \leq q \leq 1-1/\sqrt{\log\log x}} (\frac{e^{k}}{\log x} \int_{-1/2}^{1/2} \E|F_{k}(\frac{1}{2} - \frac{k}{\log x} + it)|^2 dt )^q \nonumber \\
& = & \sup_{1-2/\sqrt{\log\log x} \leq q \leq 1-1/\sqrt{\log\log x}} (\frac{e^{k}}{\log x} \E|F_{k}(\frac{1}{2} - \frac{k}{\log x})|^2 )^q \nonumber \\
& \ll & 1 . \nonumber
\end{eqnarray}
Inserting this above implies that $R(\delta) \ll 1$ for all $1/\sqrt{\log\log x} \leq \delta \leq 1/6$, which yields our upper bound theorem.
\end{proof}

\subsection{Proof of Key Proposition 1}
By H\"{o}lder's inequality, to prove Key Proposition 1 it will suffice to show that
$$ \E \textbf{1}_{\mathcal{G}(k)} \int_{-1/2}^{1/2} |F_{k}(1/2 - \frac{k}{\log x} + it)|^2 dt \ll \frac{\log x}{e^{k}} C \min\{1, \frac{1}{(1-q)\sqrt{\log\log x}} \} , $$
uniformly for $0 \leq k \leq \mathcal{K} = \lfloor \log\log\log x \rfloor$ and $2/3 \leq q \leq 1$.

We can upper bound the left hand side by
$$ \int_{-1/2}^{1/2} \E \textbf{1}_{\mathcal{G}(k,t)} |F_{k}(1/2 - \frac{k}{\log x} + it)|^2 dt , $$
where $\mathcal{G}(k,t)$ denotes the event that
$$ \Biggl(\frac{\log x}{e^{j+1}} e^{g(x,j)} \Biggr)^{-1} \leq \prod_{l = j}^{\lfloor \log\log x \rfloor - B - 2} |I_{l}(1/2 - \frac{k}{\log x} + it(l))| \leq \frac{\log x}{e^{j+1}} e^{g(x,j)} $$
for all $k \leq j \leq \log\log x - B - 2$. This is an upper bound because $\mathcal{G}(k)$ is the event that $\mathcal{G}(k,t)$ holds for all $|t| \leq 1/2$. Furthermore, since the law of $f(n)$ is the same as the law of $f(n)n^{-it}$ we have
$$ \E \textbf{1}_{\mathcal{G}(k,t)} |F_{k}(1/2 - \frac{k}{\log x} + it)|^2 = \E \textbf{1}_{\mathcal{H}(k,t)} |F_{k}(1/2 - \frac{k}{\log x})|^2 , $$
where $\mathcal{H}(k,t)$ denotes the event that
\begin{eqnarray}
\Biggl(\frac{\log x}{e^{j+1}} e^{g(x,j)} \Biggr)^{-1} \leq \prod_{l = j}^{\lfloor \log\log x \rfloor - B - 2} |I_{l}(\frac{1}{2} - \frac{k}{\log x} + i(t(l) - t))| & \leq & \frac{\log x}{e^{j+1}} e^{g(x,j)} \nonumber \\
& \leq & e^{\lfloor \log\log x \rfloor - (B+1) - j} e^{B+1 + g(x,j)} \nonumber
\end{eqnarray}
for all $k \leq j \leq \log\log x - B - 2$.

Now if we take logarithms, and recall that
\begin{eqnarray}
g(x,j) & = & 2\log\log(\frac{\log x}{e^{j+1}}) + C\min\{\sqrt{\log\log x}, \frac{1}{1-q} \} \nonumber \\
& \leq & 2\log(\lfloor \log\log x \rfloor - (B+1) - j) + 2\log(B+2) + C\min\{\sqrt{\log\log x}, \frac{1}{1-q} \} , \nonumber
\end{eqnarray}
we see $\mathcal{H}(k,t)$ is an event of the form treated in Proposition 5, taking $n = \lfloor \log\log x \rfloor - (B + 1) - k$; and $\sigma = - \frac{k}{\log x}$; and $t_m = t(\lfloor \log\log x \rfloor - (B+1) - m) - t$ for all $m$; and
$$ a = C\min\{\sqrt{\log\log x}, \frac{1}{1-q} \} + (B+1) + 2\log(B+2) , \;\;\;\;\; h(n) = 2\log n . $$
We may check that these parameters do satisfy the condition $|\sigma| \leq \frac{1}{e^{B+n+1}}$, and (using the approximation bound \eqref{approxdist}) that $|t_m| \leq \frac{1}{m^{2/3} e^{B+m+1}}$ for all $m$. So Proposition 5 implies\footnote{Here the lower bound we impose on our product corresponds to a lower bound $-(a + j + h(j))$ in Proposition 5. This is a more stringent bound than the barrier $-a - Bj$ allowed there, so the Proposition 5 upper bound is certainly still applicable.} that
$$ \frac{\E \textbf{1}_{\mathcal{H}(k,t)} |F_{k}(1/2 - \frac{k}{\log x})|^2}{\E |F_{k}(1/2 - \frac{k}{\log x})|^2} = \tilde{\p}(\mathcal{H}(k,t)) \ll \min\{1, \frac{a}{\sqrt{n}}\} \ll \min\{1, C, \frac{C}{(1-q)\sqrt{\log\log x}} \} . $$

Finally, Key Proposition 1 follows by combining the above display with the fact that
\begin{eqnarray}
\E |F_{k}(1/2 - \frac{k}{\log x})|^2 & = & \exp\{\sum_{p \leq x^{e^{-(k+1)}}} \frac{1}{p^{1 - 2k/\log x}} + O(1) \} \nonumber \\
& = & \exp\{\sum_{p \leq x^{e^{-(k+1)}}} \frac{1}{p} + O(\sum_{p \leq x^{e^{-(k+1)}}} \frac{k \log p}{p \log x} + 1) \} \ll \frac{\log x}{e^k} , \nonumber
\end{eqnarray}
which follows from \eqref{stprodusual} and the Mertens and Chebychev estimates for sums over primes (see e.g. Theorem 2.7 of Montgomery and Vaughan~\cite{mv}).
\qed

\subsection{Proof of Key Proposition 2}
By the union bound, we see $\p(\mathcal{G}(k) \; \text{fails})$ is
\begin{eqnarray}
& \leq & \sum_{k \leq j \leq \log\log x - B - 2} \p\Biggl(\prod_{l=j}^{\lfloor \log\log x \rfloor - B - 2} |I_{l}(\frac{1}{2} - \frac{k}{\log x} + it(l))| > \frac{\log x}{e^{j+1}} e^{g(x,j)} \; \text{for some} \; |t| \leq \frac{1}{2} \Biggr) \nonumber \\
&& + \sum_{k \leq j \leq \log\log x - B - 2} \p\Biggl(\prod_{l=j}^{\lfloor \log\log x \rfloor - B - 2} |I_{l}(\frac{1}{2} - \frac{k}{\log x} + it(l))|^{-1} > \frac{\log x}{e^{j+1}} e^{g(x,j)} \; \text{for some} \; |t| \leq \frac{1}{2} \Biggr) \nonumber
\end{eqnarray}
Let us temporarily write $\Sigma_{1}$, $\Sigma_{2}$ to denote these two sums.

We concentrate on trying to bound $\Sigma_{1}$. Because of the definition of $t(l)$, we see the probability inside the sum is at most as large as with the set $|t| \leq 1/2$ replaced by the set
$$ \mathcal{T}(x,j) := \Biggl\{\frac{n}{((\log x)/e^{j+1}) \log((\log x)/e^{j+1})} : |n| \leq ((\log x)/e^{j+1}) \log((\log x)/e^{j+1}) \Biggr\} . $$
Since we are now dealing with a discrete set of points, we can apply the union bound again to obtain that
$$ \Sigma_{1} \leq \sum_{k \leq j \leq \log\log x - B - 2} \sum_{t(j) \in \mathcal{T}(x,j)} \p\Biggl(\prod_{l=j}^{\lfloor \log\log x \rfloor - B - 2} |I_{l}(1/2 - \frac{k}{\log x} + it(l))| > \frac{\log x}{e^{j+1}} e^{g(x,j)} \Biggr) , $$
and by Chebychev's inequality this is all
$$ \leq \sum_{k \leq j \leq \log\log x - B - 2} \sum_{t(j) \in \mathcal{T}(x,j)} \frac{1}{(\frac{\log x}{e^{j+1}} e^{g(x,j)})^2} \E \prod_{l=j}^{\lfloor \log\log x \rfloor - B - 2} |I_{l}(1/2 - \frac{k}{\log x} + it(l))|^2 . $$

Finally, since $f$ is independent on distinct primes, and (in the Steinhaus case) its law is the same as the law of $f(n)n^{it}$ for any fixed $t$, the above is
\begin{eqnarray}
& = & \sum_{k \leq j \leq \log\log x - B - 2} \sum_{t(j) \in \mathcal{T}(x,j)} \frac{1}{(\frac{\log x}{e^{j+1}} e^{g(x,j)})^2}  \prod_{l=j}^{\lfloor \log\log x \rfloor - B - 2} \E |I_{l}(1/2 - \frac{k}{\log x})|^2 \nonumber \\
& \ll & \sum_{k \leq j \leq \log\log x - B - 2} e^{-2g(x,j)} \frac{\log((\log x)/e^{j+1})}{\frac{\log x}{e^{j+1}}}  \prod_{l=j}^{\lfloor \log\log x \rfloor - B - 2} \E |I_{l}(1/2 - \frac{k}{\log x})|^2 . \nonumber
\end{eqnarray}
As noted in \eqref{stprodusual} (and in the proof of Key Proposition 1), Lemma 1 implies the product of expectations is $\ll \exp\{\sum_{p \leq x^{e^{-(j+1)}}} \frac{1}{p^{1 - 2k/\log x}}\} \ll \frac{\log x}{e^{j+1}}$, and inserting this and the fact that $g(x,j) :=  C\min\{\sqrt{\log\log x}, \frac{1}{1-q} \} + 2\log\log(\frac{\log x}{e^{j+1}})$ we obtain
$$ \Sigma_{1} \ll \sum_{k \leq j \leq \log\log x - B - 2} e^{-2C\min\{\sqrt{\log\log x}, \frac{1}{1-q} \} - 3\log\log(\frac{\log x}{e^{j+1}}) } \ll e^{-2C\min\{\sqrt{\log\log x}, \frac{1}{1-q} \}} . $$

One can bound $\Sigma_{2}$ in exactly the same way, since (as remarked following \eqref{stprodusual}) one has the same estimate for $\E |I_{l}(1/2 - \frac{k}{\log x})|^{-2}$ as for $\E |I_{l}(1/2 - \frac{k}{\log x})|^2$. This completes the proof of Key Proposition 2.
\qed

\subsection{The upper bound in the Rademacher case}
Let $x$ be large, and let $0 \leq k \leq \mathcal{K} = \lfloor \log\log\log x \rfloor$ and $2/3 \leq q \leq 1$. We note immediately that, by H\"{o}lder's inequality, we have
\begin{eqnarray}
|| \int_{|t| \leq 1/\sqrt{\log\log x}} |F_{k}(1/2 - \frac{k}{\log x} + it)|^2 dt ||_{q} & \leq & \int_{|t| \leq 1/\sqrt{\log\log x}} \E|F_{k}(1/2 - \frac{k}{\log x} + it)|^2 dt \nonumber \\
& \ll & \frac{\log x}{e^{k} \sqrt{\log\log x}} . \nonumber
\end{eqnarray}
Here we used the estimate $\E|F_{k}(1/2 - \frac{k}{\log x} + it)|^2 \ll \exp\{\sum_{p \leq x^{e^{-(k+1)}}} \frac{1}{p^{1 - 2k/\log x}}\} \ll (\log x)/e^{k}$, which follows from \eqref{radprodusual} and from standard estimates for sums over primes (as in the proofs of Key Propositions 1 and 2). Similarly, when $|N| \geq (\log\log x)^2 $ we have
\begin{eqnarray}
\frac{1}{|N|^{1/4}} || \int_{N-1/2}^{N+1/2} |F_{k}(\frac{1}{2} - \frac{k}{\log x} + it)|^2 dt ||_{q} & \leq & \frac{1}{|N|^{1/4}} \int_{N-1/2}^{N+1/2} \E|F_{k}(\frac{1}{2} - \frac{k}{\log x} + it)|^2 dt \nonumber \\
& \ll & \frac{\log x}{e^{k} \sqrt{\log\log x}} . \nonumber
\end{eqnarray}
If we insert these bounds into Proposition 2, we obtain an acceptable contribution for the Theorem 2 upper bound. Therefore to finish the proof of that upper bound using Proposition 2, it will suffice to show that
$$ \frac{1}{(|N|+1)^{q/4}} \E( \int_{\substack{|t - N| \leq 1/2, \\ |t| > 1/\sqrt{\log\log x} }} |F_{k}(\frac{1}{2} - \frac{k}{\log x} + it)|^2 dt )^q \ll \Biggl( \frac{\log x}{e^{k}} \min\{1, \frac{1}{(1-q)\sqrt{\log\log x}} \} \Biggr)^q , $$
uniformly for $0 \leq k \leq \mathcal{K} = \lfloor \log\log\log x \rfloor$ and $2/3 \leq q \leq 1$ and $|N| \leq (\log\log x)^2$. The preceding reductions, dismissing very small and large $t$ from consideration, will be convenient when we come to apply Proposition 6.

\vspace{12pt}
We will concentrate on the case $N=0$, since this will reveal essentially all the important differences between the Rademacher and Steinhaus arguments. Thus for each $1/\sqrt{\log\log x} < |t| \leq 1/2$, we can define the sequence of approximations $(t(j))_{0 \leq j \leq \log\log x - 2}$ as we did in the Steinhaus case in section \ref{subsecsteinupper}, and have the same bounds \eqref{approxdist} on $|t - t(j)|$ as we did there. Let us further define $D(t) := \lceil \log(1/|t|) \rceil + (B+1)$, where $B$ is as in Proposition 6. Then with $I_l(s) = \prod_{x^{e^{-(l+2)}} < p \leq x^{e^{-(l+1)}}} (1 + \frac{f(p)}{p^s})$ denoting the $l$-th ``increment'' of the Rademacher Euler product, we will let $\mathcal{G}^{\text{Rad}}(k,t)$ denote the event that for all $k \leq j \leq \log\log x - D - 1$, we have
$$ \Biggl( \frac{\log x}{e^{j+1}} e^{g(x,j)} \Biggr)^{-1} \leq \prod_{l = j}^{\lfloor \log\log x \rfloor - D - 1} |I_{l}(1/2 - \frac{k}{\log x} + it(l))| \leq \frac{\log x}{e^{j+1}} e^{g(x,j)} , $$
where $g(x,j) :=  C\min\{\sqrt{\log\log x}, \frac{1}{1-q} \} + 2\log\log(\frac{\log x}{e^{j+1}})$ for a large constant $C$. Furthermore, we let $\mathcal{G}^{\text{Rad}}(k)$ denote the event that $\mathcal{G}^{\text{Rad}}(k,t)$ holds for all $1/\sqrt{\log\log x} < |t| \leq 1/2$. Notice this set-up is as close as possible to what we did in the Steinhaus case, the only real change being the introduction of the term $D(t)$, which ensures we will be able to apply Proposition 6.

Now we have the following two key estimates.

\begin{keyprop3}
For all large $x$, and uniformly for $0 \leq k \leq \mathcal{K} = \lfloor \log\log\log x \rfloor$ and $2/3 \leq q \leq 1$, we have
$$ \E( \textbf{1}_{\mathcal{G}^{\text{Rad}}(k)} \int_{1/\sqrt{\log\log x} < |t| \leq 1/2} |F_{k}(\frac{1}{2} - \frac{k}{\log x} + it)|^2 dt )^q \ll \Biggl( \frac{\log x}{e^{k}} C \min\{1, \frac{1}{(1-q)\sqrt{\log\log x}} \} \Biggr)^q , $$
where $\textbf{1}$ denotes the indicator function.
\end{keyprop3}

\begin{keyprop4}
For all large $x$, and uniformly for $0 \leq k \leq \mathcal{K} = \lfloor \log\log\log x \rfloor$ and $2/3 \leq q \leq 1$, we have
$$ \p(\mathcal{G}^{\text{Rad}}(k) \; \text{fails}) \ll e^{-2C \min\{\sqrt{\log\log x}, \frac{1}{1-q} \}} . $$
\end{keyprop4}

\begin{proof}[Proof of Key Proposition 3]
Almost all the details are the same as in the proof of Key Proposition 1 from the Steinhaus case. Now we cannot translate the event $\mathcal{G}^{\text{Rad}}(k,t)$ by shifting $t$ to 0, as we did in the Steinhaus case, but there is no need to do so because we formulated our tilted probability estimates for $\tilde{\p}_{t}^{\text{Rad}}$ for general $t \in \R$. We must take $a = C\min\{\sqrt{\log\log x}, \frac{1}{1-q} \} + D + 2\log(D+1)$, rather than $a = C\min\{\sqrt{\log\log x}, \frac{1}{1-q} \} + (B+1) + 2\log(B+2)$, and so Proposition 6 and \eqref{radprodusual} imply that
\begin{eqnarray}
\E \textbf{1}_{\mathcal{G}^{\text{Rad}}(k,t)} |F_{k}(1/2 - \frac{k}{\log x} + it)|^2 & \ll & \min\{1, \frac{a}{\sqrt{\log\log x}}\} \E|F_{k}(1/2 - \frac{k}{\log x} + it)|^2 \nonumber \\
& \ll & C D(t) \min\{1, \frac{1}{(1-q)\sqrt{\log\log x}}\} \frac{\log x}{e^k} . \nonumber
\end{eqnarray}
Since we have $\int_{|t| \leq 1/2} D(t) dt \ll \int_{|t| \leq 1/2} \log(1/|t|) dt \ll 1$, this bound works the same as the corresponding bound (without $D(t)$) from the proof of Key Proposition 1, and so Key Proposition 3 follows.
\end{proof}

\begin{proof}[Proof of Key Proposition 4]
The proof of Key Proposition 4 closely follows that of Key Proposition 2 from the Steinhaus case, until\footnote{Strictly speaking, one needs to be careful when applying the union bound in the proof because the upper limit $\lfloor \log\log x \rfloor - D(t) - 1$ in the product now depends on $t$. However, $D(t)$ is constant on $t$ intervals of the form $e^{-(r+1)} \leq |t| < e^{-r}$, so one can split up first according to which such interval $t$ lies in.} the final line where $\Sigma_{2}$ must be bounded. There we must invoke \eqref{radprodinv} to estimate each term $\E |I_{l}(1/2 - \frac{k}{\log x} + it(l))|^{-2}$ for $k \leq l \leq \log\log x - D - 1$, and the estimate this supplies is $\exp\{\sum_{x^{e^{-(l+2)}} < p \leq x^{e^{-(l+1)}}} \frac{1 + 2\cos(2 t(l) \log p)}{p^{1 - 2k/\log x}} + O(\frac{1}{\sqrt{x^{e^{-(l+2)}}} e^{-(l+2)} \log x}) \}$, which isn't identical to the corresponding estimate for $\E |I_{l}(1/2 - \frac{k}{\log x} + it(l))|^{2}$ (which wouldn't include the $2\cos(2t(l)\log p)$ term). However, since we always have $|t| \ll |t(l)| \ll 1$, standard estimates for sums over primes (as in the proof of Lemma 5) show that $\sum_{x^{e^{-(l+2)}} < p \leq x^{e^{-(l+1)}}} \frac{2\cos(2 t(l) \log p)}{p^{1 - 2k/\log x}} \ll \frac{1}{|t| e^{-l} \log x}$. Recalling that $l \leq \log\log x - D - 1$ and that $D = \lceil \log(1/|t|) \rceil + (B+1)$, we see these terms involving $2\cos(2 t(l) \log p)$ give a negligible contribution.
\end{proof}

Given Key Propositions 3 and 4, exactly the same argument as in the Steinhaus case in section \ref{subsecsteinupper} confirms that, as we wanted,
$$ \E( \int_{\substack{|t| \leq 1/2, \\ |t| > 1/\sqrt{\log\log x} }} |F_{k}(1/2 - \frac{k}{\log x} + it)|^2 dt )^q \ll \Biggl( \frac{\log x}{e^{k}} \min\{1, \frac{1}{(1-q)\sqrt{\log\log x}} \} \Biggr)^q . $$

To handle the integral over $|t-N| \leq 1/2$ for general $1 \leq |N| \leq (\log\log x)^2$, one follows the same argument as above, but with $\log\log x - D - 1$ replaced by $\log\log x - \lceil 2\log\log 10N \rceil - (B+1) - 1$, say. This produces an extra factor $\log\log 10N$ in the analogue of Key Proposition 3, which is more than cancelled out by the prefactor $1/(|N|+1)^{q/4}$ attached to the $N$ integral.
\qed

\section{Proofs of the lower bounds in Theorems 1 and 2}\label{secmainlower}
Recall that $F(s)$ denotes the Euler product of $f(n)$ over $x$-smooth numbers. In view of Propositions 3 and 4, the lower bound parts of Theorems 1 and 2 will essentially follow if we can prove suitable lower bounds for $\vert\vert \int_{-1/2}^{1/2} |F(1/2 + \frac{4V}{\log x} + it)|^2 dt \vert\vert_{q}^{1/2}$, where $V$ is a large constant. (We also need upper bounds for quantities like $\vert\vert \int_{-1/2}^{1/2} |F(1/2 + \frac{2V}{\log x} + it)|^2 dt \vert\vert_{q}^{1/2}$, but those will follow directly from our work in section \ref{secmainupper}.) As described in the introduction, we shall actually seek a lower bound for $\vert\vert \int_{\mathcal{L}} |F(1/2 + \frac{4V}{\log x} + it)|^2 dt \vert\vert_{q}^{1/2}$, where $\mathcal{L} \subseteq [-1/2,1/2]$ is a suitable random subset that makes things work nicely.

When we proved Key Proposition 1 (or the Rademacher version, Key Proposition 3) during our work on upper bounds, we used H\"{o}lder's inequality to replace a restricted $q$-th moment by a restricted first moment that we could estimate. This straightforward procedure is not available for lower bounds. However, we can write
$$ \E \int_{\mathcal{L}} |F(\frac{1}{2} + \frac{4V}{\log x} + it)|^2 dt = \E \Biggl( \int_{\mathcal{L}} |F(\frac{1}{2} + \frac{4V}{\log x} + it)|^2 dt \Biggr)^{\frac{q}{2-q}} \Biggl(\int_{\mathcal{L}} |F(\frac{1}{2} + \frac{4V}{\log x} + it)|^2 dt \Biggr)^{\frac{2(1-q)}{(2-q)}} , $$
so applying H\"{o}lder's inequality with exponents $2-q$ and $(2-q)/(1-q)$, we obtain that $\E \int_{\mathcal{L}} |F(\frac{1}{2} + \frac{4V}{\log x} + it)|^2 dt$ is
$$ \leq \Biggl( \E \Biggl( \int_{\mathcal{L}} |F(\frac{1}{2} + \frac{4V}{\log x} + it)|^2 dt \Biggr)^{q} \Biggr)^{\frac{1}{2-q}} \Biggl( \E \Biggl(\int_{\mathcal{L}} |F(\frac{1}{2} + \frac{4V}{\log x} + it)|^2 dt \Biggr)^{2} \Biggr)^{\frac{1-q}{2-q}} . $$
We will be able to lower bound $\E \int_{\mathcal{L}} |F(\frac{1}{2} + \frac{4V}{\log x} + it)|^2 dt$ using the Girsanov-type estimates we already proved in section \ref{secprobcalc}. The extra ingredient we require now is an upper bound for the second moment of the integral over $\mathcal{L}$, and to prove that we must perform some further probabilistic preparations.

\subsection{Further probabilistic calculations}\label{subsecfurtherprobcalc}
In Lemma 1, we computed the expectation of the second power of one Steinhaus Euler product times an imaginary power of a shifted Euler product. We now require a variant of this.

\begin{lem6}
If $f$ is a Steinhaus random multiplicative function, then for any real $t,u,v$, any real $400(1 + u^2 + v^2) \leq x \leq y$, and any real $\sigma \geq - 1/\log y$, we have
\begin{eqnarray}
&& \E \prod_{x < p \leq y} \left|1 - \frac{f(p)}{p^{1/2+\sigma}}\right|^{-(2+iu)} \left|1 - \frac{f(p)}{p^{1/2+\sigma+it}}\right|^{-(2+iv)} \nonumber \\
& = & \exp\{\sum_{x < p \leq y} \frac{(1+iu/2)^2 + (1+iv/2)^2}{p^{1 + 2\sigma}} + \sum_{x < p \leq y} \frac{(2+iu)(2+iv)\cos(t\log p)}{2p^{1+2\sigma}} + T(u,v)\} , \nonumber
\end{eqnarray}
where $T(u,v) = T_{x,y,\sigma,t}(u,v)$ satisfies $|T(u,v)| \ll \frac{1 + |u|^3 + |v|^3}{\sqrt{x} \log x}$, and its partial derivatives satisfy $|\frac{\partial T(u,v)}{\partial u}| \ll \frac{1 + u^2 + v^2}{\sqrt{x} \log x}$, and $|\frac{\partial T(u,v)}{\partial v}| \ll \frac{1 + u^2 + v^2}{\sqrt{x} \log x}$, and $|\frac{\partial T(u,v)}{\partial u \partial v}| \ll \frac{1 + |u| + |v|}{\sqrt{x} \log x}$.
\end{lem6}

In particular, Lemma 6 implies that for any real $t$ and any $400 \leq x \leq y$ and any $\sigma \geq -1/\log y$; and for Steinhaus random multiplicative $f$; we have
\begin{equation}\label{sttwoprods}
\E \prod_{x < p \leq y} \left|1 - \frac{f(p)}{p^{1/2+\sigma}}\right|^{-2} \left|1 - \frac{f(p)}{p^{1/2+\sigma+it}}\right|^{-2} = \exp\{\sum_{x < p \leq y} \frac{2 + 2\cos(t\log p)}{p^{1+2\sigma}} + O(\frac{1}{\sqrt{x}\log x}) \} .
\end{equation}

\begin{proof}[Proof of Lemma 6]
The proof is an extension of the proof of Lemma 1. As there, we temporarily set $R_{p}(t) := -\Re\log(1 - \frac{f(p)}{p^{1/2+\sigma+it}})$, and then we may rewrite
\begin{eqnarray}
\left|1 - \frac{f(p)}{p^{1/2+\sigma}}\right|^{-(2+iu)} \left|1 - \frac{f(p)}{p^{1/2+\sigma+it}}\right|^{-(2+iv)} & = & \exp\{(2+iu) R_{p}(0) + (2+iv)R_{p}(t) \} \nonumber \\
= 1 + \sum_{j=1}^{\infty} \frac{((2+iu)R_{p}(0) + (2+iv)R_{p}(t))^j}{j!} . \nonumber
\end{eqnarray}

In the proof of Lemma 1, we used the Taylor expansion of the logarithm to determine that $\E R_{p}(t) = 0$. We also obtained that
$$ \E R_{p}(t)^2 = \frac{1}{2 p^{1 + 2\sigma}} + O(\frac{1}{p^{3/2+3\sigma}}) , \;\;\; \text{and} \;\;\; \E R_{p}(0) R_{p}(t) = \frac{\cos(t\log p)}{2 p^{1 + 2\sigma}} + O(\frac{1}{p^{3/2+3\sigma}}) , $$
and for $j \geq 3$ we have the trivial bound $|R_{p}(t)^j| \leq (\sum_{k=1}^{\infty} \frac{1}{p^{k(1/2+\sigma)}})^j = \frac{1}{(p^{1/2+\sigma}-1)^j}$. Furthermore, we noted there that for primes $y \geq p > x \geq 400(1+u^2+v^2)$ we have $\frac{1}{p^{1/2+\sigma}} = \frac{e^{-\sigma\log p}}{p^{1/2}} \leq \frac{e}{p^{1/2}}$, which now implies that $(4+|u|+|v|)/p^{1/2 + \sigma} \leq 3e/10$.

So putting things together, for such primes we have
\begin{eqnarray}
&& \E\left|1 - \frac{f(p)}{p^{1/2+\sigma}}\right|^{-(2+iu)} \left|1 - \frac{f(p)}{p^{1/2+\sigma+it}}\right|^{-(2+iv)} \nonumber \\
& = & 1 + \frac{(2+iu)^2 \E R_{p}(0)^2 + 2(2+iu)(2+iv)\E R_{p}(0)R_{p}(t) + (2+iv)^2 \E R_{p}(t)^2}{2} + \nonumber \\
&& + \E \sum_{j=3}^{\infty} \frac{((2+iu)R_{p}(0) + (2+iv)R_{p}(t))^j}{j!} \nonumber \\
& = & 1 + \frac{(1+iu/2)^2 + (1+iv/2)^2}{p^{1+2\sigma}} + \frac{(2+iu)(2+iv)\cos(t\log p)}{2 p^{1+2\sigma}} + O(\sum_{j=3}^{\infty} \frac{(4+|u|+|v|)^j}{j! (p^{1/2+\sigma}-1)^j} ) \nonumber \\
& = & 1 + \frac{(1+iu/2)^2 + (1+iv/2)^2}{p^{1+2\sigma}} + \frac{(2+iu)(2+iv)\cos(t\log p)}{2 p^{1+2\sigma}} + D_{p}(u,v) , \nonumber
\end{eqnarray}
where $D_{p}(u,v)$ satisfies $|D_{p}(u,v)| \ll \frac{1+|u|^3 + |v|^3}{p^{3/2+3\sigma}} \ll \frac{1+|u|^3 + |v|^3}{p^{3/2}}$, and its partial derivatives satisfy $|\frac{\partial D_{p}(u,v)}{\partial u}| \ll \frac{1 + u^2 + v^2}{p^{3/2}}$, and $|\frac{\partial D_{p}(u,v)}{\partial v}| \ll \frac{1 + u^2 + v^2}{p^{3/2}}$, and $|\frac{\partial D_{p}(u,v)}{\partial u \partial v}| \ll \frac{1 + |u| + |v|}{p^{3/2}}$.

The conclusion of Lemma 6 now follows as in the proof of Lemma 1, using the independence of $f$ on distinct primes and the standard Chebychev-type estimate $\sum_{p > x} 1/p^{3/2} \ll 1/(\sqrt{x} \log x)$.
\end{proof}

We remark that in our applications of Lemma 1, we took the shift $t$ to be rather small so that $\cos(t\log p) \approx 1$. In this setting, the expression $\sum_{x < p \leq y} \frac{1 + iu\cos(t\log p) - u^{2}/4}{p^{1 + 2\sigma}}$ from Lemma 1 is approximately the expression $\sum_{x < p \leq y} \frac{(1+iu/2)^2 }{p^{1 + 2\sigma}}$ in Lemma 6. When we come to apply Lemma 6, we will take the shift $t$ there to be somewhat large compared with $x$ so the sum $\sum_{x < p \leq y} \frac{(2+iu)(2+iv)\cos(t\log p)}{2p^{1+2\sigma}}$ is small. The important point is that this regime corresponds to the product and the shifted product in Lemma 6 behaving fairly independently.

\vspace{12pt}
Next, we want to adapt the Girsanov-type calculations from section \ref{subsecstgirsanov} to a ``tilted'' probability measure weighted by two Euler products. Thus for each $t \in \R$, and for large $x$ and $-1/100 \leq \sigma \leq 1/100$, let us define a probability measure $\tilde{\p}_{t}^{\text{St}, (2)} = \tilde{\p}_{x,\sigma,t}^{\text{St}, (2)}$ by setting
$$ \tilde{\p}_{t}^{\text{St}, (2)}(A) := \frac{\E \textbf{1}_{A} \prod_{p \leq x^{1/e}} \left|1 - \frac{f(p)}{p^{1/2+\sigma}}\right|^{-2} \left|1 - \frac{f(p)}{p^{1/2+\sigma + it}}\right|^{-2}}{\E \prod_{p \leq x^{1/e}} \left|1 - \frac{f(p)}{p^{1/2+\sigma}}\right|^{-2} \left|1 - \frac{f(p)}{p^{1/2+\sigma + it}}\right|^{-2}} $$
for each event $A$, where $f$ is a Steinhaus random multiplicative function. We write $\tilde{\E}_{t}^{\text{St}, (2)}$ to denote expectation with respect to this measure. The superscript $(2)$ in this notation reflects the fact we are weighting by two products, and is designed to be distinctive from the notation used in the one product Rademacher case in section \ref{subsecradgirsanov}.

As in the one dimensional case, the exact range of $p$ in the definition of the measure $\tilde{\p}_{t}^{\text{St}, (2)}$ doesn't matter too much, because if the event $A$ doesn't involve a particular prime then the expectation of that part of the products will factor out and cancel between the numerator and denominator. We emphasise this point because, in our later calculations, we will extract various parts of the products and we must be clear that our results about $\tilde{\p}_{t}^{\text{St}, (2)}(A)$ will still be applicable, for appropriate $A$.

As before, we shall let $(l_j)_{j=1}^{n}$ denote a {\em strictly decreasing} sequence of non-negative integers, with $l_1 \leq \log\log x - 2$, and define a corresponding increasing sequence of real numbers $(x_j)_{j=1}^{n}$ by setting $x_j := x^{e^{-(l_j + 1)}}$. And for each $l \in \N \cup \{0\}$ we set $I_l(s) := \prod_{x^{e^{-(l+2)}} < p \leq x^{e^{-(l+1)}}} (1 - \frac{f(p)}{p^s})^{-1}$, the $l$-th increment of the Euler product corresponding to Steinhaus $f$.

Using Lemma 6 (which will serve as a two dimensional characteristic function calculation under the measure $\tilde{\p}_{t}^{\text{St}, (2)}$) and a two dimensional version of the Berry--Esseen inequality, we can prove the following result, which the reader may compare with the one dimensional case in Lemma 3.

\begin{lem7}
Let the situation be as above, with the restriction that $|t| \leq 1$. Suppose that $x_1 \geq e^{C/|t|^2}$ is sufficiently large, and that $|\sigma| \leq 1/\log x_n$. Suppose further that $(u_j)_{j=1}^{n}$ and $(v_j)_{j=1}^{n}$ are any sequences of real numbers satisfying
$$ |u_j| , |v_j| \leq (1/40)(\log x_j)^{1/4} + 2 \;\;\; \forall 1 \leq j \leq n . $$

Then we have
\begin{eqnarray}
&& \tilde{\p}_{t}^{\text{St}, (2)}(u_j \leq \log|I_{l_j}(\frac{1}{2} + \sigma)| \leq u_j + \frac{1}{j^2} , \;\;\; v_j \leq \log|I_{l_j}(\frac{1}{2} + \sigma + it)| \leq v_j + \frac{1}{j^2} ,  \; \forall j \leq n) \nonumber \\
& = & \left(1+O\left(\frac{1}{x_{1}^{1/100}}\right) \right) \p(u_j \leq N_j^{1} \leq u_j + \frac{1}{j^2} , \; \text{and} \; v_j \leq N_{j}^{2} \leq v_j + \frac{1}{j^2} , \; \forall 1 \leq j \leq n) , \nonumber
\end{eqnarray}
where $(N_j^{1}, N_j^{2})_{j=1}^{n}$ is a sequence of independent bivariate Gaussian random vectors, and the components $N_j^{1}, N_j^{2}$ have mean $\sum_{x_{j}^{1/e} < p \leq x_{j}} \frac{1 + \cos(t\log p)}{p^{1 + 2\sigma}}$, variance $\sum_{x_{j}^{1/e} < p \leq x_j} \frac{1}{2p^{1 + 2\sigma}}$, and covariance $\sum_{x_{j}^{1/e} < p \leq x_{j}} \frac{\cos(t\log p)}{2p^{1 + 2\sigma}}$.

Furthermore, we also have the alternative estimate
\begin{eqnarray}
&& \tilde{\p}_{t}^{\text{St}, (2)}(u_j \leq \log|I_{l_j}(\frac{1}{2} + \sigma)| \leq u_j + \frac{1}{j^2} , \;\;\; v_j \leq \log|I_{l_j}(\frac{1}{2} + \sigma + it)| \leq v_j + \frac{1}{j^2} ,  \; \forall j \leq n) \nonumber \\
& = & \left(1+O\left(\frac{1}{\sqrt{C}}\right) \right) \p(u_j \leq N_j^{1} \leq u_j + \frac{1}{j^2} \; \forall j \leq n) \cdot \p(v_j \leq N_{j}^{2} \leq v_j + \frac{1}{j^2} \; \forall j \leq n) , \nonumber
\end{eqnarray}
in other words we may replace the covariance of $N_j^{1} , N_j^{2}$ by zero.
\end{lem7}

We should perhaps comment on some of the assumptions in Lemma 7, as compared with Lemmas 3 and 5. Note that the allowed range of $u_j , v_j$ here is a multiple of $(\log x_j)^{1/4}$, as opposed to our previous $\sqrt{\log x_j}$. The assumption that $x_1 \geq e^{C/|t|^2}$ is also stronger than previous assumptions of the shape $x_1 \geq e^{C/|t|}$. These conditions are used to deduce the second part of Lemma 7, that we may treat $N_j^1$ and $N_j^2$ as independent. Some tradeoff between the conditions is possible, so if one restricted $u_j , v_j$ further one could allow $x_1$ to be smaller. But in our application we will set things up so both conditions are anyway satisfied with room to spare.

\begin{proof}[Proof of Lemma 7]
By independence, both the probability on the left and the ones on the right factor as a product over $j$, so it will suffice to prove that
\begin{eqnarray}
&& \tilde{\p}_{t}^{\text{St}, (2)}(u_j \leq \log|I_{l_j}(\frac{1}{2} + \sigma)| \leq u_j + \frac{1}{j^2} , \; \text{and} \; v_j \leq \log|I_{l_j}(\frac{1}{2} + \sigma + it)| \leq v_j + \frac{1}{j^2}) \nonumber \\
& = & \left(1+O\left(\frac{1}{x_{j}^{1/100}}\right) \right) \p(u_j \leq N_j^{1} \leq u_j + \frac{1}{j^2} , \; \text{and} \; v_j \leq N_{j}^{2} \leq v_j + \frac{1}{j^2}) \nonumber \\
& = & \left(1+O\left(\frac{1}{\sqrt{C e^j}}\right) \right) \p(u_j \leq N_j^{1} \leq u_j + \frac{1}{j^2}) \cdot \p(v_j \leq N_{j}^{2} \leq v_j + \frac{1}{j^2}) \nonumber
\end{eqnarray}
for all $1 \leq j \leq n$.

By Lemma 6, the characteristic function $\tilde{\E}_{t}^{\text{St}, (2)} e^{iu\log|I_{l_j}(1/2 + \sigma)| + iv\log|I_{l_j}(1/2 + \sigma + it)|}$ is
$$ = \exp\{\sum_{x_{j}^{1/e} < p \leq x_j} \frac{iu - \frac{u^{2}}{4} + iv - \frac{v^{2}}{4}}{p^{1 + 2\sigma}} + \sum_{x_{j}^{1/e} < p \leq x_{j}} \frac{(2iu + 2iv - uv)\cos(t\log p)}{2p^{1+2\sigma}} + T(u,v) - T(0,0) \} . $$
Without the error term $T(u,v) - T(0,0)$, a standard calculation shows this would be the characteristic function of the pair $(N_{j}^{1}, N_{j}^{2})$, where $N_{j}^{1}, N_{j}^{2}$ are Gaussian random variables each having mean $\sum_{x_{j}^{1/e} < p \leq x_{j}} \frac{1 + \cos(t\log p)}{p^{1 + 2\sigma}}$ and variance $\sum_{x_{j}^{1/e} < p \leq x_{j}} \frac{1}{2p^{1 + 2\sigma}}$, and with covariance $\E N_{j}^{1} N_{j}^{2} - \E N_{j}^{1} \E N_{j}^{2}= \sum_{x_{j}^{1/e} < p \leq x_{j}} \frac{\cos(t\log p)}{2p^{1 + 2\sigma}}$. Before proceeding further, we record a few calculations we will use later. Firstly, under our conditions $|t| \leq 1$ and $x_j \geq x_1 \geq e^{C/|t|^2}$ (in fact it would suffice to have $x_1 \geq e^{C/|t|}$ at this stage) and $|\sigma| \leq 1/\log x_n$, we have (as in the proofs of Lemmas 4 and 5) that
$$ \sum_{x_{j}^{1/e} < p \leq x_{j}} \frac{\cos(t\log p)}{p^{1 + 2\sigma}} \ll \frac{1}{|t| \log x_j} \ll \frac{1}{e^j |t| \log x_1} \ll \frac{1}{C e^j} , $$
as well as $\sum_{x_{j}^{1/e} < p \leq x_j} \frac{1}{p^{1+2\sigma}} = 1 + O(\frac{1}{\log x_j} + \frac{\log x_j}{\log x_n}) = 1 + O(C^{-1} e^{-j} + e^{-(n-j)})$. Secondly, by Lemma 6 the ``error'' terms $T(u,v)$ in our characteristic functions satisfy $ |T(u,v) - T(0,0)| \ll \frac{1 + |u|^3 + |v|^3}{\sqrt{x_{j}^{1/e}} \log x_j}$, as well as $ |T(u,0) - T(0,0)| \ll \frac{|u| + |u|^3}{\sqrt{x_{j}^{1/e}} \log x_j}$ and $|T(0,v) - T(0,0)| \ll \frac{|v| + |v|^3}{\sqrt{x_{j}^{1/e}} \log x_j}$, and
$$ |T(u,v) - T(u,0) - T(0,v) + T(0,0)| = |\int_{0}^{u} \int_{0}^{v} \frac{\partial T(u,v)}{\partial u \partial v} du dv| \ll \frac{|u| |v| (1 + |u| + |v|)}{\sqrt{x_{j}^{1/e}} \log x_j} . $$

Now by the two dimensional Berry--Esseen inequality (see Sadikova's paper~\cite{sadikova}), we have
\begin{eqnarray}
&& \Biggl| \tilde{\p}_{t}^{\text{St}, (2)}(u_j \leq \log|I_{l_j}(\frac{1}{2} + \sigma)| \leq u_j + \frac{1}{j^2} , \; \text{and} \; v_j \leq \log|I_{l_j}(\frac{1}{2} + \sigma + it)| \leq v_j + \frac{1}{j^2}) - \nonumber \\
&& - \p(u_j \leq N_{j}^{1} \leq u_j + \frac{1}{j^2} , \; \text{and} \; v_j \leq N_{j}^{2} \leq v_j + \frac{1}{j^2}) \Biggr| \nonumber \\
& \ll & \int_{-x_{j}^{1/50}}^{x_{j}^{1/50}} \int_{-x_{j}^{1/50}}^{x_{j}^{1/50}} \Biggl|\frac{\Delta(u,v)}{uv}\Biggr| du dv + \int_{-x_{j}^{1/50}}^{x_{j}^{1/50}} \Biggl|\frac{\tilde{\E}_{t}^{\text{St}, (2)} e^{iu\log|I_{l_j}(1/2 + \sigma)|} - \E e^{iu N_{j}^{1}} }{u}\Biggr| du + \nonumber \\
&& + \int_{-x_{j}^{1/50}}^{x_{j}^{1/50}} \Biggl|\frac{\tilde{\E}_{t}^{\text{St}, (2)} e^{iv\log|I_{l_j}(1/2 + \sigma + it)|} - \E e^{iv N_{j}^{2}} }{v}\Biggr| dv + \frac{1}{x_{j}^{1/50}} , \nonumber
\end{eqnarray}
where
\begin{eqnarray}
\Delta(u,v) & := & \tilde{\E}_{t}^{\text{St}, (2)} e^{iu\log|I_{l_j}(1/2 + \sigma)| + iv\log|I_{l_j}(1/2 + \sigma + it)|} - \E e^{iuN_{j}^{1} + ivN_{j}^{2}} \nonumber \\
&& - \tilde{\E}_{t}^{\text{St}, (2)} e^{iu\log|I_{l_j}(1/2 + \sigma)|} \tilde{\E}_{t}^{\text{St}, (2)} e^{iv\log|I_{l_j}(1/2 + \sigma + it)|} + \E e^{iuN_{j}^{1}} \E e^{ivN_{j}^{2}} . \nonumber
\end{eqnarray}
Our expression for the characteristic function implies that $\tilde{\E}_{t}^{\text{St}, (2)} e^{iu\log|I_{l_j}(1/2 + \sigma)|}$ is equal to $e^{T(u,0) - T(0,0)} \E e^{iu N_{j}^{1}}$, so the second integral here is
$$ \leq \int_{-x_{j}^{1/50}}^{x_{j}^{1/50}} \Biggl|\frac{e^{T(u,0) - T(0,0)} - 1}{u}\Biggr| du \ll \frac{1}{\sqrt{x_j^{1/e}} \log x_j} \int_{-x_{j}^{1/50}}^{x_{j}^{1/50}} (1+u^2) du \ll \frac{1}{x_{j}^{1/50}} , $$
and the third integral may be bounded exactly similarly. To bound the first (double) integral, we split into the ranges $|uv| \leq 1$ and $|uv| > 1$. On the latter range, to bound $\Delta(u,v)$ we can simply use that
\begin{eqnarray}
\tilde{\E}_{t}^{\text{St}, (2)} e^{iu\log|I_{l_j}(1/2 + \sigma)| + iv\log|I_{l_j}(1/2 + \sigma + it)|} - \E e^{iuN_{j}^{1} + ivN_{j}^{2}} & = & \E e^{iuN_{j}^{1} + ivN_{j}^{2}} (e^{T(u,v) - T(0,0)} - 1) \nonumber \\
& \ll & \frac{1 + |u|^3 + |v|^3}{\sqrt{x_{j}^{1/e}} \log x_j} , \nonumber
\end{eqnarray}
and similarly that
\begin{eqnarray}
&& \tilde{\E}_{t}^{\text{St}, (2)} e^{iu\log|I_{l_j}(1/2 + \sigma)|} \tilde{\E}_{t}^{\text{St}, (2)} e^{iv\log|I_{l_j}(1/2 + \sigma + it)|} - \E e^{iuN_{j}^{1}} \E e^{ivN_{j}^{2}} \nonumber \\
& = & \E e^{iuN_{j}^{1}} \E e^{ivN_{j}^{2}} (e^{T(u,0) - T(0,0) + T(0,v) - T(0,0)} - 1) \ll \frac{1 + |u|^3 + |v|^3}{\sqrt{x_{j}^{1/e}} \log x_j} . \nonumber
\end{eqnarray}
To handle the more delicate range where $|uv| \leq 1$, we note that
\begin{eqnarray}
\Delta(u,v) & = & \E e^{iuN_{j}^{1}} \E e^{ivN_{j}^{2}} \Biggl(e^{-uv \sum_{x_{j}^{1/e} < p \leq x_{j}} \frac{\cos(t\log p)}{2p^{1 + 2\sigma}} + T(u,v) - T(0,0)} - e^{-uv \sum_{x_{j}^{1/e} < p \leq x_{j}} \frac{\cos(t\log p)}{2p^{1 + 2\sigma}}} \nonumber \\
&& - e^{T(u,0) - T(0,0) + T(0,v) - T(0,0)} + 1 \Biggr) . \nonumber
\end{eqnarray}
We can write $e^{T(u,0) - T(0,0) + T(0,v) - T(0,0)} = e^{T(u,v) - T(0,0)} e^{-(T(u,v) - T(u,0) - T(0,v) + T(0,0))}$, which is $= e^{T(u,v) - T(0,0)} + O(\frac{|u| |v| (1 + |u| + |v|)}{\sqrt{x_{j}^{1/e}} \log x_j})$, using our earlier estimation of $T(u,v) - T(u,0) - T(0,v) + T(0,0)$. So we have that $|\Delta(u,v)|$ has order at most
$$ |e^{T(u,v) - T(0,0)} - 1||e^{-uv \sum_{x_{j}^{1/e} < p \leq x_{j}} \frac{\cos(t\log p)}{2p^{1 + 2\sigma}}} - 1| + \frac{|u| |v| (1 + |u| + |v|)}{\sqrt{x_{j}^{1/e}} \log x_j} \ll \frac{|u| |v| (1 + |u|^3 + |v|^3)}{\sqrt{x_{j}^{1/e}} \log x_j} , $$
and obtain overall that $\int_{-x_{j}^{1/50}}^{x_{j}^{1/50}} \int_{-x_{j}^{1/50}}^{x_{j}^{1/50}} |\frac{\Delta(u,v)}{uv}| du dv \ll \frac{x_{j}^{1/10}}{\sqrt{x_{j}^{1/e}} \log x_j} \ll \frac{1}{x_{j}^{1/20}} \ll \frac{1}{x_{j}^{1/50}}$.

To finish the proof, let us (for concision) temporarily set $\mu_{j} = \sum_{x_{j}^{1/e} < p \leq x_{j}} \frac{1 + \cos(t\log p)}{p^{1 + 2\sigma}}$, and $\sigma_j := \sqrt{\sum_{x_{j}^{1/e} < p \leq x_{j}} \frac{1}{2p^{1 + 2\sigma}}}$, and $r_j = \frac{\sum_{x_{j}^{1/e} < p \leq x_{j}} \frac{\cos(t\log p)}{2p^{1 + 2\sigma}}}{\sum_{x_{j}^{1/e} < p \leq x_{j}} \frac{1}{2p^{1 + 2\sigma}}}$, and recall we calculated earlier that $\sigma_j \asymp 1$ and $\mu_j \asymp 1$ and $r_j \ll 1/(|t| \log x_j) \ll 1/C$. Then we may note that
\begin{eqnarray}
&& \p(u_j \leq N_{j}^{1} \leq u_j + \frac{1}{j^2} , \; \text{and} \; v_j \leq N_{j}^{2} \leq v_j + \frac{1}{j^2}) \nonumber \\
& = & \p(\frac{u_j - \mu_j}{\sigma_j} \leq \frac{N_{j}^{1} - \mu_{j}}{\sigma_j} \leq \frac{u_j + 1/j^2 - \mu_j}{\sigma_j} , \; \text{and} \; \frac{v_j - \mu_j}{\sigma_j} \leq \frac{N_{j}^{2} - \mu_j}{\sigma_j} \leq \frac{v_j + 1/j^2 - \mu_j}{\sigma_j}) \nonumber \\
& = & \int_{(u_j - \mu_j)/\sigma_j}^{(u_j - \mu_j + 1/j^2)/\sigma_j} \int_{(v_j - \mu_j)/\sigma_j}^{(v_j - \mu_j + 1/j^2)/\sigma_j} \frac{1}{2\pi \sqrt{1 - r_j^2}} e^{-(x^2 - 2r_{j}xy + y^2)/2(1-r_j^2)} dx dy . \nonumber 
\end{eqnarray}
In particular, since we have $|u_j|, |v_j| \leq (\log x_j)^{1/4}$ we see this probability is $\gg \frac{1}{j^4} e^{-O(\sqrt{\log x_j})}$. Thus, as in the one dimensional case in Lemma 3, we may replace the absolute error term $1/x_{j}^{1/50}$ from the Berry--Esseen theorem by a multiplier $1 + O(1/x_j^{1/100})$, say. To obtain the second estimate in Lemma 7, we simply note that
$$ \frac{e^{-(x^2 - 2r_{j}xy + y^2)/2(1-r_j^2)}}{2\pi \sqrt{1 - r_j^2}} = \frac{e^{-(x^2 + y^2)/2}}{2\pi \sqrt{1 - r_j^2}} e^{O(r_j \sqrt{\log x_j})} = \frac{e^{-(x^2 + y^2)/2}}{2\pi}(1 + O(\frac{1}{|t| \sqrt{\log x_j}})) , $$
and here $\frac{1}{|t| \sqrt{\log x_j}} \ll \frac{1}{|t| \sqrt{e^j \log x_1}} \ll \frac{1}{\sqrt{C e^j}}$ because of our assumption that $x_1 \geq e^{C/|t|^2}$.
\end{proof}

Now we can swiftly deduce the following result, which is what we shall need for our lower bound proof and is a two dimensional analogue of Proposition 5.

\begin{prop7}
There is a large natural number $B$ such that the following is true.

Let $t \in \R$ satisfy $|t| \leq 1$, and let $D \geq 2\log(1/|t|) + (B+1)$ be any natural number. Let $n \leq \log\log x - D$ be large, and define the decreasing sequence $(l_j)_{j=1}^{n}$ of non-negative integers by $l_j := \lfloor \log\log x \rfloor - D - j$. Suppose also that $|\sigma| \leq \frac{1}{e^{D+n}}$.

Then uniformly for any large $a$ and any function $h(n)$ satisfying $|h(n)| \leq 10\log n$, and with $I_{l}(s)$ denoting the increments of the Euler product corresponding to a Steinhaus random multiplicative function (as before), we have
$$ \tilde{\p}_{t}^{\text{St}, (2)}(-a -Bj \leq \sum_{m=1}^{j} \log|I_{l_m}(\frac{1}{2} + \sigma)| , \sum_{m=1}^{j} \log|I_{l_m}(\frac{1}{2} + \sigma + it)| \leq a + j + h(j) \; \forall j \leq n) \asymp \min\{1,\frac{a}{\sqrt{n}}\}^2 . $$
\end{prop7}

\begin{proof}[Proof of Proposition 7]
Using Lemma 7, together with a slicing argument exactly as in the deduction of Lemma 4, we may show the following: under the hypotheses of Lemma 7, but now assuming that $-(1/80)(\log x_j)^{1/4} \leq u_j \leq v_j \leq (1/80)(\log x_j)^{1/4}$ for all $1 \leq j \leq n$, we have
\begin{eqnarray}
&& \tilde{\p}_{t}^{\text{St}, (2)}(u_j \leq \sum_{m=1}^{j} \log|I_{l_m}(\frac{1}{2} + \sigma)| , \sum_{m=1}^{j} \log|I_{l_m}(\frac{1}{2} + \sigma + it)|  \leq v_j ,  \; \forall 1 \leq j \leq n) \nonumber \\
& \asymp & \p(u_j - j + O(1) \leq \sum_{m=1}^{j} G_m \leq v_j - j + O(1) \; \forall 1 \leq j \leq n)^2 , \nonumber
\end{eqnarray}
where $(G_m)_{m=1}^{j}$ are a sequence of independent Gaussian random variables, each having mean zero and variance $\sum_{x_{m}^{1/e} < p \leq x_m} \frac{1}{2p^{1 + 2\sigma}}$.

Now under the hypotheses of Proposition 7, we have $\log x_1 = e^{-(l_1 + 1)} \log x \geq \frac{e^{B+1}}{|t|^2}$ and $|\sigma| \leq \frac{1}{e^{D+n}} \leq \frac{1}{\log x_n}$. This means the hypotheses of Lemma 7 will be satisfied, provided the constant $B$ is fixed sufficiently large, so we may use the probability estimate in the previous paragraph. The conclusion of Proposition 7 now follows exactly as in the one dimensional case in Proposition 5.
\end{proof}

\subsection{The lower bound in the Steinhaus case}
Let $B$ be the large fixed natural number from Proposition 5 (which we may assume, without loss of generality, is the same as the number $B$ from Proposition 7). For each $t \in \R$, let $L(t) = L_{x,q,V}(t)$ denote the event that
\begin{eqnarray}
\Biggl(\frac{\log x}{e^{j+1}} \Biggr)^{-B} e^{-\min\{\sqrt{\log\log x}, 1/(1-q)\}} & \leq & \prod_{l = j}^{\lfloor \log\log x \rfloor - B - 2} |I_{l}(1/2 + \frac{4V}{\log x} + it)| \nonumber \\
& \leq & \frac{\log x}{e^{j+1}} e^{\min\{\sqrt{\log\log x}, 1/(1-q)\} - 2\log\log(\frac{\log x}{e^{j+1}})} \nonumber
\end{eqnarray}
for all $\lfloor \log V \rfloor + 3 \leq j \leq \log\log x - B - 2$. Here the quantity $V$ will ultimately be fixed as another large constant, but initially we allow any $1 \leq V \leq (\log x)^{1/100}$, say. Let $\mathcal{L}$ denote the random subset of points $|t| \leq 1/2$ at which $L(t)$ occurs.

The following is the key restricted second moment estimate that we shall need.
\begin{keyprop5}
With the foregoing notation, and uniformly for all large $x$ and $2/3 \leq q \leq 1$ and $1 \leq V \leq (\log x)^{1/100}$, we have
$$ \E \Biggl(\int_{\mathcal{L}} |F(1/2 + \frac{4V}{\log x} + it)|^2 dt \Biggr)^{2} \ll e^{2\min\{\sqrt{\log\log x}, 1/(1-q)\}} \Biggl(\frac{\log x}{V(1 + (1-q)\sqrt{\log\log x})} \Biggr)^2 . $$
\end{keyprop5}

We remark that the factor $e^{2\min\{\sqrt{\log\log x}, 1/(1-q)\}}$ here, which may look rather alarming, will in fact be harmless because, when raised to the power $1-q$, it becomes $\ll 1$. This should become clear imminently, when we deduce the Theorem 1 lower bound.

\begin{proof}[Proof of the lower bound in Theorem 1, assuming Key Proposition 5]
As argued at the beginning of section \ref{secmainlower}, H\"{o}lder's inequality implies that
$$ \E \Biggl( \int_{\mathcal{L}} |F(1/2 + \frac{4V}{\log x} + it)|^2 dt \Biggr)^{q} \geq \frac{\Biggl( \E \int_{\mathcal{L}} |F(1/2 + \frac{4V}{\log x} + it)|^2 dt \Biggr)^{2-q}}{\Biggl( \E \Biggl(\int_{\mathcal{L}} |F(1/2 + \frac{4V}{\log x} + it)|^2 dt \Biggr)^{2} \Biggr)^{1-q}} . $$

In the Steinhaus case, translation invariance in law implies that we have a simplified expression for the numerator, namely
$$ \Biggl( \int_{-1/2}^{1/2} \E \textbf{1}_{L(t)} |F(1/2 + \frac{4V}{\log x} + it)|^2 dt \Biggr)^{2-q} = \Biggl( \E \textbf{1}_{L(0)} |F(1/2 + \frac{4V}{\log x})|^2 \Biggr)^{2-q} . $$
We can apply Proposition 5 from section \ref{subsecstgirsanov} here, taking $n = \lfloor \log\log x \rfloor - (B+1) - (\lfloor \log V \rfloor + 3)$ and $t_j \equiv 0$, and $a = \min\{\sqrt{\log\log x}, 1/(1-q)\} + O(1)$  and $h(n) = -2\log n$. Using also \eqref{stprodusual} and the standard Mertens prime number estimates, which imply that $\E |F(1/2 + \frac{4V}{\log x})|^2$ is $\gg \exp\{\sum_{p \leq x} \frac{1}{p^{1 + 8V/\log x}} \} = \exp\{\sum_{p \leq x^{1/V}} \frac{1}{p^{1 + 8V/\log x}} + O(1) \} \gg (\log x)/V$, we deduce our numerator is
$$ \gg \Biggl( \frac{1}{1 + (1-q)\sqrt{\log\log x}} \E |F(1/2 + \frac{4V}{\log x})|^2 \Biggr)^{2-q} \gg \Biggl( \frac{\log x}{V(1 + (1-q)\sqrt{\log\log x})} \Biggr)^{2-q} . $$
Inserting the estimate from Key Proposition 5 to upper bound the denominator, and taking $2q$-th roots, we deduce overall that
$$ \vert\vert \int_{-1/2}^{1/2} |F(1/2 + \frac{4V}{\log x} + it)|^2 dt \vert\vert_{q}^{1/2} \gg \sqrt{\frac{\log x}{V(1 + (1-q)\sqrt{\log\log x})}} . $$

Meanwhile, if we argue as we did in section \ref{subsecsteinupper} when proving the upper bound part of Theorem 1 (specifically when handling $F_{k}(s)$ with $k = \lfloor \log V \rfloor$), we have
$$ \vert\vert \int_{-1/2}^{1/2} |F(1/2 + \frac{2V}{\log x} + it)|^2 dt \vert\vert_{q}^{1/2} \ll \sqrt{\frac{\log x}{V(1 + (1-q)\sqrt{\log\log x})}} . $$
Substituting these two bounds into Proposition 3, and choosing $V$ to be a sufficiently large fixed constant that the term $C/e^V$ there kills off the effect of the implicit constants, we obtain
$$ || \sum_{n \leq x} f(n) ||_{2q} \gg \sqrt{\frac{x}{V(1 + (1-q)\sqrt{\log\log x})}} \gg \sqrt{\frac{x}{1 + (1-q)\sqrt{\log\log x}}} . $$
This is the lower bound claimed in Theorem 1.
\end{proof}

\begin{proof}[Proof of Key Proposition 5]
Expanding out and recalling the definition of $\mathcal{L}$, we find $\E (\int_{\mathcal{L}} |F(1/2 + \frac{4V}{\log x} + it)|^2 dt )^{2}$ is
$$ = \int_{-1/2}^{1/2} \int_{-1/2}^{1/2} \E \textbf{1}_{L(t)} |F(1/2 + \frac{4V}{\log x} + it)|^2 \textbf{1}_{L(s)} |F(1/2 + \frac{4V}{\log x} + is)|^2 ds dt . $$
In the Steinhaus case, we can use translation invariance in law to simplify this by shifting $t$ to 0, and replacing $s$ by $s-t$. This yields that
$$ \E \Biggl(\int_{\mathcal{L}} |F(1/2 + \frac{4V}{\log x} + it)|^2 dt \Biggr)^{2} \leq \int_{-1}^{1} \E \textbf{1}_{L(0)} |F(1/2 + \frac{4V}{\log x})|^2 \textbf{1}_{L(t)} |F(1/2 + \frac{4V}{\log x} + it)|^2 dt . $$

Now the point here is that given a shift of size $t$, we expect the parts of the Euler products on primes roughly $\leq e^{1/|t|}$ (or on primes $\leq x$, if $|t| < 1/\log x$) to behave in almost the same way, and the parts on larger primes to behave almost independently. So our strategy is to use the fact that, by definition, when the event $L(t)$ occurs we have
\begin{equation}\label{ltoccurs}
\prod_{\substack{l = \max\{\lfloor \log V \rfloor + 3, \\ \lfloor \log(|t|\log x) \rfloor\}}}^{\lfloor \log\log x \rfloor - B - 2} |I_{l}(\frac{1}{2} + \frac{4V}{\log x} + it)|^2 \ll \Biggl( \min\{\frac{\log x}{V (\log\log x)^2} , \frac{1}{|t|\log^2(2/|t|)}\} e^{\min\{\sqrt{\log\log x}, \frac{1}{1-q}\}} \Biggr)^2 ,
\end{equation}
and then use our probabilistic estimates to bound the remaining part
\begin{equation}\label{remintegrand}
\E \textbf{1}_{L(0)} |F(1/2 + \frac{4V}{\log x})|^2 \textbf{1}_{L(t)} \frac{|F(1/2 + \frac{4V}{\log x} + it)|^2}{\prod_{\substack{l = \max\{\lfloor \log V \rfloor + 3, \\ \lfloor \log(|t|\log x) \rfloor\}}}^{\lfloor \log\log x \rfloor - B - 2} |I_{l}(1/2 + \frac{4V}{\log x} + it)|^2}
\end{equation}
of the integrand.

It turns out that when $|t| \leq 1/(\log x)^{1/3}$, say, we can afford to take a crude approach thanks to the saving $1/(\log\log x)^4$ in \eqref{ltoccurs}, and now throw away the indicator functions $\textbf{1}_{L(0)}, \textbf{1}_{L(t)}$ from \eqref{remintegrand}. Having done this, using the independence of $f(p)$ on distinct primes we find \eqref{remintegrand} is
$$ \leq \E \prod_{\substack{l = \max\{\lfloor \log V \rfloor + 3, \\ \lfloor \log(|t|\log x) \rfloor\}}}^{\lfloor \log\log x \rfloor - B - 2} |I_{l}(1/2 + \frac{4V}{\log x})|^2 \E \frac{|F(1/2 + \frac{4V}{\log x})|^2 |F(1/2 + \frac{4V}{\log x} + it)|^2}{\prod_{\substack{l = \max\{\lfloor \log V \rfloor + 3, \\ \lfloor \log(|t|\log x) \rfloor\}}}^{\lfloor \log\log x \rfloor - B - 2} |I_{l}(1/2 + \frac{4V}{\log x})|^2 |I_{l}(1/2 + \frac{4V}{\log x} + it)|^2} . $$
Now using \eqref{sttwoprods}, the second expectation here is equal to
$$ \exp\{\sum_{p \leq x^{e^{-(\lfloor \log\log x \rfloor - B)}}} \frac{2 + 2\cos(t\log p)}{p^{1+8V/\log x}} + \sum_{x^{e^{-(\max\{\lfloor \log V \rfloor + 3, \lfloor \log(|t|\log x) \rfloor\} + 1)}} < p \leq x} \frac{2 + 2\cos(t\log p)}{p^{1+8V/\log x}} + O(1) \} , $$
and remembering that $B$ is a fixed constant, and using standard Mertens and Chebychev type estimates for sums over primes, this is
$$ \ll \exp\{\sum_{\min\{x^{1/V}, e^{1/|t|}\} < p \leq x} \frac{2 + 2\cos(t\log p)}{p^{1+8V/\log x}}\} \ll \exp\{\sum_{\min\{x^{1/V}, e^{1/|t|}\} < p \leq x^{1/V}} \frac{2 + 2\cos(t\log p)}{p}\} . $$
Again, since we may assume here that $p > e^{1/|t|}$, standard estimates for sums over primes (as in the proof of Lemma 5) show the overall contribution from the $\cos(t\log p)$ sum is $\ll 1$, so we finally have a bound $\ll \exp\{\sum_{\min\{x^{1/V}, e^{1/|t|}\} < p \leq x^{1/V}} \frac{2}{p}\} \ll (\max\{1, \frac{|t| \log x}{V}\})^2$. Meanwhile, \eqref{stprodusual} implies that the first expectation $\E \prod_{\substack{l = \max\{\lfloor \log V \rfloor + 3, \\ \lfloor \log(|t|\log x) \rfloor\}}}^{\lfloor \log\log x \rfloor - B - 2} |I_{l}(1/2 + \frac{4V}{\log x})|^2$ is equal to $\exp\{\sum_{p \leq \min\{x^{1/V}, e^{1/|t|}\}} \frac{1}{p^{1 + 8V/\log x}} + O(1) \} \ll \min\{\frac{\log x}{V}, \frac{1}{|t|} \}$. So putting together \eqref{ltoccurs} with our above upper bounds for \eqref{remintegrand}, we obtain
\begin{eqnarray}
&& \int_{|t| \leq 1/(\log x)^{1/3}} \E \textbf{1}_{L(0)} |F(1/2 + \frac{4V}{\log x})|^2 \textbf{1}_{L(t)} |F(1/2 + \frac{4V}{\log x} + it)|^2 dt \nonumber \\
& \ll & \frac{e^{2\min\{\sqrt{\log\log x}, 1/(1-q)\}}}{(\log\log x)^4} \int_{|t| \leq 1/(\log x)^{1/3}} \Biggl( \min\{\frac{\log x}{V} , \frac{1}{|t|}\} \Biggr)^3 \Biggl(\max\{1, \frac{|t| \log x}{V}\} \Biggr)^2 dt \nonumber \\
& \ll & \frac{e^{2\min\{\sqrt{\log\log x}, \frac{1}{1-q}\}} \log^{2}x}{V^2 (\log\log x)^4}  \int_{|t| \leq 1/(\log x)^{1/3}} \min\{\frac{\log x}{V} , \frac{1}{|t|}\} dt \ll \frac{e^{2\min\{\sqrt{\log\log x}, \frac{1}{1-q}\}} \log^{2}x}{V^2 (\log\log x)^3} , \nonumber
\end{eqnarray}
which is more than acceptable for Key Proposition 5.

To handle the contribution from $1/(\log x)^{1/3} < |t| \leq 1$, we follow the same approach, but need to deal with the indicator functions $\textbf{1}_{L(0)}, \textbf{1}_{L(t)}$ in \eqref{remintegrand} more carefully because our saving $1/\log^{4}(2/|t|)$ from \eqref{ltoccurs} is no longer so great. We temporarily set $D = D(t) = \lceil 2\log(1/|t|) \rceil + (B+1)$. If the event $L(0)$ occurs then, comparing the definition of $L(0)$ with $j = \lfloor \log\log x \rfloor - D$ and with general $j$ (and recalling that $B$ is an absolute constant), we must in particular have
\begin{eqnarray}
\Biggl(\frac{\log x}{e^{j+1}} \Biggr)^{-B} e^{-2\min\{\sqrt{\log\log x}, \frac{1}{1-q}\}} |t|^2 \log^2(\frac{2}{|t|}) & \ll & \prod_{l = j}^{\lfloor \log\log x \rfloor - D - 1} |I_{l}(\frac{1}{2} + \frac{4V}{\log x})| \nonumber \\
& \ll & \frac{\log x}{e^{j+1}} e^{- 2\log\log(\frac{\log x}{e^{j+1}})} \frac{e^{2\min\{\sqrt{\log\log x}, 1/(1-q)\}}}{|t|^{2B}} \nonumber
\end{eqnarray}
for all $\lfloor \log V \rfloor + 3 \leq j \leq \log\log x - D - 1$. Similarly, if the event $L(t)$ occurs then we must have the same bounds with $|I_{l}(1/2 + \frac{4V}{\log x})|$ replaced by $|I_{l}(1/2 + \frac{4V}{\log x} + it)|$. Let us write $R(t)$ for the event that one has these bounds for both $|I_{l}(1/2 + \frac{4V}{\log x})|$ and $|I_{l}(1/2 + \frac{4V}{\log x} + it)|$. Then we can upper bound \eqref{remintegrand} by
$$ \E \textbf{1}_{R(t)} |F(1/2 + \frac{4V}{\log x})|^2 \frac{|F(1/2 + \frac{4V}{\log x} + it)|^2}{\prod_{l = \lfloor \log(|t|\log x) \rfloor}^{\lfloor \log\log x \rfloor - B - 2} |I_{l}(1/2 + \frac{4V}{\log x} + it)|^2} . $$
Now here we have $\log\log x - D - 1 < \lfloor \log(|t|\log x) \rfloor$. Since $f(p)$ is independent on distinct primes, and so the event $R(t)$ is independent of $\prod_{l = \lfloor \log(|t|\log x) \rfloor}^{\lfloor \log\log x \rfloor - B - 2} |I_{l}(1/2 + \frac{4V}{\log x})|^2$, we may pull $\E \prod_{l = \lfloor \log(|t|\log x) \rfloor}^{\lfloor \log\log x \rfloor - B - 2} |I_{l}(1/2 + \frac{4V}{\log x})|^2$ out from the above and use \eqref{stprodusual} to show it is $= \exp\{\sum_{p \leq e^{1/|t|}} \frac{1}{p^{1 + 8V/\log x}} + O(1) \} \ll 1/|t|$. So employing our notation from section \ref{subsecfurtherprobcalc}, we have that \eqref{remintegrand} is
\begin{eqnarray}
& \ll & (1/|t|) \E \textbf{1}_{R(t)} \frac{|F(1/2 + \frac{4V}{\log x})|^2}{\prod_{l = \lfloor \log(|t|\log x) \rfloor}^{\lfloor \log\log x \rfloor - B - 2} |I_{l}(1/2 + \frac{4V}{\log x})|^2} \frac{|F(1/2 + \frac{4V}{\log x} + it)|^2}{\prod_{l = \lfloor \log(|t|\log x) \rfloor}^{\lfloor \log\log x \rfloor - B - 2} |I_{l}(1/2 + \frac{4V}{\log x} + it)|^2} \nonumber \\
& = & \frac{1}{|t|} \tilde{\p}_{t}^{\text{St}, (2)}(R(t)) \E \frac{|F(1/2 + \frac{4V}{\log x})|^2}{\prod_{l = \lfloor \log(|t|\log x) \rfloor}^{\lfloor \log\log x \rfloor - B - 2} |I_{l}(1/2 + \frac{4V}{\log x})|^2} \frac{|F(1/2 + \frac{4V}{\log x} + it)|^2}{\prod_{l = \lfloor \log(|t|\log x) \rfloor}^{\lfloor \log\log x \rfloor - B - 2} |I_{l}(1/2 + \frac{4V}{\log x} + it)|^2} . \nonumber
\end{eqnarray}
We may deploy Proposition 7 to bound $\tilde{\p}_{t}^{\text{St}, (2)}(R(t))$, taking $n = \lfloor \log\log x \rfloor - D - (\lfloor \log V \rfloor + 3) \gg \log\log x$ and $a = 2\min\{\sqrt{\log\log x}, 1/(1-q)\} + O(\log(2/|t|))$ and $h(n) = -2\log n$, and obtaining that
$$ \tilde{\p}_{t}^{\text{St}, (2)}(R(t)) \ll \Biggl(\frac{1}{1 + (1-q)\sqrt{\log\log x}} + \frac{\log(2/|t|)}{\sqrt{\log\log x}} \Biggr)^2 \ll \frac{\log^{2}(2/|t|)}{(1 + (1-q)\sqrt{\log\log x})^2} . $$
Furthermore, as we did earlier we may use \eqref{sttwoprods} to obtain that the expectation in the preceding display is $\ll (\frac{|t| \log x}{V} )^2$.

Finally, putting together \eqref{ltoccurs} with the above upper bounds for \eqref{remintegrand} yields that
\begin{eqnarray}
&& \int_{1/(\log x)^{1/3} < |t| \leq 1} \E \textbf{1}_{L(0)} |F(1/2 + \frac{4V}{\log x})|^2 \textbf{1}_{L(t)} |F(1/2 + \frac{4V}{\log x} + it)|^2 dt \nonumber \\
& \ll & \frac{e^{2\min\{\sqrt{\log\log x}, 1/(1-q)\}}}{(1 + (1-q)\sqrt{\log\log x})^2} \int_{1/(\log x)^{1/3} < |t| \leq 1} \Biggl( \frac{1}{|t| \log^2(2/|t|)} \Biggr)^2 \frac{\log^{2}(2/|t|)}{|t|} \Biggl(\frac{|t| \log x}{V} \Biggr)^2 dt \nonumber \\
& \ll & \frac{e^{2\min\{\sqrt{\log\log x}, 1/(1-q)\}} \log^{2}x}{V^2 (1 + (1-q)\sqrt{\log\log x})^2}  \int_{1/(\log x)^{1/3} < |t| \leq 1} \frac{1}{|t| \log^{2}(2/|t|)} dt \ll \frac{e^{2\min\{\sqrt{\log\log x}, \frac{1}{1-q}\}} \log^{2}x}{V^2 \log\log x} . \nonumber
\end{eqnarray}
Because of the factor $\log^{2}(2/|t|)$ (which we remark descends from the term $-2\log\log(\frac{\log x}{e^{j+1}})$ in the original definition of the set $\mathcal{L}$), the integral here is $\ll 1$ and we obtain the upper bound claimed in Key Proposition 5.
\end{proof}

\subsection{The lower bound in the Rademacher case}
As in our upper bound arguments, there is not too much difference between the Rademacher and Steinhaus cases, except the former is more notationally complicated because we no longer have ``translation invariance in law'' of the Euler products. We will only sketch the changes that are needed to deduce the Rademacher lower bound.

In place of Lemma 6, in the Rademacher case one needs to compute
$$ \E \prod_{x < p \leq y} \left|1 + \frac{f(p)}{p^{1/2+\sigma+it_1}}\right|^{2+iu} \left|1 + \frac{f(p)}{p^{1/2+\sigma+i(t_1 + t_2)}}\right|^{2+iv} , $$
where $t_1 , t_2, u, v$ are real. If one does this (proceeding as in the proofs of Lemmas 2 and 6), one obtains an expression that is, up to error terms, the exponential of
\begin{eqnarray}
&& \sum_{x < p \leq y} \Biggl( \frac{(1+iu/2)^2 + (1+iv/2)^2}{p^{1+2\sigma}} + \frac{(2+iu)(2+iv)\cos(t_1 \log p) \cos((t_1 + t_2)\log p)}{p^{1+2\sigma}} + \nonumber \\
&& + \frac{((2+iu)^{2}/2 - (2+iu))\cos(2t_1 \log p) + ((2+iv)^{2}/2 - (2+iv))\cos(2(t_1 + t_2)\log p)}{2p^{1+2\sigma}} \Biggr) , \nonumber
\end{eqnarray}
which can be simplified a little by observing it is
\begin{eqnarray}
& = & \sum_{x < p \leq y} \Biggl( \frac{(1+iu/2)^2 + (1+iv/2)^2}{p^{1+2\sigma}} + \frac{(2+iu)(2+iv)\cos(t_2 \log p)}{2p^{1+2\sigma}} + \nonumber \\
&& + \frac{(2+iu)(iu/2)\cos(2t_1 \log p) + (2+iv)(iv/2)\cos(2(t_1 + t_2)\log p)}{2p^{1+2\sigma}} + \nonumber \\
&& + \frac{(2+iu)(2+iv)\cos((2t_1 + t_2)\log p)}{2p^{1+2\sigma}} \Biggr) . \nonumber
\end{eqnarray}
Notice the first two fractions here are precisely analogous to the Steinhaus case in Lemma 6, whereas the second two terms are new.

For each real $t$ we can define the event $L(t) = L_{x,q,V}(t)$ exactly as in the Steinhaus case, but now we let $\mathcal{L}$ be the random subset of points $t \in [1/3,1/2]$ (rather than $[-1/2,1/2]$) at which $L(t)$ occurs. The point of this change is that we will then always have $t_1 , t_1 + t_2 \in [1/3,1/2]$, and therefore $2/3 \leq 2t_1 , 2(t_1 + t_2), 2t_1 + t_2 \leq 1$, in our calculations, which means all the new terms in the Rademacher characteristic functions will only cause negligible changes in mean, variance and covariance as compared with the Steinhaus case (because $\sum_{x < p \leq y} \frac{\cos(2t_1 \log p)}{p^{1+2\sigma}} , \sum_{x < p \leq y} \frac{\cos(2(t_1 + t_2) \log p)}{p^{1+2\sigma}} , \sum_{x < p \leq y} \frac{\cos((2t_1 + t_2) \log p)}{p^{1+2\sigma}}$ all exhibit significant cancellation). Once one is in this situation, the Rademacher lower bound follows simply by imitating the arguments from the Steinhaus case.

\section{Proof of Corollary 2, lower bound}\label{corsec}
In this section we shall prove our last remaining result, that for all $2 \leq \lambda \leq e^{\sqrt{\log\log x}}$ we have
$$ \p(|\sum_{n \leq x} f(n)| \geq \lambda\frac{\sqrt{x}}{(\log\log x)^{1/4}}) \gg \frac{1}{\lambda^2 (\log\log x)^{O(1)}} . $$
The proof is almost the same in the Steinhaus and Rademacher cases, but for definiteness let us first think of $f(n)$ as a Steinhaus random multiplicative function. As before, $F(s)$ will denote the Euler product of $f(n)$ over $x$-smooth numbers. The proof will consist of three main steps: we will prove that it will suffice to establish a suitable lower bound for $\p(\int_{-1/2}^{1/2} |F(1/2 + \frac{\log\log x}{\log x} + it)|^2 dt \geq \beta \log x )$; then we will show it will actually suffice to prove a corresponding lower bound for a discrete maximum of $|F(1/2 + \frac{\log\log x}{\log x} + it)|^2$ (roughly speaking); and finally we will deduce that from some existing results on maxima of random processes. 

Let $\hat{\E}$ denote expectation conditional on the values $(f(p))_{p \leq \sqrt{x}}$, and let $\hat{\p}$ denote the corresponding conditional probability. Thus, if we let $P(n)$ denote the largest prime factor of $n$, then $\hat{\E}|\sum_{\substack{n \leq x, \\ P(n) > \sqrt{x}}} f(n)|$ is a function of the random variables $(f(p))_{p \leq \sqrt{x}}$, and in fact we have
$$ \hat{\E}|\sum_{\substack{n \leq x, \\ P(n) > \sqrt{x}}} f(n)| = \hat{\E}|\sum_{\sqrt{x} < p \leq x} f(p) \sum_{m \leq x/p} f(m)| \asymp \sqrt{\sum_{\sqrt{x} < p \leq x} \left|\sum_{m \leq x/p} f(m) \right|^2} , $$
by Khintchine's inequality (as in the proof of Proposition 3). Let $\mathcal{A}$ denote the event that
\begin{equation}
|\sum_{n \leq x} f(n)| \geq (1/2) \hat{\E}|\sum_{\substack{n \leq x, \\ P(n) > \sqrt{x}}} f(n)| . \nonumber
\end{equation}
We will first prove that we always (i.e. {\em for any realisation of} $(f(p))_{p \leq \sqrt{x}}$) have a uniform lower bound $\hat{\p}(\mathcal{A}) \gg 1$. Having done this, to prove the lower bound in Corollary 2 it will suffice to show that
$$ \p\Biggl(\sqrt{\sum_{\sqrt{x} < p \leq x} \left|\sum_{m \leq x/p} f(m) \right|^2} \geq \lambda\frac{\sqrt{x}}{(\log\log x)^{1/4}}\Biggr) \gg \frac{1}{\lambda^2 (\log\log x)^{O(1)}} \;\;\; \forall \; 2 \leq \lambda \leq e^{\sqrt{\log\log x}} , $$
or equivalently that
\begin{equation}\label{cormainlower}
\p\Biggl(\sum_{\sqrt{x} < p \leq x} \left|\sum_{m \leq x/p} f(m) \right|^2 \geq \beta x \Biggr) \gg \frac{1}{\beta (\log\log x)^{O(1)}} \;\;\; \forall \; 2 \leq \beta \leq e^{2\sqrt{\log\log x}} .
\end{equation}

\vspace{12pt}
To establish that $\hat{\p}(\mathcal{A}) \gg 1$ for any realisation of $(f(p))_{p \leq \sqrt{x}}$, note that (as in the proof of Proposition 3) we have
$$ \sum_{n \leq x} f(n) = \sum_{\substack{n \leq x, \\ P(n) \leq \sqrt{x}}} f(n) + \sum_{\substack{n \leq x, \\ P(n) > \sqrt{x}}} f(n) \stackrel{d}{=} \sum_{\substack{n \leq x, \\ P(n) \leq \sqrt{x}}} f(n) + \epsilon \sum_{\substack{n \leq x, \\ P(n) > \sqrt{x}}} f(n) ,$$
where $\epsilon$ is an auxiliary Rademacher random variable independent of everything else, and $\stackrel{d}{=}$ denotes equality in distribution. By the triangle inequality we have
$$ 2|\sum_{\substack{n \leq x, \\ P(n) > \sqrt{x}}} f(n)| \leq |\sum_{\substack{n \leq x, \\ P(n) \leq \sqrt{x}}} f(n) + \sum_{\substack{n \leq x, \\ P(n) > \sqrt{x}}} f(n)| + |\sum_{\substack{n \leq x, \\ P(n) \leq \sqrt{x}}} f(n) - \sum_{\substack{n \leq x, \\ P(n) > \sqrt{x}}} f(n)| , $$
so for at least one of the values of $\epsilon$ we must have $|\sum_{\substack{n \leq x, \\ P(n) \leq \sqrt{x}}} f(n) + \epsilon \sum_{\substack{n \leq x, \\ P(n) > \sqrt{x}}} f(n)| \geq |\sum_{\substack{n \leq x, \\ P(n) > \sqrt{x}}} f(n)|$. So if we let $\mathcal{B}$ denote the event that
$$ |\sum_{\substack{n \leq x, \\ P(n) > \sqrt{x}}} f(n)| \geq (1/2) \hat{\E}|\sum_{\substack{n \leq x, \\ P(n) > \sqrt{x}}} f(n)| , $$
then we have $\hat{\p}(\mathcal{A}) \geq (1/2)\hat{\p}(\mathcal{B})$. Meanwhile, we can lower bound $\hat{\p}(\mathcal{B})$ in a standard (Paley--Zygmund type) way using the Cauchy--Schwarz inequality, noting that
\begin{eqnarray}
\hat{\p}(\mathcal{B}) \geq \frac{(\hat{\E} \textbf{1}_{\mathcal{B}} |\sum_{\substack{n \leq x, \\ P(n) > \sqrt{x}}} f(n)|)^2}{\hat{\E}|\sum_{\substack{n \leq x, \\ P(n) > \sqrt{x}}} f(n)|^2} & = & \frac{(\hat{\E}|\sum_{\substack{n \leq x, \\ P(n) > \sqrt{x}}} f(n)| - \hat{\E} \textbf{1}_{\mathcal{B} \; \text{fails}} |\sum_{\substack{n \leq x, \\ P(n) > \sqrt{x}}} f(n)|)^2}{\hat{\E}|\sum_{\substack{n \leq x, \\ P(n) > \sqrt{x}}} f(n)|^2} \nonumber \\
& \gg & \frac{(\hat{\E}|\sum_{\substack{n \leq x, \\ P(n) > \sqrt{x}}} f(n)| )^2}{\hat{\E}|\sum_{\substack{n \leq x, \\ P(n) > \sqrt{x}}} f(n)|^2} . \nonumber
\end{eqnarray}
Writing $\sum_{\substack{n \leq x, \\ P(n) > \sqrt{x}}} f(n) = \sum_{\sqrt{x} < p \leq x} f(p) \sum_{m \leq x/p} f(m)$ and applying Khintchine's inequality again, we obtain that the conditional expectation (squared) in the numerator, and the conditional expectation in the denominator, have the same order. So we have shown $\hat{\p}(\mathcal{A}) \geq (1/2)\hat{\p}(\mathcal{B}) \gg 1$, as we wanted.

\vspace{12pt}
To establish \eqref{cormainlower}, we again begin similarly as in the proof of Proposition 3, noting (with $X = e^{\sqrt{\log x}}$) that $\sum_{\sqrt{x} < p \leq x} \left|\sum_{m \leq x/p} f(m) \right|^2$ is
\begin{eqnarray}
\geq \sum_{\sqrt{x} < p \leq x} \frac{\log p}{\log x} \left|\sum_{m \leq x/p} f(m) \right|^2 & = & \frac{1}{\log x} \sum_{\sqrt{x} < p \leq x} \log p \frac{X}{p} \int_{p}^{p(1+1/X)} \left|\sum_{m \leq x/p} f(m) \right|^2 dt \nonumber \\
& \geq & \frac{1}{\log x} \sum_{\sqrt{x} < p \leq x} \log p \frac{X}{p} \int_{p}^{p(1+1/X)} \frac{1}{2} \left|\sum_{m \leq x/t} f(m) \right|^2 dt - \nonumber \\
&& - \frac{1}{\log x} \sum_{\sqrt{x} < p \leq x} \log p \frac{X}{p} \int_{p}^{p(1+1/X)} \left|\sum_{x/t < m \leq x/p} f(m) \right|^2 dt . \nonumber
\end{eqnarray}
The subtracted term here satisfies the first moment bound
\begin{eqnarray}
&& \frac{1}{\log x} \sum_{\sqrt{x} < p \leq x} \log p \frac{X}{p} \int_{p}^{p(1+1/X)} \E\left|\sum_{x/t < m \leq x/p} f(m) \right|^2 dt \nonumber \\
& \ll & \frac{1}{\log x} \sum_{\sqrt{x} < p \leq x} \log p \frac{X}{p} \int_{p}^{p(1+1/X)} (\frac{x}{pX} + 1) dt \ll \frac{x}{\log x} , \nonumber
\end{eqnarray}
so the probability it is larger than $x/\sqrt{\log x}$ is $\ll \frac{1}{\sqrt{\log x}}$. Since $x/\sqrt{\log x}$ is negligible compared with the term $\beta x$ in \eqref{cormainlower}, and $\frac{1}{\sqrt{\log x}}$ is negligible compared with our target lower bound $\frac{1}{\beta (\log\log x)^{O(1)}}$, we may ignore this subtracted term. Meanwhile, still broadly following the proof of Proposition 3, we have that
\begin{eqnarray}
&& \frac{1}{\log x} \sum_{\sqrt{x} < p \leq x} \log p \frac{X}{p} \int_{p}^{p(1+1/X)} \left|\sum_{m \leq x/t} f(m) \right|^2 dt \gg \frac{1}{\log x} \int_{\sqrt{x}}^{x} \left|\sum_{m \leq x/t} f(m) \right|^2 dt \nonumber \\
& = & \frac{x}{\log x} \int_{1}^{\sqrt{x}} \left|\sum_{m \leq z} f(m) \right|^2 \frac{dz}{z^{2}} \geq \frac{x}{\log x} \int_{1}^{\sqrt{x}} \left|\sum_{\substack{m \leq z, \\ m \; \text{is} \; x \; \text{smooth}}} f(m) \right|^2 \frac{dz}{z^{2 + 2\log\log x/\log x}} . \nonumber
\end{eqnarray}
We would like to complete the integral to the range $\int_{1}^{\infty}$, so we can apply Harmonic Analysis Result 1 and have a lower bound $\gg \frac{x}{\log x} \int_{-1/2}^{1/2} |F(1/2 + \frac{\log\log x}{\log x} + it)|^2 dt$, where $F$ is the Euler product of $f(n)$ over $x$-smooth numbers. It turns out that the shift $(\log\log x)/\log x$ that we introduced makes it acceptable to perform this completion, as we have the first moment bound
$$ \frac{x}{\log x} \int_{\sqrt{x}}^{\infty} \E \left|\sum_{\substack{m \leq z, \\ x \; \text{smooth}}} f(m) \right|^2 \frac{dz}{z^{2 + 2\log\log x/\log x}} \leq \frac{x}{\log x} \int_{\sqrt{x}}^{\infty} \frac{dz}{z^{1 + 2\log\log x/\log x}} \ll \frac{x}{\log x \log\log x} , $$
so $\int_{\sqrt{x}}^{\infty}$ is negligible as in our discussion of the subtracted term above.

\vspace{12pt}
Reviewing the arguments in the previous paragraph, we see that to establish \eqref{cormainlower} it will now suffice to show that
\begin{equation}\label{corprodlower}
\p\Biggl(\int_{-1/2}^{1/2} |F(1/2 + \frac{\log\log x}{\log x} + it)|^2 dt \geq \beta \log x \Biggr) \gg \frac{1}{\beta (\log\log x)^{O(1)}} \;\;\; \forall \; 2 \leq \beta \leq e^{2\sqrt{\log\log x}} .
\end{equation}
Now the proof begins to differ from what we have done previously. We would like to say that the integral is essentially the same as $\sum_{|k| \leq \frac{\log x - 1}{2}} \frac{1}{\log x} |F(1/2 + \frac{\log\log x}{\log x} + i\frac{k}{\log x})|^2$, so we could lower bound the probability on the left by the probability that the maximum of the $|F(1/2 + \frac{\log\log x}{\log x} + i\frac{k}{\log x})|^2$ is $\gg \beta \log^{2}x$. It seems technically tricky to establish an approximation quite like that, but since we only want a lower bound and since the exponential function is convex we can exploit Jensen's inequality, obtaining that
\begin{eqnarray}
&& \int_{-1/2}^{1/2} |F(\frac{1}{2} + \frac{\log\log x}{\log x} + it)|^2 \geq \sum_{|k| \leq \frac{\log x - 1}{2}} \int_{-\frac{1}{2\log x}}^{\frac{1}{2\log x}} |F(\frac{1}{2} + \frac{\log\log x}{\log x} + i\frac{k}{\log x} + it)|^2 dt \nonumber \\
& = & \frac{1}{\log x} \sum_{|k| \leq \frac{\log x - 1}{2}} \Biggl( \log x \int_{-\frac{1}{2\log x}}^{\frac{1}{2\log x}} \exp\{2\log|F(1/2 + \frac{\log\log x}{\log x} + i\frac{k}{\log x} + it)|\} dt \Biggr) \nonumber \\
& \geq & \frac{1}{\log x} \sum_{|k| \leq \frac{\log x - 1}{2}} \exp\Biggl\{2\log x \int_{-\frac{1}{2\log x}}^{\frac{1}{2\log x}} \log|F(1/2 + \frac{\log\log x}{\log x} + i\frac{k}{\log x} + it)| dt \Biggr\} . \nonumber
\end{eqnarray}

The exponential in the sum here will behave, for our purposes, in exactly the same way as $|F(1/2 + \frac{\log\log x}{\log x} + i\frac{k}{\log x})|^2$, but we may calculate precisely that it is
\begin{eqnarray}
&& \exp\Biggl\{- 2\log x \int_{-\frac{1}{2\log x}}^{\frac{1}{2\log x}} \sum_{p \leq x} \Re\log(1 - \frac{f(p)}{p^{1/2 + \frac{\log\log x}{\log x} + i\frac{k}{\log x} + it}}) dt \Biggr\} \nonumber \\
& = & \exp\Biggl\{2 \sum_{p \leq x} \log x \int_{-\frac{1}{2\log x}}^{\frac{1}{2\log x}} (\Re \frac{f(p) p^{- i\frac{k}{\log x} - it}}{p^{1/2 + \frac{\log\log x}{\log x}}} + \Re \frac{(f(p) p^{- i\frac{k}{\log x} - it})^2}{2 p^{1 + \frac{2\log\log x}{\log x}}} + O(\frac{1}{p^{3/2}}) ) dt \Biggr\} . \nonumber
\end{eqnarray}
We can simplify this expression a bit further by noting that $\log x \int_{- \frac{1}{2\log x}}^{\frac{1}{2\log x}}e^{-2it\log p} dt = \log x \int_{- \frac{1}{2\log x}}^{\frac{1}{2\log x}} (1 + O(|t|\log p)) dt = 1 + O((\log p)/\log x)$, and that $\sum_{p \leq x} \frac{\log p}{p^{1 + \frac{2\log\log x}{\log x}} \log x} = O(1)$ by Chebychev's estimate, as is $\sum_{p \leq x} \frac{1}{p^{3/2}}$. Thus the exponential is
$$ = \exp\Biggl\{2\Biggl( \sum_{p \leq x} \log x \int_{-\frac{1}{2\log x}}^{\frac{1}{2\log x}} \Re \frac{f(p) p^{- i\frac{k}{\log x} - it}}{p^{1/2 + \frac{\log\log x}{\log x}}} dt + \sum_{p \leq x} \Re \frac{(f(p) p^{- i\frac{k}{\log x}})^2}{2p^{1 + \frac{2\log\log x}{\log x}}} \Biggr) + O(1) \Biggr\} , $$
and to prove \eqref{corprodlower} (and therefore also \eqref{cormainlower}) it will suffice to show that
\begin{eqnarray}
&& \p\Biggl(\max_{|k| \leq \frac{\log x - 1}{2}} \Biggl( \sum_{p \leq x} \log x \int_{-\frac{1}{2\log x}}^{\frac{1}{2\log x}} \Re \frac{f(p) p^{- i\frac{k}{\log x} - it}}{p^{1/2 + \frac{\log\log x}{\log x}}} dt + \sum_{p \leq x} \Re \frac{(f(p) p^{- i\frac{k}{\log x}})^2}{2p^{1 + \frac{2\log\log x}{\log x}}} \Biggr) \geq \frac{\log\beta}{2} + \log\log x \Biggr) \nonumber \\
&& \gg \frac{1}{\beta (\log\log x)^{O(1)}} \;\;\; \forall \; 2 \leq \beta \leq e^{2\sqrt{\log\log x}} . \nonumber
\end{eqnarray}

As a final simplification, let us note that $|\sum_{p \leq \log^{10}x} \frac{(f(p) p^{- i\frac{k}{\log x}})^2}{p^{1 + \frac{2\log\log x}{\log x}}}| \leq \sum_{p \leq \log^{10}x} \frac{1}{p} = \log\log\log x + O(1)$ by the Mertens estimate, and
$$ \E\Biggl| \sum_{\log^{10}x < p \leq x} \frac{(f(p) p^{- i\frac{k}{\log x}})^2}{p^{1 + \frac{2\log\log x}{\log x}}} \Biggr|^2 = \sum_{\log^{10}x < p \leq x} \frac{1}{p^{2 + \frac{4\log\log x}{\log x}}} \ll \frac{1}{\log^{10}x} , $$
so we have $\p(\max_{|k| \leq \frac{\log x - 1}{2}} |\sum_{\log^{10}x < p \leq x} \Re \frac{(f(p) p^{- i\frac{k}{\log x}})^2}{2p^{1 + \frac{2\log\log x}{\log x}}}| \geq 1) \ll 1/\log^{9}x$, by the union bound and Chebychev's inequality. Thus, since $1/\log^{9}x$ is negligible compared with our target probability lower bound $\frac{1}{\beta (\log\log x)^{O(1)}}$; and since replacing $\frac{\log\beta}{2} + \log\log x$ by $\frac{\log\beta}{2} + \log\log\log x + O(1) + \log\log x$ (essentially replacing $\beta$ by $\beta e^{2\log\log\log x + O(1)}$) doesn't change the form of that lower bound; we can omit the second sum $\sum_{p \leq x} \Re \frac{(f(p) p^{- i\frac{k}{\log x}})^2}{2p^{1 + \frac{2\log\log x}{\log x}}}$ entirely when trying to prove the lower bound above. Furthermore, since the contributions from different primes are independent, and for any given $k$ the symmetry of the $f(p)$ implies that $\p(\sum_{p \leq \log^{10}x} \log x \int_{-\frac{1}{2\log x}}^{\frac{1}{2\log x}} \Re \frac{f(p) p^{- i\frac{k}{\log x} - it}}{p^{1/2 + \frac{\log\log x}{\log x}}} dt \geq 0) = 1/2$, we are free to omit the part of the first sum over primes $p \leq \log^{10}x$ as well. So to prove \eqref{corprodlower} it will suffice to show
\begin{eqnarray}\label{corfinallower}
&& \p\Biggl(\max_{|k| \leq \frac{\log x - 1}{2}} \Biggl( \sum_{\log^{10}x < p \leq x} \log x \int_{-\frac{1}{2\log x}}^{\frac{1}{2\log x}} \Re \frac{f(p) p^{- i\frac{k}{\log x} - it}}{p^{1/2 + \frac{\log\log x}{\log x}}} dt \Biggr) \geq \frac{\log\beta}{2} + \log\log x \Biggr) \nonumber \\
&& \gg \frac{1}{\beta (\log\log x)^{O(1)}} \;\;\; \forall \; 2 \leq \beta \leq e^{2\sqrt{\log\log x}} .
\end{eqnarray}

\vspace{12pt}
The lower bound \eqref{corfinallower} is the form in which we shall complete our proof. The idea is that the random variables here should behave, for the purpose of these tail probabilities, roughly like $\log x$ {\em independent} random variables each having a $N(0,(1/2)\log\log x)$ distribution, and therefore the probability that their maximum is larger than $\frac{\log\beta}{2} + \log\log x$ should be roughly
\begin{eqnarray}
\log x \cdot \p(N(0,(1/2)\log\log x) \geq \frac{\log\beta}{2} + \log\log x) & \gg & \log x \cdot \frac{e^{-(\frac{\log\beta}{2} + \log\log x)^{2}/\log\log x}}{\sqrt{\log\log x}} \nonumber \\
& \gg & \frac{1}{\beta \sqrt{\log\log x}} , \nonumber
\end{eqnarray}
on our range of $\beta$.


To prove this, using a multivariate central limit theorem we may replace the collection of sums $( \sum_{\log^{10}x < p \leq x} \log x \int_{-\frac{1}{2\log x}}^{\frac{1}{2\log x}} \Re \frac{f(p) p^{- i\frac{k}{\log x} - it}}{p^{1/2 + \frac{\log\log x}{\log x}}} dt )_{|k| \leq \frac{\log x - 1}{2}}$ by a collection of Gaussian random variables $(X(k))_{|k| \leq \frac{\log x - 1}{2}}$ with the same means, variances and covariances. Furthermore, the mean is clearly 0 for each $k$, and for any $j$ and $k$ the covariance is
\begin{eqnarray}
&& \E \Biggl( \sum_{\log^{10}x < p \leq x} \log x \int_{-\frac{1}{2\log x}}^{\frac{1}{2\log x}} \Re \frac{f(p) p^{- i\frac{k}{\log x} - it}}{p^{1/2 + \frac{\log\log x}{\log x}}} dt \Biggr) \Biggl( \sum_{\log^{10}x < p \leq x} \log x \int_{-\frac{1}{2\log x}}^{\frac{1}{2\log x}} \Re \frac{f(p) p^{- i\frac{j}{\log x} - it}}{p^{1/2 + \frac{\log\log x}{\log x}}} dt \Biggr) \nonumber \\
& = & \sum_{\log^{10}x < p \leq x} \frac{\E (\log x \int_{-\frac{1}{2\log x}}^{\frac{1}{2\log x}} \Re f(p) p^{- i\frac{k}{\log x} - it} dt)(\log x \int_{-\frac{1}{2\log x}}^{\frac{1}{2\log x}} \Re f(p) p^{- i\frac{j}{\log x} - it} dt)}{p^{1+\frac{2\log\log x}{\log x}}} \nonumber \\
& = & \sum_{\log^{10}x < p \leq x} \frac{\E (\log x \int_{-\frac{1}{2\log x}}^{\frac{1}{2\log x}} \Re f(p) p^{- it} dt)(\log x \int_{-\frac{1}{2\log x}}^{\frac{1}{2\log x}} \Re f(p) p^{- i\frac{(j-k)}{\log x} - it} dt)}{p^{1+\frac{2\log\log x}{\log x}}}  \nonumber
\end{eqnarray}
since in the Steinhaus case $f(p)p^{-ik/\log x}$ has the same distribution as $f(p)$. Notice this means that the covariance is a function of $j-k$ (and $x$) only, in other words our random variables are {\em stationary}. Calculating further, using that $\E (\Re f(p)p^{-it_1})(\Re f(p) p^{- i\frac{(j-k)}{\log x} - it_2})$ is $(1/2)\cos((\frac{j-k}{\log x}+(t_2 - t_1))\log p) = (1/2)\cos(\frac{j-k}{\log x}\log p) + O(|t_1 - t_2|\log p)$, we find the covariance\footnote{Later we will also need to know that the covariance is a decreasing function of $|j-k|$, on a suitable range. To do this, we note that if we differentiate,
$$ \frac{d}{dh} \sum_{\log^{10}x < p \leq x} \frac{\log^{2}x \int \int \cos((h+(t_2 - t_1))\log p)}{2p^{1+\frac{2\log\log x}{\log x}}} = - \sum_{\log^{10}x < p \leq x} \frac{\log p \log^{2}x \int \int \sin((h+(t_2 - t_1))\log p)}{2p^{1+\frac{2\log\log x}{\log x}}} . $$
Then $- \sum_{\log^{10}x < p \leq x} \frac{\log p \sin(h\log p)}{p^{1+\frac{2\log\log x}{\log x}}} = \Im \sum_{p} \frac{\log p}{p^{1+\frac{2\log\log x}{\log x} + ih}} + O(\log\log x) = - \Im \frac{\zeta'}{\zeta}(1 + \frac{2\log\log x}{\log x} + ih) + O(\log\log x)$, where $\zeta(s)$ denotes the Riemann zeta function. In particular, if $|h|$ is small then it is known that $- \frac{\zeta'}{\zeta}(1 + \frac{2\log\log x}{\log x} + ih) = \frac{1}{(2\log\log x)/\log x + ih} + O(1)$. See Theorem 6.7 of Montgomery and Vaughan~\cite{mv}, for example. So provided $(\log\log x)/\log x \leq h \leq 1/(500\log\log x)$, say, the derivative will be negative. This translates into a range $\log\log x \leq |j-k| \leq (\log x)/(500\log\log x)$.} is
$$ \sum_{\log^{10}x < p \leq x} \frac{\log^{2}x \int_{-\frac{1}{2\log x}}^{\frac{1}{2\log x}} \int_{-\frac{1}{2\log x}}^{\frac{1}{2\log x}} \cos((\frac{j-k}{\log x}+(t_2 - t_1))\log p) dt_1 dt_2}{2p^{1+\frac{2\log\log x}{\log x}}} = \sum_{\log^{10}x < p \leq x} \frac{\cos(\frac{(j-k)\log p}{\log x})}{2p^{1+\frac{2\log\log x}{\log x}}} + O(1) . $$
Here the final equality also used Chebychev's bound $\sum_{p \leq x} \frac{\log p}{p} \ll \log x$.

Now using a Chebychev type bound again, the contribution to the covariance sum from primes $p > x^{1/\log\log x}$ is $\ll \frac{\log\log x}{\log x} \sum_{p \leq x} \frac{\log p}{p^{1+\frac{2\log\log x}{\log x}}} \ll \frac{\log\log x}{\log x} \sum_{n \leq x} \frac{1}{n^{1+\frac{2\log\log x}{\log x}}} \ll 1$. When $p \leq x^{1/\log\log x}$, writing $\frac{1}{p^{1+\frac{2\log\log x}{\log x}}} = \frac{1 + O((\log p \log\log x)/\log x)}{p}$ and applying the Chebychev bound $\sum_{p \leq x^{1/\log\log x}} \frac{\log p}{p} \ll \frac{\log x}{\log\log x}$, we see we can replace $\frac{1}{p^{1+\frac{2\log\log x}{\log x}}}$ by $\frac{1}{p}$ in our expression for the covariance. Thus when $|j-k| \leq \log\log x$, and in particular for the variance where $j=k$, the above is
\begin{eqnarray}
& = & \sum_{\log^{10}x < p \leq x^{1/\log\log x}} \frac{1 + O((\frac{(j-k)\log p}{\log x})^2)}{2p} + O(1) = \sum_{\log^{10}x < p \leq x^{1/\log\log x}} \frac{1}{2p} + O(1) \nonumber \\
& = & (1/2)(\log\log x - 2\log\log\log x) + O(1) , \nonumber
\end{eqnarray}
using the Mertens estimate. Similarly, when $\log\log x < |j-k| \leq (\log x)/10\log\log x$, say, then standard prime number estimates (see the proof of Lemma 5, or section 6.1 of Harper~\cite{harpergp}) imply that $ \sum_{x^{1/|j-k|} < p \leq x^{1/\log\log x}} \frac{\cos(\frac{(j-k)\log p}{\log x})}{2p} \ll 1$, so the covariance is
$$ \sum_{\log^{10}x < p \leq x^{1/|j-k|}} \frac{1 + O((\frac{(j-k)\log p}{\log x})^2)}{2p} + O(1) = (1/2)(\log\log x - \log|j-k| - \log\log\log x) + O(1) . $$

Collecting everything together, we see that for any large parameter $E \in \N$ the probability on the left hand side of \eqref{corfinallower} is
$$ \p\Biggl(\max_{|k| \leq \frac{\log x - 1}{2}} X(k) \geq \frac{\log\beta}{2} + \log\log x \Biggr) \geq \p\Biggl(\max_{|j| \leq \frac{\log x}{1000 E (\log\log x)^2}} \frac{X(jE\lfloor \log\log x \rfloor)}{\sqrt{\E X(jE \lfloor \log\log x \rfloor)^2}} \geq u \Biggr) , $$
where $u = u(\beta,x)$ is given by
$$ u := \frac{(1/2)\log\beta + \log\log x}{\sqrt{\E X(jE \lfloor \log\log x \rfloor)^2}} = \frac{(1/2)\log\beta + \log\log x}{\sqrt{(1/2)(\log\log x - 2\log\log\log x) + O(1)}} \geq \sqrt{2\log\log x} . $$
Recall here that $\E X(jE \lfloor \log\log x \rfloor)^2$ is the same for all $j$ and $E$, because our random variables are stationary, so $u(\beta,x)$ doesn't depend on $j,E$. We also always have $u(\beta,x) \ll \sqrt{\log\log x}$, on our range $2 \leq \beta \leq e^{2\sqrt{\log\log x}}$. We will choose the parameter $E$, and explain the reason for introducing it, in a moment. We can ease notation a little by setting $Z(j) := X(jE\lfloor \log\log x \rfloor)/\sqrt{\E X(jE \lfloor \log\log x \rfloor)^2}$. Then the $Z(j)$ are mean zero, variance one, stationary Gaussian random variables whose covariance function $r(m) := \E Z(j)Z(j+m)$ satisfies
$$ 0 \leq r(m) = \frac{\log\log x - \log(mE\lfloor \log\log x \rfloor) - \log\log\log x + O(1)}{\log\log x - 2\log\log\log x + O(1)} \leq 1 - \frac{2\log(mE)}{u^2} , $$
for all $1 \leq m \leq (\log x)/(500 E(\log\log x)^2)$.

At this point, having seen that $r(m)$ is non-negative and decreasing on our range of $m$ we can use Theorem 1 of Harper~\cite{harpergp},obtaining that $\p(\max_{|j| \leq \frac{\log x}{1000 E (\log\log x)^2}} Z(j) \geq u)$ is
\begin{eqnarray}
& \gg & \frac{\log x}{E (\log\log x)^2} \frac{e^{-u^{2}/2}}{u} \sqrt{\frac{1-r(1)}{u^2 r(1)}} \prod_{1 \leq m \leq \frac{\log x}{500 E(\log\log x)^2}} \Phi(u\sqrt{1-r(m)}(1 + O(\frac{1}{u^{2}(1-r(m))}))) \nonumber \\
& \gg & \frac{\sqrt{\log E}}{E} \frac{\log x}{(\log\log x)^{7/2}} e^{-u^{2}/2} \prod_{1 \leq m \leq \frac{\log x}{500 E(\log\log x)^2}} \Phi(\sqrt{u^{2}(1-r(m))}(1 + O(\frac{1}{u^{2}(1-r(m))}))) , \nonumber
\end{eqnarray} 
where $\Phi$ denotes the standard normal cumulative distribution function. Here we see, by writing $(1/2)\log\beta + \log\log x = (1/2)\log\beta + 2\log\log\log x + (\log\log x - 2\log\log\log x)$, that
$$ e^{-u^{2}/2} = \exp\{- \frac{((1/2)\log\beta + \log\log x)^2}{\log\log x - 2\log\log\log x + O(1)}\} \asymp \frac{1}{\beta \log x (\log\log x)^2} . $$
And by our previous calculations, together with the fact that $\Phi(t) = 1 - \int_{t}^{\infty} \frac{e^{-x^{2}/2}}{\sqrt{2\pi}} dx = 1 - O(\frac{e^{-t^{2}/2}}{t})$ for $t \geq 1$, the product over $m$ is
$$ \geq \prod_{m \leq \frac{\log x}{500 E(\log\log x)^2}} \Phi(\sqrt{2\log(mE)}(1 + O(\frac{1}{\log(mE)}))) = \prod_{m \leq \frac{\log x}{500 E(\log\log x)^2}} \Biggl(1 - O(\frac{1}{mE\sqrt{\log(mE)}}) \Biggr) . $$
Now this product is $= \exp\{-\sum_{m \leq \frac{\log x}{500 E(\log\log x)^2}} O(\frac{1}{mE\sqrt{\log(mE)}})\}$, so if we choose $E = \sqrt{\log\log x}$, say, then the sum is uniformly bounded and so the product is $\gg 1$. This is the step where the spacing parameter $E$ is important to obtain a good bound. Putting everything together again, we have shown that the left hand side of \eqref{corfinallower} is
$$ \gg \frac{\sqrt{\log E}}{E} \frac{\log x}{(\log\log x)^{7/2}} \cdot \frac{1}{\beta \log x (\log\log x)^2} = \frac{1}{\beta (\log\log x)^{O(1)}} , $$
as we wanted.
\qed

\vspace{12pt}
The Rademacher case is exactly the same as the Steinhaus case as far as \eqref{corprodlower}, and the subsequent application of Jensen's inequality. At that point things change a little because the Rademacher Euler product is slightly different than the Steinhaus one, and specifically in place of the term $2\sum_{p \leq x} \log x \int_{-\frac{1}{2\log x}}^{\frac{1}{2\log x}} \Re \frac{(f(p) p^{- i\frac{k}{\log x} - it})^2}{2 p^{1 + \frac{2\log\log x}{\log x}}} dt$ one has the negative of that term, and because $f(p)^2 \equiv 1$ in the Rademacher case that becomes the deterministic quantity $- \sum_{p \leq x} \log x \int_{-\frac{1}{2\log x}}^{\frac{1}{2\log x}} \Re \frac{p^{- 2i\frac{k}{\log x} - 2it}}{p^{1 + \frac{2\log\log x}{\log x}}} dt$. But this is simply $- \log x \int_{-\frac{1}{2\log x}}^{\frac{1}{2\log x}} \Re \log\zeta(1 + \frac{2\log\log x}{\log x} + i\frac{2k}{\log x} + 2it) dt + O(1)$, and if one restricts attention to the range $9(\log x)/20 \leq k \leq (\log x)/2$ then standard estimates for the zeta function (see Theorem 6.7 of Montgomery and Vaughan~\cite{mv}) show this is all $\ll 1$, so may be discarded. Moreover the restriction to the range $9(\log x)/20 \leq k \leq (\log x)/2$ makes no essential difference to the subsequent arguments, since rather than looking at points $jE\lfloor \log\log x \rfloor$ with $|j| \leq \frac{\log x}{1000 E (\log\log x)^2}$, one can simply look at points of the form $9(\log x)/20 + jE\lfloor \log\log x \rfloor$.

In the calculation of covariances that follows \eqref{corfinallower}, one runs into the usual complication in the Rademacher case that the distribution of $f(n)n^{it}$ is not the same as the distribution of $f(n)$. Here this means that the covariances are not perfectly stationary, i.e. they are a function of $j$ and $k$ rather than only a function of $j-k$. However, they are almost (i.e. up to error terms) a function of $j-k$, and in fact exactly the same function one gets in the Steinhaus case, so one can adapt the analysis accordingly. Rather than give further details, we refer the interested reader to sections 6.1-6.2 of Harper~\cite{harpergp} for an example of such an argument in an extremely similar context (basically the same as here except that only a special choice of $u$, slightly smaller than $\sqrt{2\log\log x}$, is considered).

\appendix

\section{Proofs of the Probability Results}
We prove the results on Gaussian random walks that we stated in section \ref{secprobcalc}.

\subsection{Proof of Probability Result 1}
We shall prove that for large $a$, we have
$$ \p(\sum_{m=1}^{j} G_m \leq a + 10\log j \; \forall 1 \leq j \leq n) \ll \min\{1, \frac{a}{\sqrt{n}}\} , $$
which implies the upper bound part of Probability Result 1. The proof of the lower bound is very similar, replacing $a+10\log j$ by $a-10\log j$ and proving a lower bound $\gg \min\{1, \frac{a}{\sqrt{n}}\}$ for the probability.

We may assume that $\sqrt{n} \geq a$, otherwise the result is trivial. We shall temporarily adopt the convention that $\log_{0}n = n$, and then for $k \in \N$ let $\log_{k}n := \log(\log_{k-1}n)$ denote the $k$-fold iterated logarithm. We shall need the following results:
\begin{itemize}

\item[(A1)] constant barrier case: under the hypotheses of Probability Result 1, for any $c \geq 1$ we have $\p(\sum_{m=1}^{j} G_m \leq c \; \forall 1 \leq j \leq n) \asymp \min\{1, \frac{c}{\sqrt{n}}\}$.

\item[(A2)] under the hypotheses of Probability Result 1, for any $b$ and any $c \geq 1$ we have $\p(\max_{j \leq n} \sum_{m=1}^{j} G_m \in [b,b+c]) \ll c/\sqrt{n}$.
\end{itemize}

The result (A1) is standard, but annoyingly it seems difficult to find a reference when the $G_m$ have unequal (though comparable) variances. The continuous time analogue, with $(\sum_{m=1}^{j} G_m)_{1 \leq j \leq n}$ replaced by a Brownian motion on the time interval $[0,\sum_{m=1}^{n} \E G_m^2]$, is completely standard (see e.g. section 13.4 of Grimmett and Stirzaker~\cite{gs}), and implies the lower bound in (A1) because if Brownian motion stays below $c$ on the continuous interval, it certainly does so at the discrete points corresponding to $\sum_{m=1}^{j} G_m$. For the upper bound, it suffices to handle the case $c=1$, since we can group the $G_m$ into subsums of variance $\approx c^2$ (replacing $n$ by $n/c^2$) and then multiply through by $1/c$. This case follows as in the proof of Lemma 5.1.8 of Lawler and Limic~\cite{lawlim}.

To deduce (A2), we can let $J=J(b)$ denote the smallest $1 \leq j \leq n$ at which $\sum_{m=1}^{j} G_m \geq b$, if such $j$ exists, and then $\p(\max_{j \leq n} \sum_{m=1}^{j} G_m \in [b,b+c])$ is
\begin{eqnarray}
& = & \sum_{k \leq n/2} \p(\max_{j \leq n} \sum_{m=1}^{j} G_m \in [b,b+c], \; \text{and} \; J = k) + \p(\max_{j \leq n} \sum_{m=1}^{j} G_m \in [b,b+c], \; \text{and} \; J > \frac{n}{2}) \nonumber \\
& \leq & \sum_{k \leq n/2} \p(J=k) \p(\sum_{m=k+1}^{j} G_m \leq c \; \forall k+1 \leq j \leq n) + \p(\max_{n/2 < j \leq n} \sum_{m=1}^{j} G_m \in [b,b+c]) , \nonumber
\end{eqnarray}
since the random variables $(G_m)_{m=k+1}^{n}$ are independent of the event $\{J=k\}$. By result (A1), the first sums here are $\ll \sum_{k \leq n/2} \p(J=k) c/\sqrt{n} \leq c/\sqrt{n}$. And if we let $M := \max_{n/2 < j \leq n} \sum_{n/2 < m \leq j} G_m$, then we have $\p(\max_{n/2 < j \leq n} \sum_{m=1}^{j} G_m \in [b,b+c]) = \p(\sum_{m \leq n/2} G_m \in [b-M,b-M+c])$. Since $\sum_{m \leq n/2} G_m$ is a Gaussian random variable with variance $\asymp n$, that is independent of $M$, the probability it lies in the interval $[b-M,b-M+c]$ is also $\ll c/\sqrt{n}$.

\vspace{12pt}

\begin{problem1}
In the setting of Probability Result 1, and with the above conventions, the following is true. Uniformly for any $k \in \N \cup \{0\}$ such that the $k$-fold iterated logarithm $\log_{k}n \geq a/1000$, we have
\begin{eqnarray}
&& \p(\sum_{m=1}^{j} G_m \leq a + 10\log j \; \forall \; 1 \leq j \leq n) \nonumber \\
& = & \p(\sum_{m=1}^{j} G_m \leq a + 10\log(\min\{j, \log^{20}_{k}n\}) \; \forall \; 1 \leq j \leq n) + O(\frac{1}{\sqrt{n}} \sum_{i=1}^{k} \frac{1}{\log^{2}_{i}n}) . \nonumber
\end{eqnarray}
\end{problem1}

\begin{proof}[Proof of Probability Lemma 1]
We shall prove the Lemma by induction on $k$, uniformly for each given $a$. Since we assume that $\sqrt{n} \geq a$ we certainly do have $\log_{0}n = n \geq a/1000$, and moreover the statement of the lemma is trivial for $k=0$.

For the inductive step, suppose the condition $\log_{k+1}n \geq a/1000$ is satisfied, and that we have already established the lemma at the ``level'' $k$. Let $S_{k+1} := \sum_{m \leq \log_{k+1}^{20}n} G_m$, and let $s_{k+1}^2 := \sum_{m \leq \log_{k+1}^{20}n} \E G_{m}^2 \asymp \log_{k+1}^{20}n$ denote its variance. Then conditioning on the behaviour of $S_{k+1}$, we can write $\p(\sum_{m=1}^{j} G_m \leq a + 10\log(\min\{j, \log^{20}_{k}n\}) \; \forall \; j \leq n)$ as
\begin{eqnarray}
&& \frac{1}{\sqrt{2\pi} s_{k+1}} \int_{-\infty}^{\infty} e^{-\frac{u^{2}}{2s_{k+1}^{2}}} \p(\sum_{m=1}^{j} G_m \leq a + 10\log(\min\{j, \log^{20}_{k}n\}) \; \forall \; j \leq n \; \mid \; S_{k+1} = u) du \nonumber \\
& = & \frac{1}{\sqrt{2\pi} s_{k+1}} \int_{-\infty}^{\infty} e^{-\frac{u^{2}}{2s_{k+1}^{2}}} \p(\sum_{m=1}^{j} G_m \leq a + 10\log j \; \forall \; 1 \leq j \leq \log_{k+1}^{20}n \; \mid \; S_{k+1} = u) \cdot \nonumber \\
&& \cdot \p(\sum_{\log_{k+1}^{20}n < m \leq j} G_m \leq a + 10\log(\min\{j, \log^{20}_{k}n\}) - u \; \forall \; \log_{k+1}^{20}n < j \leq n) du , \nonumber
\end{eqnarray}
where we used the independence of $\sum_{\log_{k+1}^{20}n < m \leq j} G_m$ and $\sum_{m \leq \log_{k+1}^{20}n} G_m$ to remove the conditioning on $S_{k+1}$ from the final probability here.

Now we may restrict the integration to the range $(-\infty, a + 200\log_{k+2}n]$, since if $u$ is larger then the first conditional probability is zero because $S_{k+1}$ would be too large. Meanwhile, in the part of the integral where $-\log_{k+1}^{4}n \leq u \leq a + 200\log_{k+2}n$ we can trivially upper bound the integrand by
$$ e^{-\frac{u^{2}}{2s_{k+1}^{2}}} \p(\sum_{\log_{k+1}^{20}n < m \leq j} G_m \leq a + 10\log(\min\{j, \log^{20}_{k}n\}) + \log_{k+1}^{4}n \; \forall \; \log_{k+1}^{20}n < j \leq n) , $$
which is $\ll e^{-\frac{u^{2}}{2s_{k+1}^{2}}} (a + \log_{k+1}^{4}n)/\sqrt{n} \ll e^{-\frac{u^{2}}{2s_{k+1}^{2}}} \frac{\log_{k+1}^{4}n}{\sqrt{n}}$ by result (A1) and the fact that $\log_{k+1}n \geq a/1000$. So the total contribution from this part of the integral is
$$ \ll \frac{1}{\sqrt{2\pi} s_{k+1}} \int_{-\log_{k+1}^{4}n}^{a + 200\log_{k+2}n} e^{-\frac{u^{2}}{2s_{k+1}^{2}}} \frac{\log_{k+1}^{4}n}{\sqrt{n}} du \ll \frac{(a+\log_{k+1}^{4}n)}{s_{k+1}} \frac{\log_{k+1}^{4}n}{\sqrt{n}} \ll \frac{1}{\sqrt{n} \log_{k+1}^{2}n} . $$
This may be absorbed into the ``big Oh'' term in the statement of the lemma.

Finally, when $u \leq -\log_{k+1}^{4}n$ we have
\begin{eqnarray}
&& \p(\sum_{\log_{k+1}^{20}n < m \leq j} G_m \leq a + 10\log(\min\{j, \log^{20}_{k}n\}) - u \; \forall \; \log_{k+1}^{20}n < j \leq n) \nonumber \\
& = & \p(\sum_{\log_{k+1}^{20}n < m \leq j} G_m \leq a + 10\log(\log^{20}_{k+1}n) - u \; \forall \; \log_{k+1}^{20}n < j \leq n) + O(\frac{\log_{k+1}n}{\sqrt{n}}) , \nonumber
\end{eqnarray}
by result (A2) applied with $b =  a + 10\log(\log^{20}_{k+1}n) - u$ and $c = 10(\log(\log^{20}_{k}n) - \log(\log^{20}_{k+1}n)) \ll \log_{k+1}n$. Now we may trivially rewrite the preceding line as
$$ \p(\sum_{\log_{k+1}^{20}n < m \leq j} G_m \leq a + 10\log(\min\{j, \log^{20}_{k+1}n\}) - u \; \forall \; \log_{k+1}^{20}n < j \leq n) + O(\frac{\log_{k+1}n}{\sqrt{n}}) , $$
and the crucial point is that, by result (A1) and the fact that $u \leq -\log_{k+1}^{4}n$, the probability here is $\gg (\log_{k+1}^{4}n)/\sqrt{n}$, so we can replace the absolute error term $O(\frac{\log_{k+1}n}{\sqrt{n}})$ by a multiplicative factor $1 + O(1/\log_{k+1}^{3}n)$.

Putting everything back together, we have shown that
\begin{eqnarray}
&& \p(\sum_{m=1}^{j} G_m \leq a + 10\log(\min\{j, \log^{20}_{k}n\}) \; \forall j \leq n) \nonumber \\
& = & (1 + O(\frac{1}{\log_{k+1}^{3}n}))\p(\sum_{m=1}^{j} G_m \leq a + 10\log(\min\{j, \log^{20}_{k+1}n\}) \; \forall j \leq n) + O(\frac{1}{\sqrt{n} \log_{k+1}^{2}n}) . \nonumber
\end{eqnarray}
Finally, using result (A1) we have
$$ \frac{1}{\log_{k+1}^{3}n} \p(\sum_{m=1}^{j} G_m \leq a + 10\log(\min\{j, \log^{20}_{k+1}n\}) \; \forall j \leq n) \ll \frac{a+\log_{k+2}n}{\sqrt{n} \log_{k+1}^{3}n} \ll \frac{1}{\sqrt{n} \log_{k+1}^{2}n} , $$
which gives an acceptable error term and completes our induction.
\end{proof}

\begin{proof}[Proof of Probability Result 1, upper bound]
We apply Probability Lemma 1 with $k$ chosen as large as possible, so that $\log_{k}n \geq a/1000$ but $\log_{k+1}n < a/1000$. This yields that $\p(\sum_{m=1}^{j} G_m \leq a + 10\log j \; \forall \; 1 \leq j \leq n)$ is
\begin{eqnarray}
& \leq & \p(\sum_{m=1}^{j} G_m \leq a + 200\log_{k+1}n \; \forall \; 1 \leq j \leq n) + O(\frac{1}{\sqrt{n}} \sum_{i=1}^{k} \frac{1}{\log^{2}_{i}n}) \nonumber \\
& \ll & \p(\sum_{m=1}^{j} G_m \leq 3a \; \forall \; 1 \leq j \leq n) + \frac{1}{\sqrt{n}} \ll \frac{a}{\sqrt{n}} , \nonumber
\end{eqnarray}
by (A1).
\end{proof}

\subsection{Proof of Probability Result 2}
The upper bound part of Probability Result 2 follows immediately from Probability Result 1.

To prove the lower bound part of the result, recalling that $g(j) \leq -Bj$ and $h(j) \geq -10\log j$ we see that $\p(g(j) \leq \sum_{m=1}^{j} G_m \leq \min\{a,Bj\} + h(j) \; \forall 1 \leq j \leq n)$ is
\begin{eqnarray}
& \geq & \p(-Bj \leq \sum_{m=1}^{j} G_m \leq \min\{a,Bj\} + h(j) \; \forall 1 \leq j \leq n) \nonumber \\
& \geq & \p(\sum_{m=1}^{j} G_m \leq a + h(j) \; \forall 1 \leq j \leq n) \nonumber \\
&& - \sum_{k=1}^{n} \p(\sum_{m=1}^{j} G_m \leq a + h(j) \; \forall 1 \leq j \leq n , \; \text{and} \; \sum_{m=1}^{k} G_m \leq -Bk )  \nonumber \\
&& - \sum_{k=1}^{n} \p(\sum_{m=1}^{j} G_m \leq a + h(j) \; \forall 1 \leq j \leq n , \; \text{and} \; \sum_{m=1}^{k} G_m \geq Bk - 10\log k ) . \nonumber
\end{eqnarray}
By Probability Result 1, the first term here is $\gg \min\{1,\frac{a}{\sqrt{n}}\}$, and we will show the sum of the subtracted terms is smaller than this provided $B$ is fixed large enough.

Indeed, since $(1/20)k \leq \sum_{m=1}^{k} \E G_m^2 \leq 20k$, if we condition on $\sum_{m=1}^{k} G_m$ we see $\p(\sum_{m=1}^{j} G_m \leq a + h(j) \; \forall 1 \leq j \leq n , \; \text{and} \; \sum_{m=1}^{k} G_m \leq -Bk )$ is
\begin{eqnarray}
& \ll & \frac{1}{\sqrt{k}} \int_{-\infty}^{-Bk} e^{-u^{2}/40k} \p(\sum_{m=1}^{j} G_m \leq a + h(j) \; \forall 1 \leq j \leq n \; \mid \; \sum_{m=1}^{k} G_m = u) du \nonumber \\
& \leq & \frac{1}{\sqrt{k}} \int_{-\infty}^{-Bk} e^{-u^{2}/40k} \p(\sum_{m=k+1}^{j} G_m \leq a + h(j) - u \; \forall k+1 \leq j \leq n ) du . \nonumber
\end{eqnarray}
By Probability Result 1, provided that $k \leq n/2$ (so that $n-(k+1) \gg n$) this is all
$$ \ll \frac{1}{\sqrt{k}} \int_{-\infty}^{-Bk} e^{-u^{2}/40k} \min\{1,\frac{a+|u|}{\sqrt{n}}\} du \ll \min\{1,\frac{a}{\sqrt{n}}\} \sqrt{k} e^{-B^{2}k/40} . $$
If $k > n/2$, we can upper bound the probability in the integral trivially by 1 and just use the fact that $\frac{1}{\sqrt{k}} \int_{-\infty}^{-Bk} e^{-u^{2}/40k} du \ll e^{-B^{2}k/40} \ll e^{-B^{2}n/80}$. Summing over $1 \leq k \leq n$, the contribution from all these terms will indeed be small compared with $\min\{1,\frac{a}{\sqrt{n}}\}$, provided $B$ is fixed large enough.

Similarly, we have
\begin{eqnarray}
&& \p(\sum_{m=1}^{j} G_m \leq a + h(j) \; \forall 1 \leq j \leq n , \; \text{and} \; \sum_{m=1}^{k} G_m \geq Bk - 10\log k) \nonumber \\
& \ll & \frac{1}{\sqrt{k}} \int_{Bk-10\log k}^{a+h(k)} e^{-u^{2}/40k} \p(\sum_{m=1}^{j} G_m \leq a + h(j) \; \forall 1 \leq j \leq n \; \mid \; \sum_{m=1}^{k} G_m = u) du \nonumber \\
& \leq & \frac{1}{\sqrt{k}} \int_{Bk-10\log k}^{a+h(k)} e^{-u^{2}/40k} \p(\sum_{m=k+1}^{j} G_m \leq a + h(j) - Bk + 10\log k \; \forall k+1 \leq j \leq n ) du , \nonumber
\end{eqnarray}
which is $\ll \min\{1,\frac{a}{\sqrt{n}}\} e^{-(Bk-10\log k)^{2}/40k}$ when $k \leq n/2$, and is $\ll e^{-(Bk-10\log k)^{2}/40k} \ll e^{-B^{2}n/100}$ (say) when $k > n/2$. Again summing over $1 \leq k \leq n$, the contribution from all these terms will be negligible compared with $\min\{1,\frac{a}{\sqrt{n}}\}$ if $B$ is large enough.
\qed

\vspace{12pt}
\noindent {\em Acknowledgements.} The author would like to thank Kevin Ford, for a discussion about Gaussian random walks; and Maksym Radziwi\l\l \, and Eero Saksman, for their comments on a draft of this paper. He also thanks the Heilbronn Institute for Mathematical Research for support during their May 2016 workshop on ``Extrema of Logarithmically Correlated Processes, Characteristic Polynomials, and the Riemann Zeta Function'', several sessions of which led him to start thinking about multiplicative chaos.


\begin{thebibliography}{99}

\bibitem{abh} L.-P. Arguin, D. Belius, A. J. Harper. Maxima of a randomized Riemann zeta function, and branching random walks. {\em Ann. Appl. Probab.}, \textbf{27},  no. 1, pp 178-215. 2017

\bibitem{barknsw} J. Barral, A. Kupiainen, M. Nikula, E. Saksman, C. Webb. Basic properties of critical lognormal multiplicative chaos. {\em Ann. Probab.}, \textbf{43}, no. 5, pp 2205-2249. 2015

\bibitem{berestycki} N. Berestycki. An elementary approach to Gaussian multiplicative chaos. Preprint available online at \url{http://arxiv.org/abs/1506.09113}

\bibitem{BondarenkoSeip} A. Bondarenko, K. Seip. Helson's problem for sums of a random multiplicative function. {\em Mathematika}, \textbf{62},  no. 1, pp 101-110. 2016

\bibitem{chatsound} S. Chatterjee, K. Soundararajan. Random multiplicative functions in short intervals. {\em Int. Math. Res. Not.}, pp 479-492. 2012

\bibitem{duprsvnorm} B. Duplantier, R. Rhodes, S. Sheffield, V. Vargas. Renormalization of critical Gaussian multiplicative chaos and KPZ relation. {\em Comm. Math. Phys.}, \textbf{330},  no. 1, pp 283-330. 2014

\bibitem{gs} G. R. Grimmett, D. R. Stirzaker. {\em Probability and Random Processes.} Third edition, published by Oxford University Press. 2001

\bibitem{gut} A. Gut. {\em Probability: A Graduate Course.} Second edition, published by Springer Texts in Statistics. 2013

\bibitem{halasz} G. Hal\'{a}sz. On random multiplicative functions. In {\em Hubert Delange Colloquium, (Orsay, 1982)}. {\em Publications Math\'{e}matiques d'Orsay}, \textbf{83}, pp 74-96. Univ. Paris XI, Orsay. 1983

\bibitem{harpergp} A. J. Harper. Bounds on the suprema of Gaussian processes, and omega results for the sum of a random multiplicative function. {\em Ann. Appl. Probab.}, \textbf{23}, no. 2, pp 584-616. 2013

\bibitem{harperlimits} A. J. Harper. On the limit distributions of some sums of a random multiplicative function. {\em Journal f\"{u}r die reine und angewandte Mathematik}, \textbf{678}, pp 95-124. 2013

\bibitem{harperrmfhigh} A. J. Harper. Moments of random multiplicative functions, II: High moments. {\em In preparation.}

\bibitem{hnr} A. J. Harper, A. Nikeghbali, M. Radziwi\l\l. A note on Helson's conjecture on moments of random multiplicative functions. In {\em Analytic number theory}, pp 145-169, Springer, Cham. 2015.

\bibitem{heaplindqvist} W. Heap, S. Lindqvist. Moments of random multiplicative functions and truncated characteristic polynomials.  Preprint available online at \url{http://arxiv.org/abs/1505.03378}

\bibitem{helson} H. Helson. Hankel Forms. {\em Studia Math.}, \textbf{198}, no. 1, pp. 79-84, 2010

\bibitem{hough} B. Hough. Summation of a random multiplicative function on numbers having few prime factors. {\em Math. Proc. Cambridge Philos. Soc.}, \textbf{150}, pp 193-214. 2011

\bibitem{hushi} Y. Hu, Z. Shi. Minimal position and critical martingale convergence in branching random walks, and directed polymers on disordered trees. {\em Ann. Probab.}, \textbf{37}, no. 2, pp 742-789. 2009

\bibitem{tenenbaum} Y-K. Lau, G. Tenenbaum, J. Wu. On mean values of random multiplicative functions. {\em Proc. Amer. Math. Soc.}, \textbf{141}, pp 409-420. 2013. Also see \url{www.iecl.univ-lorraine.fr/~Gerald.Tenenbaum/PUBLIC/Prepublications_et_publications/RMF.pdf} for some corrections to the published version.

\bibitem{lawlim} G. F. Lawler, V. Limic. {\em Random walk: a modern introduction.} First edition, published by Cambridge University Press. 2010

\bibitem{mv} H. L. Montgomery, R. C. Vaughan. {\em Multiplicative Number Theory I: Classical Theory.} First edition, published by Cambridge University Press. 2007

\bibitem{ng} N. Ng. The distribution of the summatory function of the M\"{o}bius function. {\em Proc. London. Math. Soc.}, \textbf{89}, no. 3, pp. 361-389, 2004.

\bibitem{rhodesvargas} R. Rhodes, V. Vargas. Gaussian multiplicative chaos and applications: A review. {\em Probab. Surveys}, \textbf{11}, pp 315-392. 2014

\bibitem{sadikova} S. M. Sadikova. Two-dimensional analogues of an inequality of Esseen with applications to the Central Limit Theorem. {\em Theory of Probability and Its Applications}, \textbf{11}, pp 325-335. 1966

\bibitem{saksmanseipintmeans} E. Saksman, K. Seip. Integral means and boundary limits of Dirichlet series. {\em Bull. Lond. Math. Soc.}, \textbf{41}, no. 3, pp 411-422. 2009

\bibitem{saksmanseip} E. Saksman, K. Seip. Some open questions in analysis for Dirichlet series. In {\em Recent progress on operator theory and approximation in spaces of analytic functions}, pp 179-191, Contemp. Math., 679, Amer. Math. Soc., Providence, RI. 2016

\bibitem{saksmanwebb} E. Saksman, C. Webb. Multiplicative chaos measures for a random model of the Riemann zeta function. Preprint available online at \url{http://arxiv.org/abs/1604.08378}

\bibitem{saksmanwebb2} E. Saksman, C. Webb. The Riemann zeta function and Gaussian multiplicative chaos: statistics on the critical line. Preprint available online at \url{http://arxiv.org/abs/1609.00027}

\bibitem{soundmoments} K. Soundararajan. Moments of the Riemann zeta function. {\em Ann. Math.}, \textbf{170}, pp 981-993. 2009

\bibitem{weber} M. J. G. Weber. $L^{1}$-Norm of Steinhaus chaos on the polydisc. {\em Monatsh. Math.}, \textbf{181}, no. 2, pp 473-483. 2016

\bibitem{wintner} A. Wintner. Random factorizations and Riemann's hypothesis. {\em Duke Math. J.}, \textbf{11}, pp 267-275. 1944


\end{thebibliography}
\end{document}